\documentclass[11pt]{article}
\usepackage{graphicx}
\usepackage{amssymb}
\def\reals{\mathbb R}
\def\ints{\mathbb Z}

\def\complexes{\mathbb C}
\def\semi{\hbox{ $\times $ \kern-.972em \raise.12719em\hbox{ $_{^|}$}  }}
\newtheorem{theorem}{Theorem}
\newtheorem{lemma}{Lemma}[section]

\newtheorem{corollary}{Corollary}[section]
\newtheorem{remark}{Remark}[section]

\newtheorem{openproblem}{Open Problem}
\newtheorem{example}{Example}[section]

\def\ms{\bigskip}
\def\ms{\medskip}
\def\semi{\hbox{ $\times $ \kern-.972em \raise.12719em\hbox{ $_{^|}$}  }}
\def\mapright#1{\smash{\mathop{\longrightarrow}\limits^{#1}}}
\def\pf {\noindent {\bf Proof:} \ }
\def\endpf{$\|$ \bigskip}
\def\be{\begin{enumerate}}
\def\ee{\end{enumerate}}
\def\bi{\begin{itemize}}
\def\ei{\end{itemize}}

\def\cA{{\cal A}}
\def\cB{{\cal B}}
\def\cC{{\cal C}}
\def\cD{{\cal D}}
\def\cE{{\cal E}}

\def\cK{{\cal K}}
\def\cL{{\cal L}}
\def\cM{{\cal M}}

\def\cO{{\cal O}}

\def\cR{{\cal R}}

\def\cT{{\cal T}}

\def\cY{{\cal Y}}
\def\sgn{\Sigma_n}
\def\axis{\bf A}
\def\semi{\hbox{ $\times $ \kern-.972em \raise.12719em\hbox{ $_{^|}$}}}

\begin{document}

\markboth{Joan S. Birman and Tara E. Brendle }{BRAIDS: A SURVEY}
\title{BRAIDS: A SURVEY}
\author{Joan S. Birman
\thanks {The first author acknowledges partial support from the U.S.National Science Foundation under  grant number 0405586.}
\\e-mail  jb@math.columbia.edu  \and Tara E. Brendle
\thanks{The second author is partially supported by a VIGRE postdoc under NSF
grant number 9983660 to Cornell University.}
\\e-mail brendle@math.cornell.edu}
\date{December 2, 2004}
\maketitle

\begin{abstract}
This article is about Artin's braid group ${\bf B}_n$ and its role in
knot theory.  We set ourselves two goals: (i) to provide enough of
the essential background so that our review would be accessible to
graduate students, and (ii) to focus on those parts of the subject
in which major progress was made, or interesting new proofs of
known results were discovered, during the past 20 years. A central
theme that we try to develop is to show ways in which structure
first discovered in the braid groups generalizes to structure in
Garside groups, Artin groups and surface mapping class groups.
However, the  literature is extensive, and for reasons of space
our coverage necessarily omits many very interesting developments.
Open problems are noted and so-labelled, as we encounter them.
A guide to computer software is given together with a 10 page bibliography.
 \end{abstract}

 \

\tableofcontents
\newpage
\section{Introduction}
\label{section:introduction}

\begin{figure}[htpb!]
\centerline{\includegraphics[scale=.9, bb=82 425 458 513] {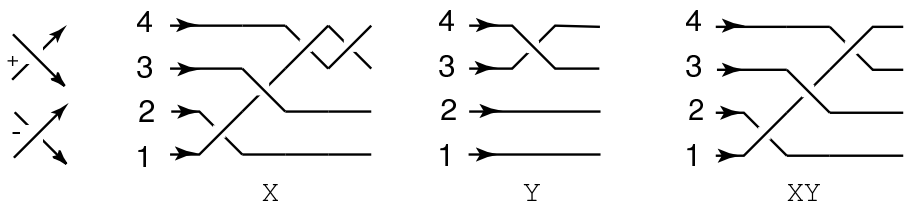}}
% \centerline{\includegraphics[scale=.9] {braid.eps}}
\caption{ Examples of 4-braids $X, Y$ and their product $XY$ }
\label{figure:braids}
\end{figure}

In a review article, one is obliged
to begin with definitions. Braids can be defined by very simple
pictures such as the ones in Figure \ref{figure:braids}. Our
braids are illustrated as oriented from left to right, with the
strands numbered $1,2,\dots,n$ from bottom to top.  Whenever it is more convenient,
we will also think of braids `vertically', i.e., oriented from top
to bottom, with the strands numbered $1,2,\dots,n$ from left to right.
Crossings are suggested as they are in a picture of a highway
overpass on a map. The identity braid has a canonical
representation in which two strands never cross. Multiplication of
braids is by juxtaposition, concatenation, isotopy and rescaling.

Pictures like the ones in Figure \ref{figure:braids} give an
excellent intuitive feeling for the braid group, but one that
quickly becomes  complicated when one tries to pin down details.
What is the ambient space? (It is a slice $\reals^2\times I$ of
3-space.)  Are admissible isotopies constrained to $\reals^2\times
I$? (Yes, we cannot allow isotopies in which the strands are
allowed to loop over the initial points.) Does isotopy mean
level-preserving isotopy? (No, it will not matter if we allow more
general isotopies in $\reals^2\times I$, as long as  strands don't
pass through one-another.)  Are  strands allowed to
self-intersect? (No, to allow self-intersections would give the
`homotopy braid group', a proper homomorphic image of the group
that is our primary focus.) Can we replace the ambient space by
the product of a more general surface and an interval, for example
a sphere and an interval? (Yes, it will become obvious shortly how
to modify the definition.)

We will bypass these questions and other related ones by giving several more
sophisticated definitions. In $\S$\ref{subsection:bBn and bPn via
configuration spaces},  $\S$\ref {subsection:definition of bBn via
generators and relations} and $\S$\ref {subsection:bBn and bPn as
mapping class groups} we will define the braid group ${\bf B}_n$ and
pure braid group ${\bf P}_n$ in three distinct ways.  We will give a
proof that two of them yield the same group.  References to the
literature establish the isomorphism in the remaining case. In
$\S$\ref{subsection:braiding is fundamental} we will demonstrate
the universality of `braiding' by describing four examples which
show how braids have played a role in parts of mathematics which
seem far away from knot theory.

\subsection{${\bf B}_n$ and ${\bf P}_n$ via configuration spaces}
\label{subsection:bBn and bPn via configuration spaces} We define
the topological concept of a braid and of a group of braids
via the notion of a configuration space.  This approach is nice
because it gives, in a concise way, the appropriate equivalence
relations and the group law, without any fuss.

The \underline{configuration space} of $n$ points on the complex plane $\complexes$ is:
$$ \cC_{0,{\hat n}} = \cC_{0,{\hat n}}(\complexes) = \{(z_1,\dots,z_n) \in \complexes \times \dots \times \complexes  \  |  \  z_i
\neq z_j \ {\rm if} \  i \neq j \}.$$
A point on $\cC_{0,{\hat n}}$ is denoted   by a vector ${\vec z} = (z_1,\dots,z_n)$.   The
symmetric group
acts freely on
$\cC_{0,{\hat n}}$, permuting the coordinates in each $\vec{z} \in \cC_{0,{\hat n}}$. The orbit space of the action is
$\cC_{0,n} = \cC_{0,{\hat n}}/\Sigma_n$  and the orbit space projection is $\tau:\cC_{0,{\hat n}} \to
\cC_{0,n}$.  Choosing a fixed base point $\vec{p} = (p_1,\dots,p_n)$, we define the
\underline{pure braid group} ${\bf P}_n$ on $n$ strands and the
\underline{braid group} ${\bf B}_n$ on $n$ strands to be the fundamental groups:
$${\bf P}_n = \pi_1({\cC_{0,{\hat n}}, \vec p}), \ \ \ \ \ \ {\bf B}_n =
\pi_1(\cC_{0,n},\tau({\vec p})).$$

At first encounter $\pi_1(\cC_{0,n},\tau({\vec p}))$ doesn't look
as if it has much to do with braids as we illustrated them in
Figure \ref{figure:braids}, but in fact there is a simple
interpretation  which reveals the intuitive picture.  While the
manifold $\cC_{0,n}$ has dimension $2n$, the fact that the points
$z_1,\dots,z_n$ are pairwise distinct allows us to think of a
point ${\vec z} \in \cC_{0,n}$ as a set of $n$ distinct points on
$\complexes$. An element of $\pi_1(\cC_{0,n},\tau({\vec p}))$ is
then represented by a loop which lifts uniquely to a path ${\vec
g}:I \to \cC_{0,{\hat n}}$, where $\vec{g} = \langle g_1,\dots,g_n \rangle$
consists of $n$ coordinate functions $g_i:I \to \complexes$ satisfying
$g_i(t) \neq g_j(t)$ if $i \neq j, \ t \in [0,1],$ also
$\vec{g}(0) = \vec{g}(1) = \vec{p}$, the base point.   The graph
of the $n$ simultaneous functions $g_1,\dots,g_n$ is a
\underline{geometric (pure) braid}. The appropriate equivalence
relation on geometric braids is captured by simultaneous homotopy
of the $n$ simultaneous paths, rel their endpoints, in the
configuration space.

The group ${\bf B}_n$ is the group which Artin set out to investigate
in 1925 in his seminal paper \cite{Artin1925} (however he defined
it in a less concise way); in the course of his investigations he
was led almost immediately to study ${\bf P}_n$. Indeed, the two braid
groups are related in a very simple way.  Let $\tau_*$ be the
homomorphism on fundamental groups which is induced by the orbit
space projection  $\tau: \cC_{0,{\hat n}} \to \cC_{0,n}$. Observe
that the orbit space projection  is a regular $n!$-sheeted
covering space
 projection, with $\sgn$
as the group of covering translations.  From this is follows that
the group ${\bf P}_n$ is a subgroup of index $n!$ in ${\bf B}_n$, and there is a short exact sequence:

\begin{equation}
 \label{equation:bbn, bpn and symmetric gp}
1 \to {\bf P}_n \mapright{\tau_*} {\bf B}_n
\to
\Sigma_n \longrightarrow 1
\end{equation}

\subsection{${\bf B}_n$ and ${\bf P}_n$ via generators and relations}
\label{subsection:definition of bBn via generators and relations}
\begin{figure}[htpb!]
\centerline{\includegraphics[scale=.7, bb=63 382 444 523] {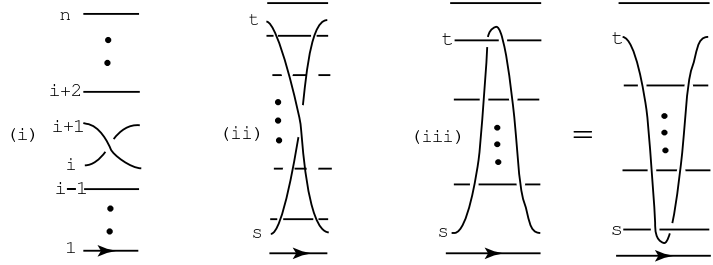}}
%\centerline{\includegraphics[scale=.7] {braids-elem.eps}}
\caption{(i) The elementary braid $\sigma_i$. (ii) The elementary
braid $\sigma_{s,t}$. (iii) The pure braid $A_{s,t} =
A_{s,t} = \sigma_{s,t}^2$.  } \label{figure:elem-bds}
\end{figure}

We give two definitions of the group ${\bf B}_n$  by generators and relations.  The \underline {classical presentation} for ${\bf B}_n$  first appeared in \cite{Artin1925}.   We record it now, and will refer back to it many times later. It has generators $\sigma_1,\dots,\sigma_{n-1}$ and defining relations:
\begin{equation}
\label{equation:classical presentation}
\sigma_i \sigma_k = \sigma_k \sigma_i \  \ {\rm if} \  \ |i-k| \geq 2, \  \  \sigma_i \sigma_{i+1} \sigma_i =  \sigma_{i+1} \sigma_i \sigma_{i+1}.
\end{equation}
The elementary braid $\sigma_i$ is depicted in sketch (i) of
Figure \ref{figure:elem-bds}.

Many years after Artin did his fundamental work,  Birman, Ko and
Lee discovered a new presentation, which enlarged the set of
generators to a more symmetric set.  Let  $\sigma_{s,t} =
(\sigma_{t-1}\cdots\sigma_{s+1})\sigma_s(\sigma_{s+1}^{-1}\cdots\sigma_{t-1}^{-1})$,
where $1\leq s < t \leq n$. We define $\sigma_{s,t} =
\sigma_{t,s}$, and adopt the convention that whenever it is
convenient to do so we will write the smaller subscript first. The
\underline{new presentation}  has generators $\{\sigma_{s,t}, \ 1
\leq s < t \leq n \}$ and defining relations:
\begin{eqnarray}
\label{equation:new presentation}
  \sigma_{s,t}\sigma_{q,r} &=&  \sigma_{q,r}\sigma_{s,t}   \  \  \  \  \    \  \  \  \  \   \  \  \  \  \  {\rm if}  \  \   (t-r)(t-q)(s-r)(s-q)>0,  \\
    \sigma_{s,t} \sigma_{r,s} &=&  \sigma_{r,t}\sigma_{s,t} =  \sigma_{r,s} \sigma_{r,t}   \  \  { \rm   if} \  \  \  \
 1 \leq  r <  s < t \leq n. \nonumber
\end{eqnarray}
See sketch (ii) of Figure \ref{figure:elem-bds} for a picture of
$\sigma_{s,t}$, and  \cite{BKL1998} for a proof that
(\ref{equation:classical presentation}) and (\ref{equation:new
presentation}) define the same group.  Both presentations will be
needed in our work.  Note that the new generators include the old ones as a proper subset, since $\sigma_i = \sigma_{i,i+1}$ for each $i = 1,2,\dots, n-1.$

By (\ref{equation:bbn, bpn and symmetric gp}) the pure braid group ${\bf P}_n$ has index $n!$ in ${\bf B}_n$.  Let $A_{s,t} = A_{t,s} = \sigma_{s,t}^2$. (See sketch (iii) of Figure
\ref{figure:elem-bds}). The symmetry $A_{s,t} = A_{t,s}$ can be seen by tightening the $s^{th}$ strand at the expense of loosening the $t^{th}$ strand.   It is proved in \cite{Artin1925}  and also in \cite{F-vB1962} that ${\bf P}_n$ has a presentation with generators $A_{r,s}, \ \ 1\leq r < s \leq n $ and defining relations:

\begin{eqnarray}
\label{equation:presentation bPn}
A_{r,s}^{-1}A_{i,j}A_{r,s} & = & A_{i,j}  \ \
{\rm if} \ \  1\leq r<s<i<j\leq n \ \ {\rm or} \ \  \ \  1\leq i<r<s<j\leq n  \cr
& = & A_{r,j}A_{i,j}A_{r,j}^{-1}    \ \   {\rm if}  \  \  1\leq r<s=i<j\leq n, \cr
& = & (A_{i,j}A_{s,j})A_{i,j} (A_{i,j}A_{s,j})^{-1} \ \   {\rm if} \ \  1\leq r=i<s<j\leq n,\cr
&  =  &  (A_{r,j}, A_{s,j}  A_{r,j}^{-1} A_{s,j}^{-1} ) A_{i,j} (A_{r,j}, A_{s,j}  A_{r,j}^{-1} A_{s,j}^{-1} )^{-1} \ \   {\rm  if} \ \  1\leq r<i<s<j\leq n.
 \end{eqnarray}
The relations in (\ref{equation:presentation bPn}) come from the existence of a split short exact sequence, for every $k=2,\dots,n$:
 \begin{equation}
 \label{equation:Pn, Pn-1 and Fn-1}
 \{1\} \to  {\bf F}_{n-1} \to {\bf P}_n \stackrel{\pi_n^\star}\to {\bf P}_{n-1} \to \{1\}.
 \end{equation}
The map $\pi_n^\star$ is defined by pulling out the last braid strand, and the image of ${\bf P}_{n-1}$ under its inverse embeds ${\bf P}_{n-1}$ in  ${\bf P}_n$, as the subgroup generated by pure braids on the first $n-1$ strands. The free subgroup $F_{n-1}$ is generated by the braids $A_{1,n}, A_{2,n},\dots, A_{n-1,n}$.  The pure braid group ${\bf P}_2$ is infinite cyclic and generated by $A_{1,2}$.  Inducting on $n$, the structure of ${\bf P}_n$ via a sequence of semi-direct products is uncovered.

\subsection {${\bf B}_n$ and ${\bf P}_n$ as mapping class groups}
\label{subsection:bBn and bPn as mapping class groups}

Our announced goal in this review was to concentrate on areas
where there have been new developments in recent years.  While it
has been known for a very long time that Artin's braid group is
naturally isomorphic to the mapping class group of an n-times
punctured disc, people have asked us many times for a simple proof
of this fact.  We do not know of a simple one in the literature,
therefore we present one here. In this case `simple' does not mean
intuitive and based upon first principles, rather it means  using
machinery which is normally available to a graduate student who
has the tools learned in a first year graduate course in topology,
and is preparing to begin research.

Let $S = S_{g,b,n}$ denote a 2-manifold of genus $g$ with $b$
boundary components and $n$ punctures, and let
$\textrm{Diff}^+(S)$ denote the groups of all  orientation
preserving diffeomorphisms of $S$. Observe that we may assign the
compact open topology to $\textrm{Diff}^+(S)$, making it into a
topological group.  The \underline{mapping class group} $\cM =
\cM_{g,b,n}$ of $S$ is $\pi_0(\textrm{Diff}^+(S))$, that is, the
quotient of $\textrm{Diff}^+(S)$ modulo its subgroup of all
diffeomeorphisms of $S$ which are isotopic to the identity rel
$\partial S$. We allow diffeomeorphisms in $\textrm{Diff}^+(S)$ to
permute the punctures, writing $\textrm{Diff}^+(S_{g,b,{\hat n}})$
if they are to be fixed pointwise.  Our interest in this article is mainly
in the special case of $\cM_{0,1,n}$.

\begin{theorem}
\label{theorem:isomorphism between bbn and MCG of punctured disc}
There are natural isomorphisms:
$$ {\bf B}_n \cong \cM_{0,1,n}   \ \ \  {\rm and} \ \ \
{\bf P}_n \cong \cM_{0,1,{\hat n}}  $$
\end{theorem}

\pf We begin with an intuitive description of  how to pass from
diffeomorphisms to geometric braids and back again. Choose any $h
\in \textrm{Diff}^+( S_{0,1,n})$. While $h$ is in general  not
isotopic to the identity, its image $i(h)$ in $\textrm{Diff}^+(
S_{0,1,0})$ under the inclusion map $i:\textrm{Diff}^+( S_{0,1,n})
\to \textrm{Diff}^+( S_{0,1,0})$ is, because $\textrm{Diff}^+(
S_{0,1,0}) = \{ 1 \}$.  Let $h_t$ denote the isotopy.  If the
punctures in $S_{0,1,n}$ are at $(p_1,p_2,\dots,p_n),$ then the
$n$ paths $(h_t(p_1),h_t(p_2),\dots,h_t(p_n))$ defined by the
traces of the points $(p_1,p_2,\dots,p_n)$ under the isotopy sweep
out  a braid in $S_{0,1,0}\times [0,1]$, and the equivalence class of this braid is the image
of the mapping class $[h]$ in the braid group ${\bf B}_n$.

It's a little bit harder to understand the inverse isomorphism,
from the braid group to the mapping class group. One chooses a
geometric braid and imagines it as being located in a slice of
3-space, with the bottom endpoints of the $n$ braid strands (which
are oriented top to bottom) as being pinned to the distinguished
points $p_1,\dots,p_n$ on the punctured disc. If one is very
careful the $n$ braid strings can be laid down on the punctured
disc so that they become $n$ non-intersecting simple arcs, each of
which begins and ends at a base point. One then constructs a
homeomorphism of the punctured disc to itself  in such a way that
the trace of the isotopy to the identity is the given set of $n$
non-intersecting simple arcs.

To prove the theorem, we begin by establishing the isomorphism
between  ${\bf P}_n$ and $\cM_{0,1,{\hat n}}$.  A good general
reference for the underlying mathematics is Chapter 6 of the
textbook \cite{Bredon}.   As previously noted $\textrm{Diff}^+(
S_{0,1,0})$ is a topological group.  Also, $\textrm{Diff}^+(
S_{0,1,{\hat n}})$ is a closed subgroup of $\textrm{Diff}^+(
S_{0,1,0})$.   The \underline{evaluation map} \ \  $\cE:
\textrm{Diff}^+( S_{0,1,0}) \to \cC_{0,{\hat n}}$ is defined  by
$\cE (h) = (h(p_1),\dots,h(p_n))$.  It is clear that $\cE$ is
continuous with respect to the compact open topology on
$\textrm{Diff}^+( S_{0,1,{\hat n}})$ and the subspace topology for
$\cC_{0,{\hat n}} \subset
\complexes\times\complexes\times\cdots\times\complexes$. The
topological group $\textrm{Diff}^+( S_{0,1,0})$ acts
$n$-transitively on the disc in the sense that  if
$(p_1,\dots,p_n)$ are $n$ distinct points and $(w_1,\dots,w_n)$
are $n$ others then there is an $h\in \textrm{Diff}^+(S_{0,1,0})$
such that $h(p_i) = w_i, \ i=1,\dots,n$. Observe that if $h\in
\textrm{Diff}^+(S_{0,1,{\hat n}})$ then $(h(p_1),\dots,h(p_n)) =
(p_1,\dots,p_n)$ and if $h,h'\in \textrm{Diff}^+(S_{0,1,0})$ with
$\cE(h) = \cE(h')$, then $h,h'$ are in the same left coset of
$\textrm{Diff}^+( S_{0,1,{\hat n}})$ in $\textrm{Diff}^+(
S_{0,1,0})$.    In this situation it is shown in \cite{Steenrod}, part I, Sections
7.3 and 7.4, that the 3-tuple $({\cal E}, \ {\rm Diff}^+( S_{0,1,0}), \
 \cC_{0,{\hat n}})$ is a fiber space,
with total space Diff$^+( S_{0,1,0})$, base space $\cC_{0,{\hat n}},$  projection $\cE$ and fiber Diff$^+(S_{0,1,\hat{n}})$.  (It is a good exercise for a graduate student to prove this directly by constructing the required local product structure in an explicit manner.)  The long exact sequence of homotopy groups of a fibration then gives the following exact sequence of groups
and homomorphisms, where we focus on the range that is of
interest:
$$\dots \to \pi_1( \textrm{Diff}^+( S_{0,1,0}))
\stackrel{\cE_*}\to\pi_1(\cC_{0,{\hat n}})
\stackrel{\partial_*}\to \pi_0(\textrm{Diff}^+( S_{0,1,{\hat n}}))
\stackrel{i_*}\to \pi_0( \textrm{Diff}^+( S_{0,1,0}))
\stackrel{\cE_*} \to \dots $$
The two end groups are trivial. The
left middle group is ${\bf P}_n$ and the right middle group is
$\cM_{0,1,{\hat n}}$. The isomorphism of Theorem
\ref{theorem:isomorphism between bbn and MCG of punctured disc} is
$\partial_*$.  Tracing through the mathematics one finds that in
fact its inverse is the map that we described right after we
stated the theorem.  The assertion about ${\bf P}_n$ is therefore true.
The proof for ${\bf B}_n$ can then be completed by  comparing the short
exact sequence (\ref{equation:bbn, bpn and
symmetric gp}), which says that ${\bf B}_n$ is a finite extension of
${\bf P}_n$ with quotient the symmetric group with a related short
exact sequence for the mapping class groups:
\begin{equation}
\label{equation:MCG gp,pure MCG and symmetric gp}
1 \to \cM_{0,1,{\hat n}} \mapright{j_*} \cM_{0,1,n} \to \Sigma_n \longrightarrow 1
\end{equation}

\noindent Thus $\cM_{0,1,n}$ is a finite extension of
$\cM_{0,1,{\hat n}}$ with quotient the symmetric group $\Sigma_n$.
Comparing corresponding groups in the short exact sequences
( \ref{equation:bbn, bpn and symmetric gp}) and
(\ref{equation:MCG gp,pure MCG and symmetric gp}), we see that the
first two and last two are isomorphic. The  5-Lemma then shows
that the middle ones are too. This completes the proof of Theorem
\ref{theorem:isomorphism between bbn and MCG of punctured
disc}.\endpf

\subsection {Some examples where braiding appears in mathematics, unexpectedly}
\label{subsection:braiding is fundamental}
We discuss, briefly, a variety of examples, outside of knot theory, where `braiding' is an essential
aspect of a mathematical or physical problem.

\subsubsection{Algebraic geometry}
\label{subsubsection:braids in algebraic geometry}
Configuration spaces and the braid group appear in a natural way in algebraic geometry.  Consider the complex
polynomial
$$(X-z_1)(X-z_2)\dots(X-z_n) = X^n+a_1X^{n-1}+\dots+ a_{n-1}X + a_n$$
of degree $n$ with $n$ distinct complex roots $z_1,\dots,z_n$.
The coefficients $a_1,\dots,a_n$ are the elementary symmetric polynomials in
 $\{z_1,\dots,z_n\}$, and so we get a continuous map $\complexes^n \to \complexes^n$ which takes roots to
coefficients.    Two points have the same image if and only if they differ by a permutation,
 so we get the same identification as in the quotient map $\tau: \cC_{0,{\hat n}} \to \cC_{0,n}$, in quite a
different way.  Since we are requiring that our polynomial have $n$ distinct roots, a  point
 $\{a_1,\dots,a_n\}$ is in the image of $\vec{z}$  under the root-to-coefficient map if and only
if the polynomial
$X^n+a_1X^{n-1}+\dots+a_n$ has $n$ distinct roots, i.e. if and only if its coefficients avoid the points where
the discriminant $$\Delta=\prod_{i<j}(z_i-z_j)^2,$$ expressed as a polynomial
 in $\{a_1,\dots,a_n\}$, vanishes.   Thus
$\cC_{0,n}(\complexes)$ can be interpreted as the complement in
$\complexes^n$ of the algebraic hypersurface
 defined by the equation $\Delta = 0$, where $\Delta$ is rewritten as a
polynomial in the coefficients $a_1,\dots,a_n$. (For example, the polynomial $X^2 +a_1X + a_2$ has distinct roots
precisely when $a_1^2 - 4a_2 = 0$).   In
this setting the base point $\tau({\vec p})$  is regarded as the choice of a complex polynomial of degree $n$
which has $n$ distinct roots, and an element in the braid group is a choice of a continuous deformation of
that polynomial along a path on which two roots never coincide.  There is a substantial literature in this area,
from which we mention only one paper, by Gorin and Lin \cite{Gorin-Lin}.  We chose it because it contains a description of the commutator subgroup ${\bf B}_n'$ of the braid group, and many people have asked the first author for a reference on that over the years.  While there may be other references they are unknown to us.

\subsubsection{Operator algebras}
\label{subsubsection:braiding in operator algebras}

Our next example, taken from the work of Vaughan Jones \cite{Jones1983},\cite{Jones1986}, is interesting because it shows how `braiding' can appear in disguise, so that initially one misses the connection.
  We consider the theory, in operator algebras, of `type II$_1$ factors', ordered by inclusion.  Let
$M$ denote a Von Neumann algebra, i.e. an algebra of bounded operators acting on a Hilbert space $h$. The algebra $M$ is
called a \underline{factor} if its center consists only of scalar multiples of the identity. The factor is  type II$_1$
if it admits a linear functional, called a trace,  $tr:M \to \complexes$, which satisfies the following three
conditions: (i) $tr(xy) = tr(yx) \  \forall  \ x,y \in M$, \ \   (ii)       $tr(1)=1 $, \ \   and   $tr(xx^\star) > 0$, where
$x^\star$ is the adjoint of $x$. In this situation it is known that the trace is unique, in the sense that it is the only
linear functional satisfying the  first two conditions. An old discovery of Murray and Von Neumann was that factors of type
II$_1$ provide a type of `scale' by which one can measure the dimension  of $h$. The notion of dimension which occurs
here generalizes the familiar notion of integer-valued dimensions, because for appropriate $M$ and $h$ it can be any
non-negative real number or $\infty$.

 The starting point of Jones' work was the following question: if $M_1$ is a type II$_1$ factor and if  $M_0 \subset
M_1$ is a subfactor, is there any restriction on the real numbers which occur as the ratio
$\lambda = {\rm dim}_{M_0}(h)/{\rm dim}_{M_1}(h)$ ?
The question has the flavor of questions one studies in Galois theory. On the face of it, there was no
reason to  think that $\lambda$ could not take on any value in $[1,\infty]$, so Jones' answer came as a complete
surprise. He called $\lambda$ the \underline{index} $|M_1:M_0|$ of $M_0$ in $M_1$, and proved a type of rigidity theorem
about type
II$_1$ factors and their subfactors:\ms

\noindent {\bf The Jones Index Theorem:}  $\lambda \subset [4,\infty] \cup \{4 cos
2\pi/p\}$ , where $ p \in \ints, p\geq 3$. Moreover, each real number in the continuous part of the spectrum
$[4,\infty]$ and in the discrete part $\{\{4 cos 2\pi/p\}, p\in \ints, p\geq 3\}$ is realized. \ms

What does all this have to do with braids?  To answer the question, we sketch the idea of the proof, which is to be found in \cite{Jones1983}.   Jones begins with the type II$_1$ factor
$M_1$ and the subfactor $M_0$. There is also a tiny bit of additional structure: It turns out that in this setting there
exists a map
$e_1:M_1\to M_0$, known as the \underline{conditional expectation} of $M_1$ on $M_0$. The map $e_1$ is a projection, i.e.
$e_1^2= e_1$.  His first step is to prove that the ratio $\lambda$ is independent of the choice of the Hilbert space
$h$. This allows him to choose an appropriate $h$ so that the algebra $M_2$ generated by $M_1$ and $e_1$
makes sense. He then investigates $M_2$ and proves that it is another type II$_1$ factor, which contains $M_1$ as a
subfactor, moreover $|M_2:M_1| = |M_1:M_0| = \lambda$. Having in hand another
II$_1$ factor
$M_2$ and its subfactor $M_1$, there is also a trace on
$M_2$ (which by the uniqueness of the trace) coincides with the trace on $M_1$ when it is restricted to $M_1$, and another
conditional expectation $e_2: M_2 \to M_1$. This allows Jones to iterate the construction, to build algebras
$M_1,M_2, \dots $ and from them a family of algebras $\{J_n, \ n=1,2,3,\dots\}$, where $J_n$ is generated by
$1,e_1,\dots,e_{n-1}$.

    Rewriting history a little bit in order to make the subsequent connection with braids a little more transparent,
we now replace the projections $e_k$, which are not units, by a new set of generators which are units, defining:
    $g_k = t e_k - (1-e_k)$, where      $(1-t)(1-t^{-1}) = 1/\lambda$.
The $g_k$'s generate $J_n$ because the $e_k$'s do, and we can solve for the $e_k$'s in terms of the $g_k$'s.  So  $J_n = J_n(t)$  is generated by $ 1,g_1,...,g_{n-1} $
and we have a tower of algebras,        $J_1(t) \subset  J_2(t) \subset \dots$, ordered by inclusion.
The parameter $t$, which replaces the index $\lambda$, is the quantity now under investigation.
It's woven into the construction of the tower.  The algebra  $J_n(t)$  has defining relations:
\begin{eqnarray}
\label{equation:Jones algebra relations}
    g_i  g_k =g_k g_i \  \ {\rm if} \  \ |i-k| \geq 2, \  \  g_i g_{i+1} g_i =  g_{i+1} g_i g_{i+1}, \ \  g_i^2 = (t-1)g_i + t , \\
 1 + g_i + g_{i+1} + g_ig_{i+1} + g_{i+1}g_i +  + g_ig_{i+1}g_i   =  0. \nonumber
\end{eqnarray}
Of course there are braids lurking in the background.  If we rename the $g_i's$, replacing $g_i$ by
$\sigma_i$, and declare the $\sigma_i's$ to be generators of a group, then the first
two relations are defining relations in the group algebra $\complexes{\bf B}_n$.  The algebra $J_n(t)$ is thus a homomorphic image of the group algebra of the braid group.   Loosely speaking, braids are encountered in Operator Algebras because they encode the way in which each type II$_1$ factor $M_i$ acts on its subfactor $M_{i-1}$.   Braiding is thus involved in defining the associated extensions.   We shall see later, in $\S$\ref{subsection:additional representations} that a similar action, via group extensions,  can be used to define representations of ${\bf B}_n$.

To see the connection with knots and links, recall that since $M_n$ is type
II$_1$  it supports a unique trace, and since $J_n$ is a subalgebra it does too, by restriction. This trace is known as
a \underline{Markov trace}, i.e. it satisfies the important property:
\begin{equation}
\label{equation:Jones trace}
tr(w g_n) = f(t)tr(w)  \ \ \    {\rm if} \ \ \   w \in J_n,
\end{equation}
where $f(t)$ is a fixed function of $t$. Thus, for each fixed value of $f$ the trace is multiplied by a fixed scalar
when one passes from one stage of the tower to the next, if one does so by multiplying an arbitrary element of $J_n$ by
the new generator $g_n$ of $J_{n+1}$.
The Jones trace is nothing more or less than the 1-variable Jones polynomial \cite{Jones1987} associated
to the knot or link which is obtained from the closed braid.  We will have more to say about all this in $\S$\ref{subsubsection:Hecke algebras and polynomial invariants of knots}.

\subsubsection{Homotopy groups of spheres}
\label{subsubsection:braids in homotopy theory}
As before, let ${\bf P}_{n+1}$ denote the pure braid group on $n+1$ strands. For each $i = 1,\dots,n+1$ there is a natural homomorphism $p_i:{\bf P}_{n+1}\to {\bf P}_n$, defined by pulling out the $i^{th}$ strand. The group
of \underline{Brunnian} braids is
$$BR_{n+1} = \cap_{i=1}^{i=n+1} {\rm kernel}(p_i),$$
i.e. a braid is in $BR_{n+1}$  if and only if, on pulling out any
strand, it becomes the identity braid on $n$ strands.  Brunnian
braids have received some attention in knot theory.

Braids have played a role in homotopy theory for many years, most
particularly in the work of F. Cohen and his students (see for
example \cite{BCWW}), but during the past few years the connection
was sharpened when it was discovered that there is an embedding of
a free group $F_n$ in $P_{n+1}$ with the property that a
well-defined quotient of $BR_{n+1}\cap F_n$ (a little bit too
complicated to describe here) is isomorphic to $\pi_{n+1}S^2$.  It
remains to be seen whether new knowledge about the unidentified
higher homotopy groups of spheres can be obtained through the
methods of \cite{BCWW}.

\subsubsection{Robotics}
\label{subsubsection:braids in robotics} Our fourth example is an
application of configuration spaces to robotics. It shows the
braid group popping up in an unexpected way (until you realize how
natural it is).   Robots, or AGVs (automatic guided vehicles),
are required to travel across a factory floor that contains many
obstacles, en route to a goal position (e.g. a loading dock or an
assembly workstation). The problem is to design a control system
which insures that the AGVs not collide with the obstacles, or
with each other, and complete the task with efficiency with regard
to various work functionals. Here is how configuration spaces
appear: The underlying space in this simple example is the workspace floor $X$, from
which a finite set $\cO$ of obstacles are to be removed. The
configuration space of $n$ non-colliding AGVs is then precisely
$\cC_{0,n}(X - \cO)$.  More generally,  $X - \cO$ is replaced by a finite graph $Y$ , and the the braid group
${\bf B}_n$ by the braid group $\pi_1(\cC_{0,n}(Y))$ of the graph.  There is a vast literature on this subject;
we suggest \cite{Ghrist} by R. Ghrist, as a starter.

\subsubsection{Public key cryptography}
\label{subsubsection:braid groups and public key cryptography}

In this example braids are important for rather different reasons
than they were in our earlier examples.  In our earlier examples the
underlying phenomenon which was being investigated involved actual
braiding, albeit sometimes in a concealed way.  In the example
that we now describe particular properties of the braid groups
${\bf B}_n, \ n=1,2,3,\dots$, rather than the actual interweaving
of braid strands,  are used in a clever way to construct a new
method for encrypting data.

The problem which is the focus of `public key cryptography' will
be familiar to everyone: the security of our online
communications, for example our credit card purchases, our ATM
transactions, our cell phone conversations and a host of other
transactions that have become a part of everyday life in the
$21^{st}$ century,  The basic problem is to encrypte or translate
a secret message into a code that can be sent safely over a public
system such as the internet, and decoded at the receiving end by
the use of a secret piece of information known only to the sender
and the recipient, the `key'.  The problem that must then be
solved is to establish a private key that will be known only to
the sender and the recipient, who will then be able to exchange
information over an insecure channel.
 In recent years much work has been done on certain codes which are based upon the assumption that the word problem has polynomial growth as braid index $n$ is increased, whereas the conjugacy problem does not.
But in $\S$\ref{section:the word and conjugacy problems in the braid groups} we will review recent work on the word and conjugacy problems in the braid groups, and show that such an assumption seems problematic at best. See  $\S$\ref{section:the
word and conjugacy problems in the braid groups}, and in
particular the discussion in $\S$\ref{subsection:Braid groups and
public key cryptography} below.

\

\noindent {\bf Acknowledgements:} We thank Tahl Nowik, who
suggested, during a course that the first author gave on mapping
class groups, that techniques she had used for other purposes
could be adapted to give the proof we presented here of Theorem
\ref {theorem:isomorphism between bbn and MCG of punctured disc}.
We also thank  Robert Bell, David Bessis, Nathan Broaddus, John Cannon, Ruth Charney, Fred
Cohen, Patrick Dehornoy, Roger Fenn, Daan Krammer, Lee Mosher,
Luis Paris, Richard Stanley, Morwen Thistlethwaite and Bert Wiest for their help in
chasing down  facts and references, and helping us fill gaps in
our knowledge. We  would particularly like to thank all the
students who attended Math 661 in the Fall of 2003 at Cornell
University, who were gracious guinea pigs for large parts of this
article, especially Heather Armstrong, whose careful attention and
diligence significantly improved the manuscript, and Bryant Adams,
who suggested the proof of Lemma~\ref{lemma:step3} to us.

\newpage

\section{From knots to braids}

\label{section:from knots to braids}  In this chapter we will explore, for the benefit of readers who are new to the subject,  the foundations of the close relationship between knots and
braids. We will first describe the straightforward process of
obtaining a knot or a link from a given braid by  `closing' the braid. This leads us directly to formulate two
fundamental questions about knots and braids. First, is it always
possible to transform a given knot into a closed braid? This
question will be answered in the affirmative in Theorem
\ref{theorem:Alexander's theorem}, first proved by Alexander in
1928 in \cite{Alexander1923}. The correspondence between knots and
braids is clearly not one-to-one (for example, conjugate braids
yield equivalent knots), leading naturally to the second question:
which closed braids represent the same knot type? That question is
addressed in Theorem \ref{theorem:MT}, first formulated by A. Markov in \cite{Markov1935}, which gives `moves' relating any
two closed braid representatives of a knot or link, while
simultaneously preserving the closed braid structure.

Together, Theorems \ref{theorem:Alexander's theorem} and
\ref{theorem:MT} form the cornerstone of any study of knots via
closed braids, so we feel obliged to prove them.  Among the many
proofs that have been published of both over the years, we have
chosen ones that we like but which do not seem to have appeared in
any of the review articles that we know. The proof that we give of
Alexander's Theorem is due to Shuji Yamada \cite{Yamada1987}, with
subsequent improvements by Pierre Vogel \cite{Vogel1990}.  The
algorithm is elementary enough to be accessible to a beginner, and
has the advantage for experts of being suitable for programming.
The proof that we present of Markov's theorem is due to Pawel
Traczyk \cite{Traczyk1998}. It is relatively brief, as it assumes
Reidemeister's well-known theorem about the equivalence relation
on any two diagrams of a knot, Theorem \ref{theorem:Reidemeister} below,  building on methods introduced in the proof of Theorem \ref{theorem:Alexander's theorem}.

\subsection{Closed braids}\label{subsection:closedbraids}
\begin{figure}[htpb]
\centerline{\includegraphics[scale=.35, bb=66 306 496 604]{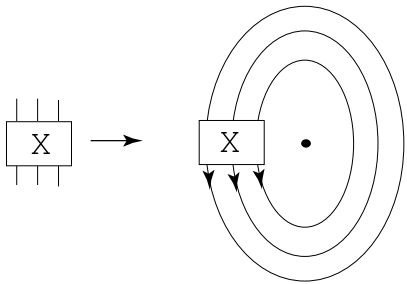}}
%\centerline{\includegraphics[scale=.35]{braid-closure.eps}}
\caption{The operation of closing a braid $X$ to form a closed braid}
\label{figure:braidclosure}
\end{figure}

For simplicity, let us begin with a planar diagram of a given
geometric braid. To obtain a knot or link, one simply `closes up'
the ends of the braid as in Figure~\ref{figure:braidclosure}.
The pre-image in $\reals^3$ of the `center point' shown in
Figure~\ref{figure:braidclosure} under the usual projection map is
called the \underline{axis} of the braid. (If one wishes to
consider the knot in $S^3$, then we include the point at infinity
so that the braid axis is an embedded $S^1$.) We then orient the
resulting knot or link in such a way that the strands of the braid
are all travelling counterclockwise about the braid axis. The knot
or link type resulting from performing this operation on a braid
$X$ is known as the \underline{closure} of $X$ and will be denoted
by $b(X)$.  The same notation may also refer to the particular
diagram as in Figure~\ref{figure:braidclosure}.

Equivalently, consider a knot $K \subset S^3$. Suppose there exists $A=
h(S^1)$ where $h$ is an embedding and $Z$ is unknotted in $S^3$
and contained in the complement of $K$. Suppose further that we
choose the point at infinity $\{ \infty \}$ to be in $A$ and,
using standard cylindrical coordinates $(\rho, \theta, z)$ on
$\reals^3$ , identify the resulting copy of $\reals \cong A - \{
\infty \}$ with the $z$-axis in $\reals^3 \cong S^3 - \{ \infty
\}$.  If we always have $d \theta / dt > 0$ as we travel about the
knot $K$ with an appropriate cylindrical parametrization, then we
say that $K$ is a \underline{closed braid} with respect to the axis $A$.   The closed braid diagram of Figure~\ref{figure:braidclosure} is then obtained by projection parallel to the direction defined by $A$ onto a plane that is orthogonal to $A$.

\subsection{Alexander's Theorem}\label{subsection:Alexander's Theorem}

As we just observed, it is a simple matter to obtain a knot or
link from a braid.  The classical theorem of J. Alexander
allows us to reverse this process, though not in a unique way:

\begin{theorem} {\bf (Alexander's Theorem \cite{Alexander1923})}
\label{theorem:Alexander's theorem}
Every knot or link in $S^3$ can be represented as a closed braid.
\end{theorem}

\pf Alexander's original proof was algorithmic, i.e. it gave an
algorithm for transforming a knot or link into closed braid form.
While it is straightforward, we do not know of any computer
program based upon it.   We shall give instead a rather different
and newer algorithm originally due to Yamada \cite{Yamada1987}, as
later improved by Vogel \cite{Vogel1990}. We like it for two
reasons: (1) It has a beautiful corollary (see Corollary \ref
{corollary:seifert cicles and braid index} below) which reveals
structure about knot diagrams that had not even been conjectured
by any of the experts before 1987, even though there was abundant
evidence of its truth;  (2) It leads, very easily, to an efficient
computer program for putting knots into braid form. In this regard
we note that when Jones was writing the manuscript
\cite{Jones1986}, which resulted in his award of the Fields medal,
he computed closed braid representatives for the 249 knots of
crossing number less than or equal to 10, constructing the first
table known to us of closed braid representatives of knots. His
list remains extremely useful to the workers in the area in 2004.
Yamada's work was not yet known when he did that work, and there
did not seem to him to be a good way to program the Alexander
method for a computer, so he calculated them one at a time by
hand.  The amount of work that was involved can only be
appreciated by the reader who is willing to try a few examples.

In order to prove Alexander's Theorem, we shall first present the
Yamada-Vogel algorithm for transforming a knot into closed braid
form in full, followed immediately by an illustrative example
(Example~\ref{example:52knot}). We shall then prove Alexander's
Theorem by showing that it is always possible to perform the steps
of the Yamada-Vogel algorithm for any given knot or link and that
the algorithm always leads to a closed braid.

The Yamada-Vogel algorithm draws on Seifert's well-known algorithm
for using a diagram of an oriented knot or link $K$ to construct a
Seifert surface for $K$ (see \cite{Rolfsen2003} or
\cite{Lickorish1997}, e.g., for a thorough treatment of Seifert's
algorithm), and we will need some related terminology. Let $C$ and
$C'$ be two oriented disjoint simple closed curves in $S^2$. Then
$C$ and $C'$ cobound an annulus $A$. We say that $C$ and $C'$ are
\underline{coherent} (or \underline{coherently oriented}) if $C$
and $C'$ represent the same element of $H_1 (A)$. Otherwise we say
that $C$ and $C'$ are \underline{incoherent}.  Following Traczyk
\cite{Traczyk1998}, we define the \underline{height} of a knot
diagram $D$, denoted $h(D)$, to be the number of distinct pairs of
incoherently oriented Seifert circles which arise from applying
Seifert's algorithm to $D$.  The height function gives us a useful
characterization of a closed braid: {\em a diagram $D$ represents
a closed braid if and only if $h(D) = 0$}.  (Recall that $D$ lives
in $S^2$.)

\

\noindent{\bf The Yamada-Vogel Algorithm}

\begin{enumerate}
\item  Let $D$ be a diagram of an oriented knot $K$.
Smooth all crossings of $D$ as in Seifert's algorithm  to obtain
$n$ Seifert circles $C_1, \ldots, C_n$.  Record each original
crossing with a signed arc: $(+)$ for a positive crossing (often called a \underline{right-handed} crossing),
$(-)$ for a negative (or \underline{left-handed}) crossing (see
Figure~\ref{figure:elem-bds}(iv)). The resulting diagram is the
\underline{Seifert picture} $S$ corresponding to the diagram $D$.
Note that any two circles joined by a signed arc in any Seifert
picture are necessarily coherent.    For an example that illustrates the construction of a Seifert diagram, see the passage from the bottom left to the top left sketches in Figure \ref{figure:52example}.
\begin{figure}[htpb!]
$$
\setlength{\unitlength}{0.05in}
\begin{picture}(0,0)(0,11)
\tiny
%Second picture
\put(-36.5,8){$+$} \put(-36.5,1.8){$+$} \put(-45,-3){$+$}
\put(-27.5,-3){$+$} \put(-36.5,-11){$+$}\put(-38,-4){$\alpha_1$}
%Third picture
\put(-9,2){$+$} \put(-9,-4){$+$} \put(-4,-10){$\alpha_2$}
\put(-5.25,-2){$-$} \put(1,-2){$+$} \put(5,4){$+$}
\put(5,-.25){$+$} \put(5,-4.5){$+$}
%Fourth picture
\put(34,4){$+$} \put(34,2){$+$} \put(34.5,-2.25){$+$}
\put(34,-5.75){$+$} \put(33,-10){$+$} \put(30.5,.25){$+$}
\put(37.5,0){$-$}\put(44.5,.25){$+$}\put(25.5,-4.5){$-$}
%Last picture
%Inner two
\put(19.5,-26.5){$+$} \put(14.5,-31.5){$-$}
%Middle five
\put(14.5,-21){$+$} \put(19,-21.5){$+$} \put(23,-28){$+$}
\put(20.5,-33.5){$+$} \put(10,-32){$+$}
%Outer two
\put(22.5,-35.5){$+$} \put(5.5,-26.5){$-$} \normalsize
\end{picture}
$$
\centerline{\includegraphics[scale=.5]{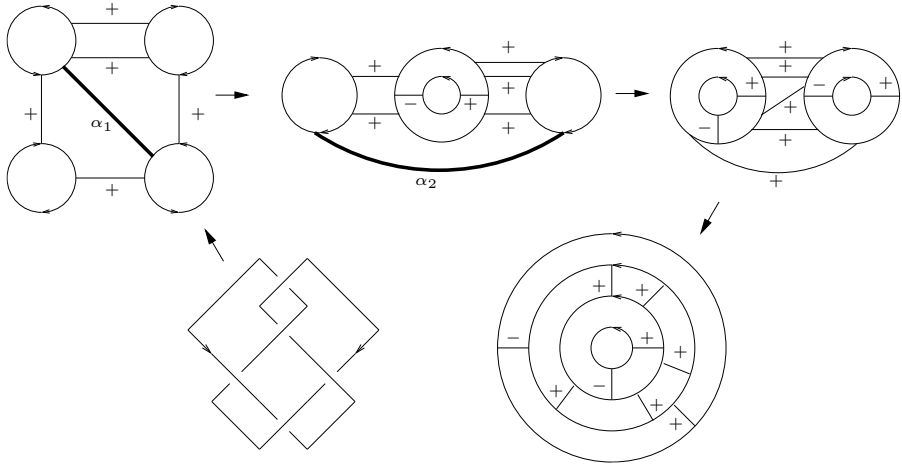}}
\caption{The Yamada-Vogel algorithm performed on the knot $5_2$.}
\label{figure:52example}
\end{figure}

\item  If $h(D) = 0$, the knot $K$ is already
in closed braid form, and we are done.
  If $h(D) > 0$, we can find a \underline{reducing arc}
$\alpha$, i.e., an arc joining an incoherent pair $C_i, C_j$ such
that $\alpha$ intersects $S$ only at its endpoints. Reducing arcs are illustrated as heavy black arcs in the example in Figure \ref{figure:52example}. A component of
$S^2\setminus S$ which admits a reducing arc is called a
\underline{defect region}. Perform a \underline{reducing move}
along $\alpha$, as shown in Figure~\ref{figure:reducingmove}, to
obtain a new Seifert picture $S'$ in which a pair of coherent
Seifert circles, $C_a$ and $C_z$, joined by two oppositely signed
arcs, replaces the incoherent pair $C_i, C_j$. The corresponding
move on the original diagram $D$ is a Reidemeister move of type II
in which we slide $C_i$ over $C_j$ in a small neighborhood of the
arc $\alpha$ to obtain a new diagram $D'$ with two new crossings.
Note that if we instead slide $C_i$ under $C_j$, we obtain the
same two new Seifert circles but the signs of the two new signed
arcs are now switched.

\begin{figure}[htpb!]
$$
\setlength{\unitlength}{0.05in}
\begin{picture}(0,0)(0,11)
\tiny \put(11.5,17){$\alpha$}\put(-2,18){$C_i$}\put(24,18){$C_j$}
\put(69,19){$+$}\put(69,12){$-$} \put(80,18){$C_z$}\put(71,15
){$C_a$} \put(67,15){$\cD_a$}\put(58,23){$\cD_z$} \normalsize
\end{picture}
\includegraphics[width=4in]{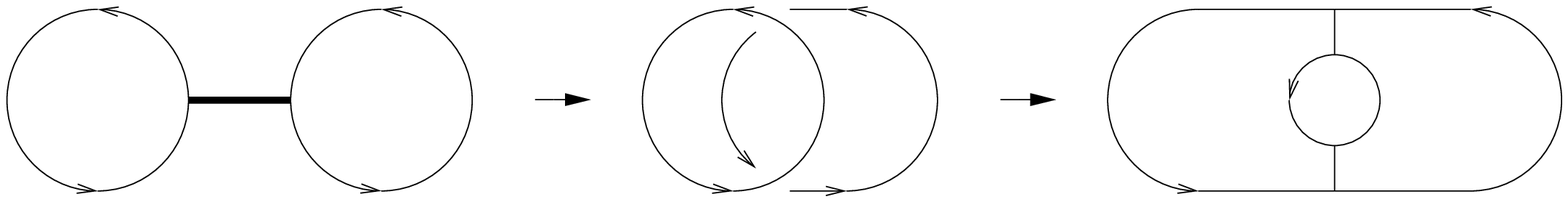}
$$
\caption{The local picture of a reducing move.}
\label{figure:reducingmove}
\end{figure}

\item  Continue performing reducing moves on incoherent pairs until
a diagram with height zero is obtained.
\end{enumerate}

\begin{example}\label{example:52knot}
{\rm  We apply the Yamada-Vogel algorithm to the diagram of the
knot $5_2$, pictured in Figure~\ref{figure:52example}.  (The
reader can check that this is the first example in the knot tables
of a knot diagram with height greater than zero.)
The first stage shows the Seifert
picture associated to the original diagram, consisting of 4
Seifert circles, with 5 signed arcs (all positive) recording the
original crossings. We see that the original knot diagram has
height 2. The figure shows a choice of reducing arc, $\alpha_1$,
joining one of the two pairs of incoherent circles.

In the third sketch of Figure~\ref{figure:52example}, we see the
new Seifert picture resulting from performing the reducing move
along $\alpha_1$.  Note that we have introduced two new crossings
of opposite sign.  We also see a new reducing arc, $\alpha_2$,
joining the only remaining pair of incoherent circles. We see in
the fourth sketch the Seifert picture with height zero resulting
from the second reducing move performed along $\alpha_2$.  At this
point, we are done, but in the final sketch we see a different
planar projection of the same Seifert picture which allows us
easily to read off a braid word associated to the knot: beginning
with the positive signed arc in the `twelve o'clock' position and
reading counter-clockwise, we see that the knot $5_2$ is
equivalent to $b(X)$, where $X = \sigma_2 \sigma_1^{-1} \sigma_2
\sigma_3^{-1} \sigma_2 \sigma_1 \sigma_2 \sigma_3 \sigma_2$.

We may learn several things from this simple example. Applying the
braid relations to the word defined by $X$, we see that $X =
\sigma_2 \sigma_1^{-1} \sigma_2 \sigma_3^{-1} \sigma_2 \sigma_1
\sigma_2 \sigma_3 \sigma_2
  = \sigma_2 \sigma_1^{-1} \sigma_2 \sigma_3^{-1}  \sigma_1 \sigma_2
\sigma_1 \sigma_3 \sigma_2
 = \sigma_2 \sigma_1^{-1} \sigma_2 \sigma_1  \sigma_3^{-1} \sigma_2
\sigma_3 \sigma_1 \sigma_2 \ \ \  = \sigma_2 \sigma_1^{-1}
\sigma_2 \sigma_1  \sigma_2 \sigma_3 \sigma_2^{-1} \sigma_1
\sigma_2 .$ Since this braid only involves $\sigma_3$ once, we may
`delete a trivial loop to get the 8-crossing 3-braid $\sigma_2
\sigma_1^{-1} \sigma_2 \sigma_1  \sigma_2  \sigma_2^{-1} \sigma_1
\sigma_2 ,$  so the algorithm did not give us minimum braid index.
The algorithm also does not give shortest words because our
8-crossing braid may be shortened to the 6-crossing 3-braid
$\sigma_2 \sigma_1^{-1} \sigma_2 \sigma_1^2 \sigma_2 .$   In fact,
6 is minimal, because the crossing number of a 3-braid knot must
be even, and this knot has no diagram with fewer than 5 crossings,
so 6 is minimal.  So we may deduce one more fact: when we use the Yamada-Vogel algorithm to change a knot which is not in closed braid form to one which is, the crossing number goes up. }
$\spadesuit$
\end{example}

To prove Alexander's Theorem, we first need to show
that a reducing move strictly decreases the height of a diagram.
This lemma is sometimes stated as `obvious' in the literature, but
the question arises frequently enough to warrant a short but
thorough argument.

\begin{lemma}\label{lemma:reducing}  Suppose a reducing move is
performed which transforms a diagram $D$ to a diagram $D'$.  Then
$h(D') = h(D) - 1$.
\end{lemma}

\noindent{\bf Proof.} Let $C_1, \ldots, C_n$ be the Seifert
circles in the Seifert picture $S$ corresponding to the diagram
$D$. Let $C_i, C_j$ be an incoherent pair. The union $C_i \cup
C_j$ separates the 2-sphere into three components: an annulus $A$
cobounded by $C_i$ and $C_j$, and two disks, $\cD_i$ and $\cD_j$
bounded by $C_i$ and $C_j$, respectively. Suppose that $A$ admits
a reducing arc $\alpha$.  A reducing move along $\alpha$ preserves
circles $C_p, p \neq i,j$ and replaces $C_i$ and $C_j$ with two
new circles, one of which necessarily bounds a disk $\cD_a$
containing no other Seifert circles in the new Seifert picture
$S'$. We denote this circle by $C_a$ (see
Figure~\ref{figure:reducingmove}). The other new circle, denoted
$C_z$, bounds a disk $\cD_z$ containing all Seifert circles
originally contained in the annulus $A$ in $S$.

\

To simplify the bookkeeping, we shall write $(C_r, C_s) = 1$ if
the pair $C_r, C_s$ is coherent, or else $(C_r, C_s) = -1$ if
$C_r, C_s$ are incoherent. Obviously, if $\{ p,q \}  \cap \{ i,j
\} = \emptyset$, then $(C_p, C_q)$ is unchanged by the reducing
move, so we need only consider the effect of the reducing move on
$(C_p, C_x)$, where $x = i$ or $x = j$ and $p \neq i,j$.   Now if
$C_p$ is contained in the annulus $A$ in $S$ (and hence in $\cD_z$
in $S'$), then clearly $(C_p,C_z) = (C_p,C_a)=(C_p,C_i) =
(C_p,C_j)$. Also, if $C_p \subset \cD_i$ in $S$, then $(C_p, C_z)
= (C_p, C_i)$ and $(C_p, C_a) = (C_p, C_j)$. Similarly, if $C_p
\subset \cD_j$ in $S$, then $(C_p, C_z) = (C_p, C_j)$ and $(C_p,
C_a) = (C_p, C_i)$.  Therefore the number of distinct incoherent
pairs $C_r, C_s$ in $S$ with $\{r,s\} \neq \{ i,j \}$ is equal to
the total number of distinct incoherent pairs in $S'$ of the form
$C_r, C_s$ with $\{r,s\} \neq \{ a,z \}$.  By construction,
however, we have replaced $(C_i, C_j) = -1$ with $(C_z, C_a) = 1$.
Thus $h(D') = h(D) -1$. \endpf

The previous lemma tells us that the Yamada-Vogel algorithm will always lead
to a diagram of height zero, i.e., a closed braid, as long as it
is always possible to perform Step 3. Therefore the following
lemma, whose proof was suggested to us by Bryant Adams, will
conclude the proof of Alexander's Theorem.

\begin{lemma}\label{lemma:step3}{\rm \cite{Yamada1987}}
Let $D$ be a knot or link diagram.  If $h(D) > 0$, then the
Seifert picture $S$ associated to $D$ contains a defect region.
\end{lemma}

\pf Each component $\cR$ of $S^2 \setminus S$ is a surface of
genus $0$ with $k \geq 1$ boundary components. Each boundary
component of $\cR$ is a union of some number of signed arcs
(possibly zero) and subarcs of Seifert circles.  We call the
collection of Seifert circles in $S$ which form part or all of a
boundary component of $\cR$ the \underline{exposed circles} of
$\cR$.

Let us now examine the possible ways in which $\cR$ could fail to
be a defect region. Certainly, if $\cR$ has one exposed circle,
then it is not a defect region. If $\cR$ has two exposed circles,
then $\cR$ is either a disk or an annulus.  In this case, if $\cR$
is a disk, then its two exposed circles are joined by at least one
signed arc; hence the two circles are coherent and $\cR$ is not a
defect region.  If $\cR$ is an annulus, either its two exposed
circles are incoherent, in which case it is a defect region, or
else they are coherent and $\cR$ is not a defect region.   If
$\cR$ has three or more exposed circles, then there is necessarily
one incoherent pair among them, and $\cR$ is a defect region.

Suppose that no component of $S^2 \setminus S$ is a defect region
and hence that each component is of one of the three types of
non-defect regions described above.
It is clear that we must have at least one region of the second
type, since otherwise $h(D)$ is clearly zero.  We can think of
such a region as lying between two nested, coherent circles joined
by at least one signed arc.

Let us now start with such a region, and try to build a diagram
with no defect regions.  We cannot add any circles in the annulus
cobounded by the two nested circles, since this necessarily gives
rise to at least one component with three or more exposed circles.
In fact, our only option is to add coherent circles which nest
with the original two circles (and as many signed arcs between
adjacent pairs as we like). However, such a diagram has height
zero. Therefore, $h(D) > 0$ implies that a defect region exists.  This completes the proof of Lemma
\ref{lemma:reducing}, and so also of Theorem \ref{theorem:Alexander's theorem}.  \endpf

The \underline{braid index} of a knot or link $K$ is the minimum
number $n$ such that there exists a braid $X \in {\bf B}_n$ whose
closure $b(X)$ represents $K$. (We note that it is also common to
refer to the \underline{index} of a braid or a closed braid,
meaning simply the number of its strands or the number of times it
travels around its axis, respectively.) It is clear that the
minimum number of Seifert circles in any diagram of a knot or link
$K$ is bounded above by the braid index of $K$. It is equally
clear from the Yamada-Vogel algorithm that the reverse inequality
holds. Thus we obtain the following corollary, which is due to Yamada  \cite{Yamada1987}. It seems remarkable
that it was not noticed long before 1987.

\begin{corollary}{\rm \cite{Yamada1987}}
\label{corollary:seifert cicles and braid index}
The minimum number of Seifert circles in any diagram of a knot or
link $K$ is equal to the braid index of $K$.
\end{corollary}
It also follows that we have a measure of the complexity of the
process of transforming a knot into closed braid form as follows.

\begin{corollary} {\rm (\cite{Traczyk1998}, \cite{Vogel1990})} Let $N$ denote the length of any
sequence of
reducing moves required to transform a diagram $D$ into closed
braid form.  Then we have:

\begin{equation}
N = h(D)  \leq  \frac{(n-1)(n-2)}{2}.
\end{equation}
where $n$ is the number of Seifert circles associated to $D$.
\end{corollary}

\begin{openproblem} {\rm It is an open problem to determine, among all regular diagrams for a given knot or link,  the minimum number of
Seifert circles that are needed.  By Corollary
\ref{corollary:seifert cicles and braid index} this is the same as
the minimum braid index, among all closed braid representatives of
a given knot or link.  We know of only one general result relating
to this problem, namely the Morton-Franks-Williams inequality of
\cite{Morton1986_2} and \cite{F-W}.  It will be discussed briefly
in $\S$\ref {subsubsection:Hecke algebras and polynomial
invariants of knots}.    The literature also contains an assorted
collection of ad-hoc techniques for determining  the braid index
of individual knots.  For example, see the methods used in
\cite{BM-stab-II} to prove that the 6-braid template in Figure
\ref{figure:templates} below actually has braid index 6, which
rests on the fact that non-trivial braid-preserving flypes always
have braid index at least 3. } $\clubsuit$
\end{openproblem}

\subsection{Markov's Theorem}
\label{subsection:Markov's Theorem}
To introduce the main goal of this section, we begin by recalling for the reader  Reidemeister's theorem, which dates from the earliest days of knot theory. It was assumed and used (as a folk theorem) long before anybody wrote down a formal statement and proof.

\begin{figure}[htpb!]
\centerline{\includegraphics[scale=.7,  bb=61 492 576 587] {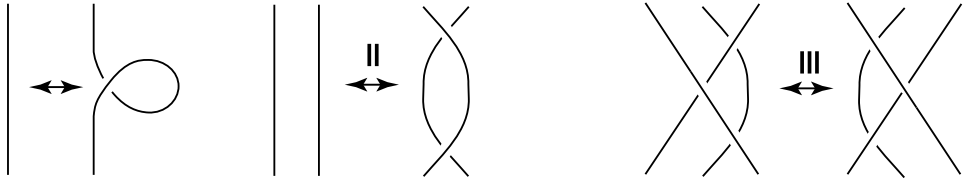}}
%\centerline{\includegraphics[scale=.7] {reidemeister.eps}}
\caption{The 3 Reidemeister moves. The 1, 2 or 3 strands in the
left sketch of each have arbitrary orientations, also we give only
one of the possible choices for the signs of the crossings, for
each move.} \label{figure:reidemeister}
\end{figure}

\begin{theorem}
 \label{theorem:Reidemeister}  {\bf (Reidemeister's Theorem)} Let $D, D^\prime$ be any two  (in general not closed braid) diagrams of the same knot or link $K$. Then there exists a sequence of diagrams $D = D_1\to D_2\to\cdots\to D_k = D'$ such that any $D_{i+1}$ in the sequence is obtained from $D_i$ by one of the three Reidemeister moves, depicted in Figure \ref{figure:reidemeister}.
\end{theorem}
\pf  We refer the reader to  \cite{Burde-Zieschang} for a complete proof.

\

Alexander's Theorem, proved in the last section, guarantees us
that closed braid representatives of a knot exist, but as
previously noted, they are certainly not unique.  Markov's
Theorem, first stated in \cite{Markov1935} with a sketch of a
proof, gives us a certain amount of control over different closed
braid representatives of the same knot. It asserts that any two
are related by a finite sequence of elementary moves and serves as
the analogue for closed braids of the Reidemeister Theorem for
knots.

\begin{figure}[htpb]
\centerline{\includegraphics[scale=.7, bb=58 421 575 571] {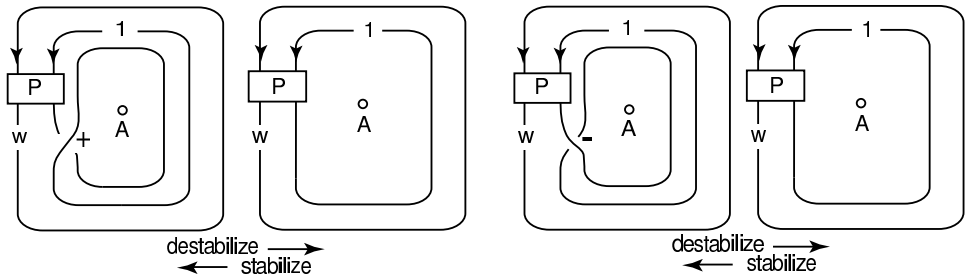}}
% \centerline{\includegraphics[scale=.7] {stab-destab.eps}}
\caption{The destabilization and stabilization moves.}
\label{figure:stab-destab}
\end{figure}

One of the moves of the Markov Theorem is \underline{braid
isotopy}. From the point of view of a topologist, braid isotopy
means isotopy of the closed braid, through braids, in the
complement of the braid axis. Morton has proved that if two braids
have closures that are braid isotopic, then they are conjugate in
${\bf B}_n$ \cite{Morton1978}. The other two moves that we need are
mutually inverse, and are illustrated in Figure
\ref{figure:stab-destab} as a move on certain  $(w+2)$-braids. We
call them \underline{destabilization} and
\underline{stabilization}, where the former decreases braid index
by one and the latter increases it by one.  The \underline{weight}
$w$ that is attached to one of the braid strands in Figure
\ref{figure:stab-destab} denotes that many `parallel' strands,
where parallel means in the framing defined by the given
projection.  The braid inside the box which is labelled P is an
arbitrary $(w+1)$-braid. Later, it will be necessary to
distinguish between positive and negative destabilizations, so we
illustrate both now.

\begin{theorem}
\label{theorem:MT}
{\bf (Markov's Theorem)} Let $X, X'$ be closed
braid representatives of the same oriented link type $K$ in oriented
3-space.  Then there exists a sequence of closed braid representatives
of $K$:
$$ X = X_1 \to X_2 \to \cdots \to X_r = X' $$
taking such that each $X_{i+1}$ is obtained from $X_i$ by either {\rm (i)} braid isotopy,
or {\rm (ii)} a single stabilization or destabilization.
\end{theorem}

We call the moves of Theorem \ref{theorem:MT} \underline{Markov
moves}, and say that closed braids that are related by a sequence
of Markov moves are \underline{Markov-equivalent}.

Forty years after Markov's theorem was announced,  the first
detailed proof was published in \cite{Birman1974}.  At least 5
essentially different proofs exist today.   See for example
\cite{Morton1986}, in which Morton gives his beautiful threading
construction for knots and braids which also yields an alternate
proof of Alexander's Theorem.   Here we shall present a proof due to Pawel Traczyk \cite{Traczyk1998}. It
begins with Reidemeister's theorem, and uses the circle of ideas
that were described in the previous section, and so it is
particularly appropriate for us.

\

\pf We are given closed braids $X, X'$ which represent the same
oriented knot type $K$.  Without loss of generality we may assume
that $X$ and $X'$ are defined by closed braid diagrams $Y, Y'$ of
height $h(Y) = h(Y') = 0$.  By Theorem \ref{theorem:Reidemeister}
we know there is a sequence of knot diagrams $Y = Y_1 \to Y_2 \to
\cdots \to Y_k = Y'$, where in general $h(Y_i) \geq 0$ for
$i=2,\dots,k-1$, such that any two diagrams in the sequence are
related by a single Reidemeister move of type I, II or III. The
first step in Traczyk's proof is to reduce the proof to sequences of knot diagrams which are related by Yamada-Vogel reducing moves:

\begin{lemma}
\label{lemma:reducingtor}
 It suffices to prove Theorem \ref{theorem:MT} for closed
braid diagrams $Y, Y'$ which are related by sequences  $Y = Y_1
\to Y_2 \to \cdots \to Y_q = Y'$ with the properties {\rm (i)} $h(Y) =
h(Y') = 0$, {\rm (ii)} $h(Y_i)>0$ for $i = 2,\dots,q-1$, and {\rm (iii)}
$Y_{i+1}$ is obtained from $Y_i$ by a single Yamada-Vogel reducing
move or the inverse of a reducing move.
\end{lemma}

\pf  We may always assume that the
diagrams $Y_2,\dots,Y_{q-1}$ have height $>0$, for if not we simply replace the given
sequences by the subsequences joining any two intermediate
diagrams of height zero.

We say that a Reidemeister move is \underline{braid-like} if the
strands that are involved in it are locally oriented in a coherent
fashion, as they would be if the diagram is a closed braid.  In
particular, any Reidemeister move of type I is braid-like. To
begin the proof of the lemma, we establish a somewhat weaker
result: we claim that we can get from $Y$ to $Y'$ via a finite
sequence of the following four types of moves and their inverses:

$\bullet$ a braid-like Reidemeister move of type I, denoted type I$^b$,

$\bullet$ a braid-like Reidemeister move of type II, denoted type
II$^b$,

$\bullet$ a braid-like Reidemeister move of type III, denoted type
III$^b$,

$\bullet$ a Yamada-Vogel  reducing move, denoted type $\cY$.

\noindent To prove the claim, it is enough to show that
non-braid-like Reidemeister moves, which we denote by the symbols
I$^{nb}$, II$^{nb}$ and III$^{nb}$, can be achieved via a finite
sequence of moves of type I$^b$, II$^b$, III$^b$ and $\cY$.  To
prove this, we examine the cases  I$^{nb}$, III$^{nb}$ and
II$^{nb}$ in that order:

\be

\item  As previously noted, any type I Reidemeister move is of type I$^b$.

\begin{figure}[htpb!]
$$
\setlength{\unitlength}{0.05in}
%\begin{picture}(0,0)(0,11)
%\small \put(10,22){$r_1$} \put(66,23){$r_2^{-1}$}\normalsize
%\end{picture}
\includegraphics[width=4in]{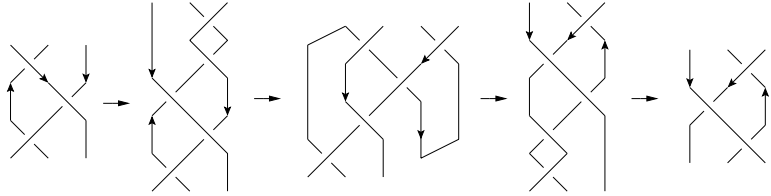}
$$
\caption{Replacing moves of type III$^{nb}$}
\label{figure:replacingIIInb}
\end{figure}

\item   A type III$^{nb}$ Reidemeister move involves three arcs of the knot
or braid.  There are many different cases, depending on the local
orientations and the signs of the 3 crossings, but they are all
similar.  One of the possible cases is given by the first and last
sketches of Figure~\ref{figure:replacingIIInb}, where an arc
passing under a crossing formed by the other two arcs is locally
oriented opposite to the other two strands.  The replacement
sequence that is given in Figure~\ref{figure:replacingIIInb} shows
that our  type III$^{nb}$ Reidemeister move can be achieved by a
sequence consisting of a type II$^{nb}$ move, an isotopy, a type
III$^b$ move and finally another type II$^{nb}$ move. We leave the
other type III$^{nb}$ cases to the reader, and we have reduced to
the case of moves of type II$^{nb}$.

\begin{figure}[htpb!]
$$
\setlength{\unitlength}{0.05in}
%\begin{picture}(0,0)(0,11)\small
%\put(8,19){$a_1$} \put(27,19){$a_2$} \put(52,19){$r_1$}
%\put(71,19){$r_2^{-1}$} \normalsize
%\end{picture}
\includegraphics[width=4in]{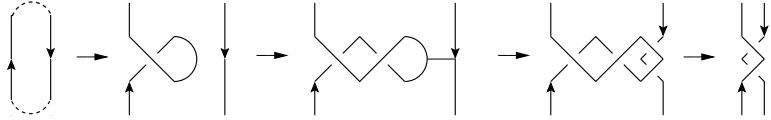}
$$
\caption{Replacing moves of type II$^{nb}$}
\label{figure:replacingIInb}
\end{figure}

\item  A move of type II$^{nb}$ may be regarded as
a move of type $\cY^\pm$ if the arcs that are involved belong to
distinct Seifert circles, so we only need to handle the case where
they are subarcs of the same Seifert circle. This is done in
Figure~\ref{figure:replacingIInb}, where it is shown that the move
can be replaced by two moves of type I$^b$ (which create two new
Seifert circles) followed by a move of  type $\cY$ and another of
type $\cY^{-1}$. This proves the claim.\ee

We are thus reduced to the case in which each diagram in the
sequence taking $Y$ to $Y^\prime$ is either type I$^b$, II$^b$,
III$^b$ or $\cY^\pm$.  To complete the proof of
Lemma~\ref{lemma:reducingtor}, let $t$ be a braid-like
Reidemeister move to be performed on diagram $Y_i$.  Suppose that
$h_i = h(Y_i)
> 0$. Then we can find a sequence of reducing moves $r_1, \ldots,
r_{h_1}$ such that $r_{h_i} \circ \cdots \circ r_1 (Y_i)$ is a
braid and such that the associated reducing arcs $\alpha_1,
\ldots, \alpha_{h_i}$ are each disjoint from the region in which
$t$ is to be performed.  Thus each reducing move $r_j$ commutes
with $t$, and we can replace $t$ with its `conjugate' $r_1^{-1}
\circ \cdots \circ r_{h_i}^{-1} \circ t \circ r_{h_i} \circ \cdots
\circ r_1$, so that $t$ is now performed at height $0$, i.e., on a
braid.

If $t$ is of type II$^b$ or III$^b$, then we are done, since a
braid-like move of type II or type III performed on a braid is a
braid isotopy.  If $t$ is of type I$^b$, it is a stabilization (up
to isotopy) only if it is performed on the braid strand `nearest'
to the braid axis. However, it is not hard to see how to realize a
type I$^b$ move on an arbitrary strand in terms of Markov moves:
simply push the strand under the others via type II$^b$ moves, and
perform the required stabilization in a neighborhood of the braid
axis. Note that to pass the resulting `kink' back under a
neighboring strand in a braid requires first a type III$^b$ move
followed by a type $\cY^{-1}$ move (the `kink' is always its own
Seifert circles, so two distinct circles are necessarily
involved). Thus we can return the strand with the `kink' in it
back to its original position by repeated applications of this
two-step process. Since we have just seen that any type II$^b$ or
type III$^b$ move can be realized  by a finite sequence of type
$\cY^\pm$ moves and braid isotopies, this means that a type I$^b$
move can also be replaced by a finite sequence of type $\cY^\pm$
moves and braid isotopies.  We can handle inverse moves of type
I$^b$ in a similar fashion

We have thus replaced our original sequence relating $Y$ to
$Y^\prime$ by a new one which is in general much longer, but which
consists entirely of Markov moves (performed, by definition, on
braid diagrams, i.e., on diagrams of height zero) and moves of
type $\cY^\pm$.  To be precise, our original sequence from $Y$ to
$Y'$ may be replaced by a sequence of the form $Y= Y_0, \ldots,
Y_{a_1}, \ldots, Y_{a_2}, \ldots,  Y_{a_n}=Y'$, where $h(Y_{a_i})
= 0$ for all $i$ and in each subsequence $Y_{a_i}, \ldots,
Y_{a_{i+1}}$, either (1) each diagram in the subsequence has
height zero and adjacent diagrams are related by a single Markov
move, or (2) all the intermediate diagrams have strictly positive
height and adjacent diagrams are related by a single move of type
$\cY$ or $\cY^{-1}$. Therefore in order to prove Markov's theorem,
it suffices to consider only sequences of the second type, and the
proof of Lemma~\ref{lemma:reducingtor} is complete. \endpf

\begin{remark}
\label{remark:Reidemeister and Markov}
{\rm The astute reader will have noticed the following:  we have eliminated Reidemeister moves completely (they will not appear in the arguments that follow), nevertheless they played an important role already.  We started with a sequence relating the given braids $X$ and $X'$ that consisted entirely of Reidemeister moves. We replaced it with a sequence of reducing moves and braid-like Reidemeister moves. The latter are in general not applied to diagrams of height zero, but we changed them to apply to diagrams of height zero. That is the moment when Traczyk's braid-like Reidemeister moves were changed to Markov moves.  The modified sequence from $X$ to $X'$ has changed to a series of subsequences, each of which starts with a closed braid and ends with a closed braid, after which the ending braid is modified (by Markov moves) to a new closed braid, which is the initial closed braid in the next subsequence.
As will be seen,  Markov moves will not be used explicitly  again until the proof of Lemma  \ref{lemma:irred-dpairs}, where they are used (without the help of Reidemeister's Theorem) to relate very special closed braid diagrams.}
\end{remark}

Consider now a sequence of diagrams $Y = Y_1, \ldots, Y_n = Y'$
satisfying the criteria of Lemma~\ref{lemma:reducingtor}.  We note
that, as in the above proof, we shall not in general distinguish
between a diagram and its associated Seifert picture.  Thus we
shall make reference to `a Seifert circle in the diagram $Y_i$',
for example, meaning a Seifert circle in the Seifert picture
associated to $Y_i$.  In fact, we can think of each circle in a
Seifert picture as forming part of the associated diagram, except
in a small neighborhood of signed arcs, which correspond to
crossings.  Since reducing arcs avoid signed arcs, there is no
ambiguity when referring to `a reducing arc in a diagram'.

We now wish to consider the graph of the height function on our
sequence.  The graph will begin and end at height zero; each
`step' in between will either take us up 1 or down 1 since we have
reduced to the case where all moves are reducing moves (or their
inverses). We will examine local maxima in the height function.
Let $Y(r), \hat{Y}, Y(s)$ be three consecutive diagrams in our
sequence such that the height function has a local maximum at $\hat{Y}$.
In other words, we have two reducing moves $r,s$ with
corresponding arcs $\alpha_r, \alpha_s$ in $\hat{Y}$ such that
reducing $\hat{Y}$ along $\alpha_r$ (resp. $\alpha_s$) results in
the diagram $Y(r)$ (resp. $Y(s)$), and it makes sense discuss
$\alpha_r\cup\alpha_s$. We will call such a triple $\{Y(r),
\hat{Y}, Y(s)\}$ a \underline{peak} in the height function of our
sequence. We define the \underline{height of the peak} to be
$h(\hat{Y})$ and define the \underline{height of the sequence} to be the
maximum value attained by the height function on the sequence, in
other words, the maximum over the height of all the peaks in the
sequence.  In order to prove Theorem~\ref{theorem:MT}, we are
going to induct on the height of the sequence.

\begin{lemma}
\label{lemma:disjoint} We may assume that the reducing arcs
involved in any peak in the height function of our sequence are
disjoint.  Further, the adjustments in our sequence of reducing
moves which are required preserve the height of the sequence.
\end{lemma}

\pf   Let $\{Y(r), \hat{Y}, Y(s)\}$ be a peak in the height
function with associated reducing arcs $\alpha_r$ and $\alpha_s$.
We may always assume that the arcs intersect transversally and
minimally. Suppose that $|\alpha_r \cap \alpha_s| = n \geq 2$.  By
smoothing out one or more of the points of intersection, we find a
new reducing arc $\alpha_{r'}$ with the same endpoints as
$\alpha_r$ such that $\alpha_r \cap \alpha_{r'} = \emptyset$ and
$|\alpha_{r'} \cap \alpha_s| < n$. We can then replace the given
peak $\{Y(r), \hat{Y}, Y(s)\}$ with two consecutive peaks $\{Y(r),
\hat{Y}, Y_{r'})$ and $(Y_{r'}, \hat{Y}, Y(s)\}$.  We call this
procedure \underline{inserting the reducing operation $r'$} at
$\hat{Y}$, and it essentially amounts to replacing one peak with
two peaks of the same height.  In this way, we continue on until
the intersection numbers of all adjacent pairs is at most $1$.

\begin{figure}[htpb!]
\centerline{\includegraphics[scale=.75, bb=243 444 404 597] {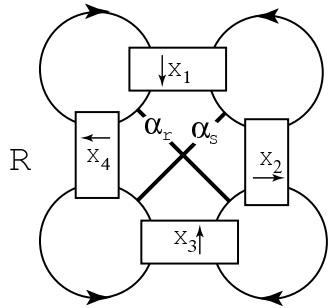}}
% \centerline{\includegraphics[scale=.75] {irreducible1.eps}}
\caption{The case where two reducing arcs at a peak  intersect
once.}
 \label{figure:irreducible1}
\end{figure}

Now suppose the arcs $\alpha_r, \alpha_s$ associated to a given
peak $\{Y(r), \hat{Y}, Y(s)\}$ have intersection number $1$. If
there exists a reducing arc $\alpha_t$ such that $\alpha_t \cap
\alpha_r = \alpha_t \cap \alpha_s = \emptyset$, then we can insert
the reducing operation $t$ at $\hat{Y}$ to produce two peaks, each
with a disjoint pair of associated reducing arcs.
Suppose that the defect region which supports $\alpha_r$ and
$\alpha_s$ contains no third reducing arc which is disjoint from
both $\alpha_r$ and $\alpha_s$.  There is only one
possible arrangement for such a defect region, shown in
Figure~\ref{figure:irreducible1}.

With this arrangement, the reducing
arcs $\alpha_r$ and $\alpha_s$ must act on four distinct Seifert
circles, possibly joined by signed arcs.   We will now examine the
region labelled $R$ which lies `outside' of the circles involved
in our defect region and the signed arcs which join them. If $R$
contains a Seifert circle, then it is easy to check that there
must be a region somewhere in the diagram with three exposed
circle. Such a region, as observed in the previous section, is
necessarily a defect region which would contain a reducing arc
disjoint from both $\alpha_r$ and $\alpha_s$.  If, on the other
hand, $R$ contains no Seifert circles, then it contains no
signed arcs either, since all possible signed arcs between the
four exposed circles of region $R$ already appear in
Figure~\ref{figure:irreducible1}. Thus we can join either pair
of diagonally opposed circles by a reducing arc in $R$.

We conclude that if a peak $\{Y(r), \hat{Y}, Y(s)\}$ has
associated arcs of intersection 1, we can always find a third
reducing arc $\alpha_t$ so that we can insert the reducing
operation $t$ at $\hat{Y}$, thereby replacing the original peak
with two peaks, each having a disjoint pair of reducing arcs.
Since the operation of inserting a reducing operation at a peak
preserves the height of the sequence, this finishes the proof of
the lemma.
\endpf

Thanks to the previous lemma, we can assume from now on that each
peak in the graph of our height function corresponds to a disjoint
pair of reducing arcs.  Before we can state the next lemma, which
concerns the peaks in the height function, we need to introduce
a few some concepts. Note that when the reducing arcs involved
in a peak $\{Y(r), \hat{Y}, Y(s)\}$ are disjoint, the
reducing moves commute, that is, they can be performed in
either order, starting with the diagram $\hat{Y}$ and resulting in
the same diagram $Y'$. Further, as long as the reducing arcs
$\alpha_r, \alpha_s$ act on 3 or 4 distinct Seifert circles, we
can perform the two reducing moves in either order with the same
result.  Observe that, since the reducing arcs $\alpha_r,
\alpha_s$ are disjoint in the diagram $\hat{Y}$, it makes sense to
talk about the arc $\alpha_s$ (resp. $\alpha_r$) in the context of
the diagram $Y(r)$ (resp. $Y(s)$) obtained by reducing $\hat{Y}$
along $\alpha_r$ (resp. $\alpha_s$). In this case we say we have a
`commuting pair' of reducing moves associated to the
peak.  If a peak $\{Y(r), \hat{Y}, Y(s)\}$ has a commuting pair,
then we may replace it by a `valley', that is, a subsequence
$\{Y(r), Y', Y(s)\}$ where $h(Y') = h(Y) - 2$ and $Y' = Y(s \circ
r) = Y(r \circ s)$ is the result of reducing $Y(r)$ along
$\alpha_s$, or equivalently reducing $Y(s)$ along $\alpha_r$. Thus
we can eliminate any peak corresponding to a commuting pair (such
a peak necessarily has height at least 2).

In the case where two reducing arcs at a peak act on the same 2
circles, then after one move is performed, the second Reidemeister
move will no longer be a reducing move; we call this a
`non-commuting pair' of reducing moves. Let $\{Y(r),
\hat{Y}, Y(s)\}$ be a peak corresponding to a non-commuting pair
of reducing arcs, and let $C_1, C_2$ be the two Seifert circles
involved. Suppose there is a reducing arc $\alpha_t$ such that
$\alpha_t \cap \alpha_r = \alpha_t \cap \alpha_s = \emptyset$ and
such that $t$ involves a circle other than $C_1$ or $C_2$, then we
can insert $t$ at $\hat{Y}$ to replace our peak $\{Y(r), \hat{Y},
Y(s)\}$ with two new peaks with commuting pairs of reducing
moves:$\{Y(r), \hat{Y}, Y(t) \}$ and $\{Y(t),\hat{Y},Y(s)\}$.  As
above, we now replace each peak with a `valley': $\{Y(r), Y',
Y(t)\}$ and $\{Y(t), Y'', Y(s)\}$, respectively, where $Y' = Y(t
\circ r) = Y(r \circ t)$ is the diagram resulting from reducing
$Y(r)$ by $t$ (or equivalently, from reducing $Y(t)$ by $r$, and
$Y''$ is the diagram resulting from reducing $Y(s)$ by $t$ (or
equivalently, from reducing
$Y(t)$  by $s$). Again, this implies
that the height of the original non-commuting peak $\{Y(r),
\hat{Y}, Y(s)\}$ was at least 2. Thus, we can replace such a peak
with peaks of strictly smaller height and repeat this process
until all peaks either have height 1 or do not admit such a
reducing arc $\alpha_t$ as above; we call a peak of the latter
type \underline{irreducible}.

We can now state the next lemma:

\begin{lemma}\label{lemma:noncommutingpeaks}
We may assume that any peak in the height function of our sequence
either has height 1 or is irreducible.
\end{lemma}

\pf  The proof is clear. In the discussion that preceded the statement
of the Lemma we defined a peak to be irreducible in such a way that it
subsumed all possibilities which did not allow us to reduce to height 1.
\endpf

\noindent We note that each condition of
Lemma~\ref{lemma:noncommutingpeaks} necessarily implies the
non-commuting condition.

\begin{lemma}\label{lemma:heightonepeaks}
We may assume that no peaks in the height function of our sequence
have height 1.
\end{lemma}

\pf Let $\{Y(r),\hat{Y}, Y(s)\}$ be a peak of height 1.  We recall
that height 1 implies that the two reducing arcs are
non-commutative and hence involve precisely two circles.
it is an easy exercise to show that these two circles must live
either on the `inside' or on the `outside' of a band of circles,
and that $\alpha_r$ and $\alpha_s$ are in fact equivalent as
reducing moves.  Therefore the diagram $Y(r)$ is equivalent to
$Y(s)$ and we can simply eliminate this peak from our sequence.
\endpf

It remains to deal with irreducible peaks in the height function
of our sequence.  Fortunately, it turns out that they only occur
in a very particular way.   To describe the particular way, define
a \underline{ weighted Seifert circle} in a manner which is similar
to the weights that we attach to closed braid diagrams, e.g. as in Figure
\ref{figure:stab-destab}. That is, A Seifert circle with weight $w$
attached means  a collection of $w$ coherently oriented, nested, parallel
Seifert circles.  We use the term \underline{band} for a Seifert circle
with an attached weight.

\begin{figure}[htpb!]
\centerline{\includegraphics[scale=.75, bb=56 499 534 683] {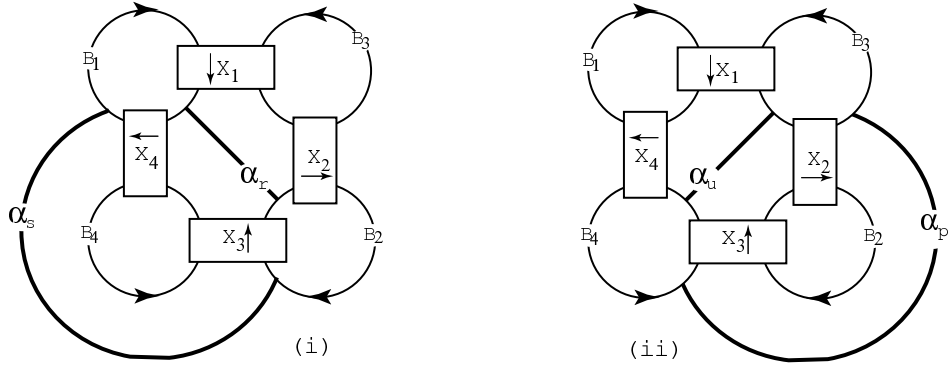}}
%\centerline{\includegraphics[scale=.75] {irreducible2.eps}}
\caption{(i) The diagram corresponding to to the irreducible peak
$\{ Y(r), \hat{Y}, Y(s) \}$. (ii) Two additional reducing arcs used
to eliminate irreducible peaks.   Each
pair $\alpha_r, \alpha_s$ and $\alpha_{p}, \alpha_{u}$  is a possible non-commuting
pair associated to the peak.
The blocks $X_1,X_2,X_3,X_4$ indicate (possibly empty) collections of
signed arcs.}
\label{figure:irreducible2}
\end{figure}

\begin{lemma}\label{lemma:fourbands}
If $\{Y(r), \hat{Y}, Y(s)\}$ is an irreducible peak in the height
function of our sequence, then the diagram $\hat{Y}$ contains at
most four bands, arranged as in Figure~\ref{figure:irreducible2}.
\end{lemma}

\pf  Let $\{Y(r),\hat{Y}, Y(s)\}$ be an irreducible peak, and let
$C_1, C_2$ be the two circles involved in the reducing arcs
$\alpha_r$ and $\alpha_s$.  See Figure \ref{figure:irreducible2}(i).
 For $i = 1,2$, let $D_i$ denote the
disk bounded by $C_i$ in $S^2$ which does not contain the reducing
arcs $\alpha_r$ and $\alpha_s$. Then $S^2 \setminus (D_1 \cup D_2
\cup \alpha_r \cup \alpha_s)$ has two components. By assumption,
neither component can contain a defect region (or else we could
find a reducing arc $\alpha_t$ as above, contradicting the
irreducibility of our peak). Thus if either component contains any
Seifert circles, the circles must form a band, oriented oppositely
to $C_1$ and $C_2$. We allow the possibility that the weight of
these bands is zero. The same reasoning shows that $D_i$ cannot
contain a defect region for $i= 1,2$ and hence that $C_i$ must be
the outer circle of a band. It is of course possible that some
braiding takes place between adjacent coherent circles, as
indicated in Figure~\ref{figure:irreducible2}(i). The braids joining the
various bands are labelled $X_i$. Thus we have a diagram of the
given form, in which each circle pictured represents a band, and
the lemma is proved.\endpf
The following lemma will allows us to replace irreducible peaks in
the height function of our sequence with peaks of strictly smaller
height.

\begin{lemma}\label{lemma:irred-dpairs}
Let $\{Y(r), \hat{Y}, Y(s)\}$ be an irreducible peak of height
$n+1$ in the height function of our sequence.  Then there exist
sequences of diagrams $Y (r)= Y^r_1 , \ldots, Y^r_n = Y(p \circ
r)$ and $Y(s) = Y^s_1, \ldots, Y^s_n = Y(u \circ s)$ such that
$Y^r_{i+1}$ (resp. $Y^s_{i+1}$) is obtained from $Y^r_i$ (resp.
$Y^s_i$) by a reducing move and such that $h(Y(p \circ r)) =
h(Y(u\circ s) = 0$ and $Y(p \circ r)$ and $Y(u \circ s))$ are
Markov equivalent.
\end{lemma}

Before discussing the proof of this lemma, we show how to use it
to prove the Markov theorem. Using
Lemmas~\ref{lemma:reducingtor}-\ref{lemma:fourbands}, we have
reduced the proof of Markov's Theorem to the situation of two
closed braid diagrams $X$ and $X'$ related by a sequence of
diagrams related by reducing moves (and their inverses) such that
the height of each intermediate diagram is strictly positive and
such that any peak in the height function of the sequence is
irreducible (of height at least 2).  Let $\{Y(r), \hat{Y}, Y(s)\}$
be an irreducible peak of height $n$ in the height function of our
sequence. By Lemma~\ref{lemma:irred-dpairs}, we can replace this
subsequence with a subsequence of strictly smaller height
(possibly including a subsequence entirely at height 0 related by
Markov moves, which we can remove from consideration as before).
  In doing
so, we create new peaks whose height is {\em strictly lower} than
the height of the peak being replaced.  If we perform this
operation at every irreducible peak, then we obtain a new sequence
relating our closed braid diagrams $Y, Y'$.

The new peaks may or may not be irreducible; in fact, their
corresponding arcs may not even be disjoint, but now we are back
in the same situation as we were before Lemma~\ref{lemma:disjoint}
except that we are starting with a sequence of lower height.  Thus
by induction on the height of the sequence (with base case
provided by Lemma~\ref{lemma:heightonepeaks}), we can replace any
sequence as described in Lemma~\ref{lemma:reducingtor} with a
sequence consisting entirely of diagrams of height zero.  This
completes the proof of Theorem~\ref{theorem:MT}, modulo the
proof of Lemma~\ref{lemma:irred-dpairs}.\endpf

\noindent{\bf Sketch of the Proof of Lemma~\ref{lemma:irred-dpairs}.}
By Lemma~\ref{lemma:fourbands}, the diagram $\hat{Y}$ has at most
four bands which we label $B_1, B_2, B_3$, and $B_4$, joined by
four (possibly trivial) braids $X_1, X_2, X_3, X_4$, as in
Figure~\ref{figure:irreducible2}(i) and (ii).
Let $w_i$ denote the weight of the
band $B_i$. Note that, e.g., the braid $X_1$ has $w_1 + w_3$
strands, and similarly for the other $X_i$. If $w_i$ or $w_j$ is
greater than 1, then a reducing arc joining $B_i, B_j$ is
understood to indicate a sequence of $w_i w_j$ reducing moves.
There are several ways to construct such a sequence; we adopt the
convention that we choose reducing moves in such a way that the
strands of one band all slide under or else all slide over the
strands of the other band.

\begin{figure}[htpb!]
\centerline{\includegraphics[scale=.80, bb=127 429 302 616] {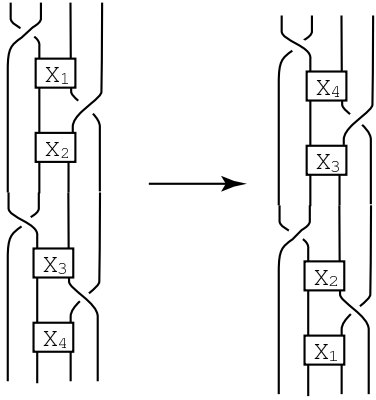}}
%\centerline{\includegraphics[scale=.80] {braids-two.eps}}
\caption{Two braids whose respective closures are the result of
performing the two pairs of reducing moves $r,p$ and $s,u$
indicated in Figure~\ref{figure:irreducible2}. \label{figure:twobraids}}
\end{figure}

Referring again to Figure~\ref{figure:irreducible2}, we first perform
the reducing move $r$ along the arc $\alpha_r$ and then reduce
again via the arc $\alpha_{p}$ involving $B_3$ and $B_4$, giving
us some number of reducing moves depending on the weights of the
bands involved.  The resulting diagram $Y(p \circ r)$ is the
closure of the first braid shown in Figure~\ref{figure:twobraids},
up to a choice of arcs passing over or under in the various
Reidemeister moves of type II.  This gives us the first sequence
of the lemma.

To get the second sequence of the lemma, we begin instead with the
reducing move $s$ along the arc $\alpha_s$ and then reduce again
via the arc $\alpha_{u}$.  The resulting diagram $Y(u \circ s)$
is the closure of the second braid shown in
Figure~\ref{figure:twobraids}.

We have found our two sequences of reducing moves, and we have
reduced the proof of Markov's Theorem to just one specific
calculation, namely, showing that the closed braids $Y(p \circ
r)$ and $Y(u \circ s)$ corresponding to the two braids of
Figure~\ref{figure:twobraids} are M-equivalent. The two braids can
be related by a sequence consisting of several braid isotopies
as well as two stabilizations and two destabilizations; the reader
is referred to \cite{Traczyk1998} for the details of this
calculation.  This concludes the proof of Lemma \ref{lemma:irred-dpairs}, and so also of
Theorem \ref{theorem:MT}.
\endpf

\begin{remark}
{\rm   Our choice of `over' or `under' in the reducing
moves $p$ and $u$ leads to possible ambiguity, but the
various diagrams which would result from different choices are all
related by
`exchange moves', which are defined and
discussed at the beginning of the next section.
We will show in the next section that an exchange move
replaces a sequence of 4 Markov moves: a braid isotopy, a
stabilization, a second braid isotopy and a destabilization.  }
\end{remark}

\

The essential groundwork has been laid regarding the connection
between braids and knots.  From this point on, all knot diagrams
will be assumed to be in closed braid form, i.e. in the form
indicated in Figure \ref{figure:braidclosure}.  In the sections
that follow we will examine many consequences.

\newpage

\section{Braid foliations}
\label{section:braid foliations} We begin our study of new results
which relate to the study of knots via closed braids by presenting
some results that use the theory of braid foliations of a Seifert
surface bounded by a knot which is represented as a closed braid.
We will develop three applications of braid foliations: The first
is Theorem \ref{theorem:mtws}, in $\S$\ref{subsection:mtws}. We
give an essentially complete proof of the `Markov theorem without
stabilization' (MTWS) in the special case of the unknot, based
upon the presentation in \cite{BirFink}.  In
$\S$\ref{subsection:MTWS} we  state the MTWS, Theorem
\ref{theorem:MTWS}, in the general case.  A full proof of that
theorem can be found in \cite{BM-stab-I}.  In the final section,
$\S$\ref{section:braids and contact structures} we give an
application of the MTWS  to contact topology.

\subsection{The Markov Theorem Without Stabilization (special case: the unknot)}
\label{subsection:mtws}

\begin{figure}[htpb!]
\centerline{\includegraphics[scale=.65,  bb=61 497 496 645] {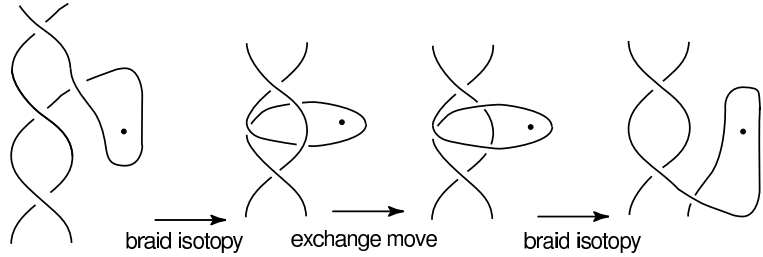}}
%\centerline{\includegraphics[scale=.65] {exchange2.eps}}
\caption{Exchange moves are the obstruction to sliding a trivial loop around a closed braid.}
\label{figure:exchange2}
\end{figure}

After describing the basic ideas about braid foliations,  we will apply them to the study of a classical problem in topology, the unknot recognition problem.
Alexander's Theorem (Theorem~\ref{theorem:Alexander's theorem})
tells us that every link $K$ may be represented as a closed
$n$-braid, for some $n$.   Markov's Theorem (Theorem
\ref{theorem:MT})  tells us how any two closed braid
representatives of the same knot or link are related.  Looking for
a way to simplify a given closed braid representative of a knot or
link systematically,  the first author and Menasco were lead to
the study of the unknot as a key example. There is an obvious
choice of a simplest representative, namely a 1-braid
representative. The Markov Theorem Without Stabilization for the
unknot, which is stated below as Theorem~\ref{theorem:mtws},
asserts that, in the special case of the unknot,  the
stabilization move of Markov's Theorem can be eliminated, at the
expense of adding the exchange move. Therefore we begin with a
discussion of the \underline{exchange move}, which is defined in
Figure \ref{figure:exchange2}, and the reason why it is so
important.

\begin{figure}[htpb!]
\centerline{\includegraphics[scale=.7,  bb=91 382 539 658]{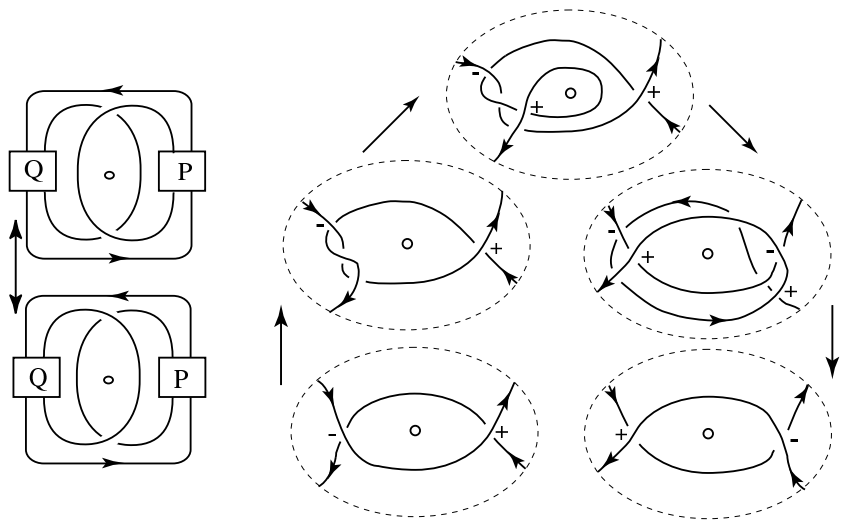}}
%\centerline{\includegraphics[scale=.7]{exchange1.eps}}
\caption{The left top and bottom sketches define the exchange move. The right sequence of 5 sketches shows how it  replaces a sequence of Markov moves which include braid isotopy, a single stabilization, additional braid isotopy and a single destabilization.}
\label{figure:exchange1}
\end{figure}
\begin{figure}[htpb!]
\centerline{\includegraphics[scale=.7, bb=49 454 530 720] {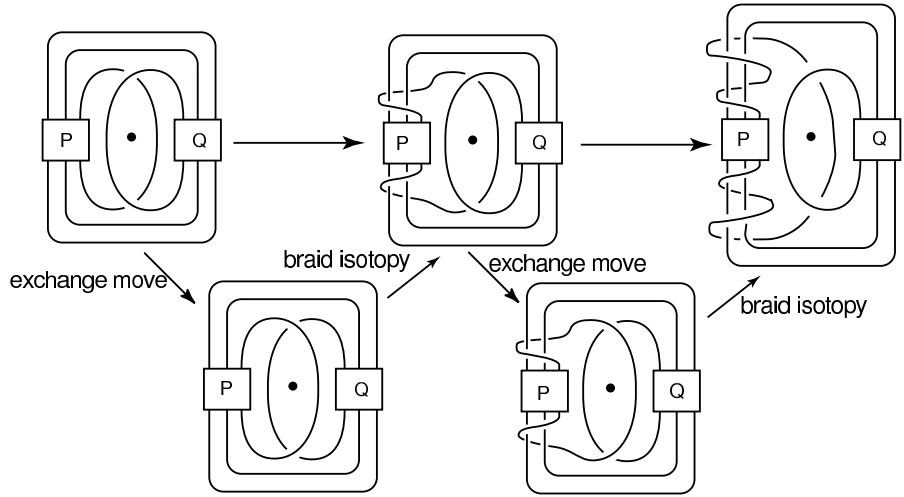}}
% \centerline{\includegraphics[scale=.7] {exchange3.eps}}
\caption{A sequence of exchange moves and  braid isotopies with non-trivial consequences.}
\label{figure:exchange3}
\end{figure}

A natural question to ask is what is accomplished by stabilization, braid isotopy and destabilization?
Figure \ref{figure:exchange2} shows that exchange moves are the
obstruction to sliding a trivial loop around a braid
\cite{Birman-Wrinkle}.   However there is more at issue than just
sliding trivial loops around a braid. As is illustrated in Figure
\ref{figure:exchange3} a sequence of exchange moves, together with
braid isotopy, can create infinitely many closed braid
representatives of a single knot type, all of the same braid
index.  This bypasses a key question: do exchange moves
actually change conjugacy class? The answer is ``yes", and in fact
the phenomenon first appears in the study of closed 4-braid
representatives of the unknot.  In $\S$\ref{section:the word and
conjugacy problems in the braid groups} we show that there is a
definitive test for proving it.

Keeping Figures \ref{figure:exchange1}, \ref{figure:exchange2} and
\ref{figure:exchange3} in mind,  we are ready to state the MTWS
for the unknot:
\begin{theorem}\label{theorem:mtws}{\rm \cite{BM-SLVCB-V}} Every closed braid representative $K$ of the
unknot ${\cal U}$ may be reduced to the standard
$1$-braid representative $U_1$, by a finite sequence of braid isotopies,  destabilizations and exchange moves.  Moreover there is a complexity function associated to closed braid representative in the sequence, such that each destabilization and exchange move is strictly complexity-reducing.
\end{theorem}

The first proof of Theorem \ref{theorem:mtws} was the one in  \cite{BM-SLVCB-V}.   A somewhat different and slicker proof can be found in \cite{BirFink}, but it requires more machinery than was necessary for present purposes.  We follow the proof in \cite{BM-SLVCB-V}.  However, parts of our presentation and most of our figures were essentially lifted (with the permission of both authors) from \cite{BirFink}, a review article on braid foliation techniques.  Our initial goal is to set up the machinery needed for the proof.

Recalling the definition of a closed braid from
$\S$\ref{subsection:closedbraids}, we adopt some additional
structure and say that $K$ is in closed $n$-braid form if there is
an unknot $A$ in $S^3 \setminus K$, and a choice of fibration $H$
of the solid torus $S^3 \setminus A$ by meridian disks, such that
$K$ intersects each fiber of $H$ transversely. Sometimes it is
convenient to replace $S^3$ by $R^3$ and to think of the fibration
$H$ as being by half-planes $\{ H_\theta;\theta \in [0, 2\pi] \}$
of constant polar angle $\theta$ through the $z$-axis. Note that
$K$ intersects each fiber $H_\theta$ in the same number of points,
that number being the \underline{index} $n$ of the closed braid
$K$. We may always assume that $K$ and $A$ can be oriented so that
$K$ travels around $A$ in the positive direction, using the right
hand rule. Now let $K$ be a closed braid representative of the
unknot, and let $\cD$ denote a disk spanned by $K$, oriented so
that the positive normal bundle to each component  has the
orientation induced by that on $K =
\partial \cD$.  In this section we will describe a set of ideas
which shows that there is a very simple method that changes the
pair $(K,\cD)$ to a planar circle that bounds a planar disc, via
closed braids, moreover there is an associated complexity function
that is strictly reducing. The ideas that we describe come from
\cite{BM-SLVCB-V}, however our main reference will be to the
review article \cite{BirFink}.

The braid axis $A$ and the fibers of $H$ will serve as a
coordinate system in 3-space in which to study $\cD$.   A singular
foliation of $\cD$ is induced by its intersection with fibers of
$H$.  A \underline{singular leaf} in the foliation is one which
contains a point of tangency with a fiber of $H$. All other leaves
are \underline{non-singular}.

It follows from standard general position arguments that the disk
$\cD$ can be chosen to be `nice' with respect to our fibration.
More precisely, we can assume the following:

\begin{enumerate}
\item[(i)] The intersections of $A$ and $\cD$ are finite in number and transverse.
\item[(ii)] There is a  neighborhood $N_A$ of $A$  in $R^3 \setminus K$ such that each component of $\cD \cap N_A$ is a disk, and each disk
is radially foliated by its arcs of intersection with fibers of $H$. There is also a
neighborhood $N_K$ of $K$ in $R^3$ such that $N_K \cap \cD$ is foliated by
arcs of intersection with fibers of $H$ which are transverse to $K$.
\item[(iii)] All but finitely many fibers $H_\theta$ of $H$ meet $\cD$
transversely, and those which do not (the \underline{singular
fibers}) are each tangent to $\cD$ at exactly one point in the
interior of both $\cD$  and $H_\theta$.  Moreover, each point of
tangency is a saddle point (with respect to the parameter
$\theta$). Finally, each singular fiber contains exactly one
singularity of the foliation, each of which is a saddle point.
\end{enumerate}
It can also be assumed (see \cite{BirFink} for details) that
\be
\item[(iv)] Each non-singular leaf is either an \underline{{\it a}-arc}, which has one endpoint on $A$ and
one on $K = \partial \cD$ or a \underline{{\it b}-arc}, which has both endpoints on $A$.

\item[(vii)]  Each $b$-arc  in a fiber $H_\theta$ separates that fiber into two components. Call
 the $b$-arc  \underline{essential} if each of these components is pierced at least  once by $K$, and
\underline{inessential} otherwise.  Then we also have that all $b$-arcs in the foliation of $\cD$ may be
assumed to be essential.
\ee

\begin{figure}[htpb!]
\centerline{\includegraphics[scale=.8,  bb=126 244 548 537] {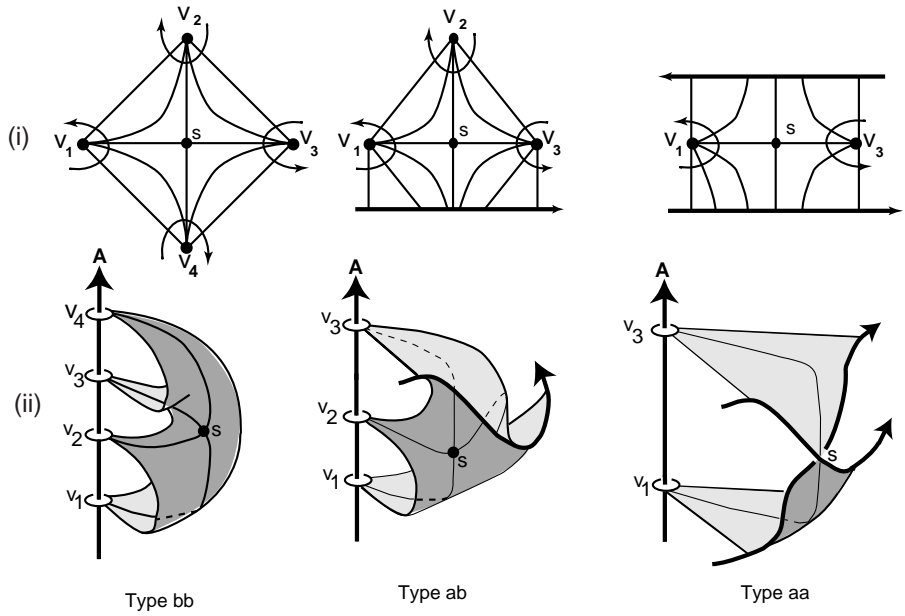}}
% \centerline{\includegraphics[scale=.8] {tiles.eps}}
\caption{(i) The three types of tiles in the decomposition of $\cD$; (ii)The canonical embedding of each type of tile.} \label{tiles}
\end{figure}

Note that $\cD$ cannot be foliated entirely by $b$-arcs, since such a surface is necessarily a $2$-sphere.
Thus, if the foliation of $\cD$ contains no singularities, $\cD$ is foliated entirely by $a$-arcs. Otherwise,
let $\it U$ be the union of all the singular leaves in the foliation of $\cD$.  Moving forward through the
fibration, we see that any singular leaf in the foliation is formed by non-singular leaves  moving together to
touch at a saddle singularity.  The three types of singular leaves which can occur are labelled $aa$, $ab$ or
$bb$, corresponding to the non-singular leaves associated to them.  Now each singular leaf $\lambda$ in $U$
has a foliated neighborhood  $N_\lambda$ in $\cD$ such that $N_\lambda \cap U = \lambda$.    According to
whether $\lambda$ has type $aa$, $ab$, or $bb$, $N_\lambda$ is one of the foliated open $2$-cells shown in
Figure~\ref{tiles}(i), with the arrows indicating the direction of increasing $\theta$.

The complement of $U$ in $\cD$ is a union $B_1 \cup B_2
\cup...\cup B_k$, where each $B_i$ is  foliated entirely by
$a$-arcs or entirely by $b$-arcs.  Choose one non-singular leaf in
each $B_i$ and declare it to be a \underline{boundary arc} of
\underline{type $a$} or \underline{type $b$} according to whether it
is an $a$-arc or $b$-arc, respectively.  Then the union of all
boundary arcs determines a \underline{tiling} of $\cD$, that is, a
decomposition into regions called \underline{tiles}, each of which
is a  foliated neighborhood of one singular leaf.  Each tile has
type $aa$, $ab$ and $bb$, according to the type of its unique
singularity. Note that a tiling of $\cD$ is a foliated
cell-decomposition. We further define the \underline{sign of a
singular point} $s$ in the foliation of $\cD$ to be positive if
the positive normal to $\cD$ points in the direction of increasing
$\theta$ in the fibration, negative otherwise. The \underline{sign
of a tile} is then defined to be the sign of its singularity.

\begin{figure}[htpb!]
\centerline{\includegraphics[scale=.9, bb=71 327 523 469] {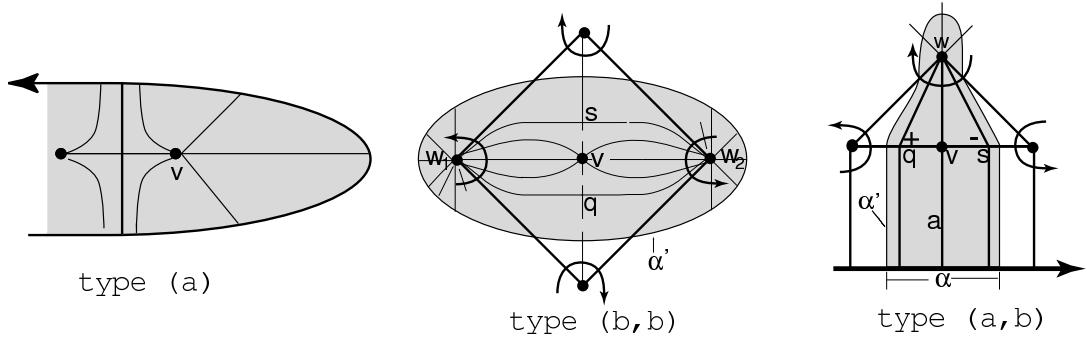}}
%\centerline{\includegraphics[scale=.9] {vertices.eps}}
\caption{A vertex $v$ of type $(a), (ab)$ and $(bb)$.}\label{figure:vertex type}
\end{figure}

The axis $A$ intersects the surface in a finite number of points, called \underline{vertices} of the tiles.
Each vertex $v$ is an endpoint of finitely many boundary arcs in the surface decomposition. Let the
\underline{type} of $v$ be the cyclic sequence $(x_1,...,x_r)$, where  each $x_i$ is either $a$ or $b$, and the  sequence lists the  types of boundary arcs meeting at $v$ in the cyclic order in which they occur in the  fibration.  See Figure \ref{figure:vertex type}.
 The \underline{valence of a vertex} $v$ is the number of distinct  tiles intersecting  at $v$. The
\underline{sign of a vertex} $v$ is the cyclic array of signs of
the tiles meeting at $v$.  See Figure \ref{figure:vertex type}.

The \underline{parity} of $v$ is said to be positive (resp.
negative) if the outward-drawn normal to the surface has the same
(resp. opposite) orientation as the braid axis at the vertex.  Thus
when we view  the positive side of the surface, the sense of
increasing $\theta$ around a vertex will be counterclockwise
(resp. clockwise) when the vertex is positive (resp. negative), as
illustrated in Figure~\ref{tiles}(i).

Finally, we note that (again, see \cite{BirFink} for a proof)
that, up to a choice of their sign, tiles of type $aa,ab,bb$ each
have a canonical embedding in 3-space, which is determined up to
an isotopy of 3-space which preserves the axis $A$ and each fiber
of $H$ setwise. The canonical embedding for each tile is shown in
Figure~\ref{tiles}(ii).

The decompositions of the disk $\cD$ which we have just described
are not unique.  We shall now describe three ways in which they
can be changed. In each of the three cases the possibility of
making the change is indicated by examining   the combinatorics of
the tiling. The change is  realized by an isotopy of the disk
which is supported  in a neighborhood  $N$ of a specified small
number of tiles of type $aa,ab$ or $bb$, leaving the decomposition
of $\cD$ unchanged outside $N$.

For a spanning disk $\cD$ and a tiling $\cT$, we denote the \underline{complexity} of the
tiling by $c(\cD, \cT)$, where we define $c(\cD, \cT) = (n, S)$, where $n$ is the braid index of the boundary and $S$ is the number of singularities in the tiling.  Setting $V_+$ and $V_-$ equal to the number of positive and negative vertices, one has $n = V_+ - V_-$, so $n$ is also determined by the tiling.
Tiled discs $(\cD, \cT)$ can then be ordered, using lexicographical ordering of  the associated pair $c(\cD, \cT)$.  Each of our three moves will replace the given disk and tiling $\cD,\cT$ with some $\cD',\cT '$, with the following effects on complexity:

\begin{itemize}
\item {\bf Change in foliation:} $(n',S') = (n,S)$.
\item {\bf Destabilization:} $(n',S') = (n - 1,S - 1)$.
\item {\bf Exchange moves (two types):} $(n',S') = (n,S - 2)$.
\end{itemize}

We focus here only on the combinatorics of tilings which admit one
of the above moves, as well as the effect of each move on the
embedding of $\cD$ and the new tiling of $\cD$ which results. For
proof and  further details, see \cite{BirFink}.

\

\noindent $\bullet${\bf Changes in foliation.}   The choice of a
foliation  of $\cD$  is not unique, and our first move involves
ways in which the surface decomposition can be changed by an
isotopy of $\cD$ or, equivalently, by an isotopy of the fibers of
$H$ keeping $\cD$ fixed. This particular change was introduced in
\cite{BM-SLVCB-IV} for $2$-spheres and was modified in
\cite{BM-SLVCB-V} for certain spanning surfaces.   In what
follows, we say that two tiles are \underline{adjacent} if they
have a common boundary arc in the given tiling.
\begin{figure}[htpb!]
\centerline{\includegraphics[scale=.8, bb=80 392 480 619] {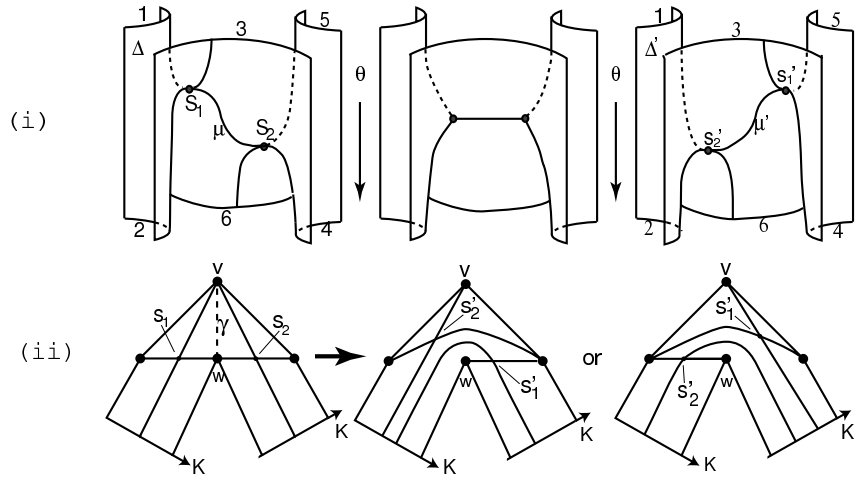}}
%\centerline{\includegraphics[scale=.8] {ch_fol.eps}}
\caption{Sketch (i) shows an isotopy of $\cD$ that induces a
change in foliation.  Sketch (ii) shows the effect of the change
on two $ab$-tiles} \label{ch_fol}
\end{figure}

 A tiling admits a change in foliation whenever there are two tiles $T_1$, $T_2$ of the same sign
adjacent at a $b$-arc.  Roughly speaking, a change in foliation is
a local isotopy of the surface which pushes two saddle points past
each other.  Locally, the disk $\cD$ is embedded as in the left
sketch in Figure~\ref{ch_fol}(i). A \underline{change in
foliation} is defined as the passage from the embedding left to
the right embedding. Figure~\ref{ch_fol}(ii)  the effect of this
move on the local foliation, while Figure~\ref{ch_fol}(iii)
illustrates the effect of a change in foliation on the tiling, in
the case of two adjacent $(ab)$-tiles.

\

\begin{figure}[htpb!]
\centerline{\includegraphics[scale=.7, bb=42 325 578 484]{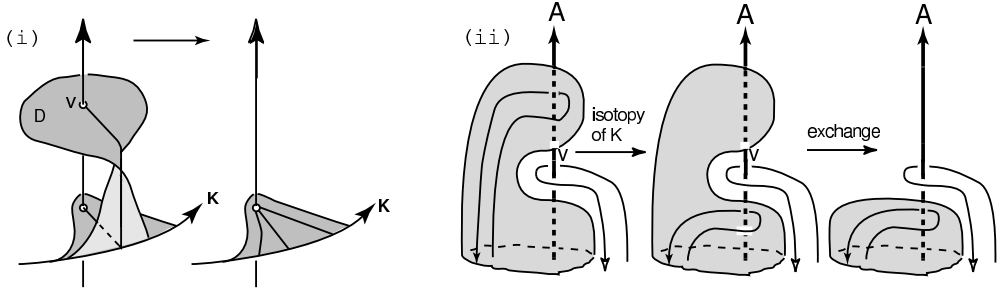}}
%\centerline{\includegraphics[scale=.7] {ch_disc.eps}}
\caption{The effect of (i) a destabilization, and (ii) an exchange move of type (bb), on the embedded
(and foliated) foliated disc
$\cD$ and on the braid $K$ that is its boundary.}\label{ch_disc}
\end{figure}

\noindent $\bullet${\bf Destabilization via a type $(a)$ vertex.}
A vertex $v$ of  valence 1 in a tiling of $\cD$ occurs when two of the edges in a single tile $T$  are
identified in $\cD$. Such a vertex must have type $(a)$, and therefore $T$ is an $aa$-tile. See Figure
\ref{figure:vertex type}.
 Since there is a canonical embedding for an $aa$-tile in 3-space (Figure~\ref{tiles}),
identifying two edges of an $aa$-tile with endpoints on a common
vertex yields the canonical embedding for $T$ shown in
Figure~\ref{ch_disc}(i). Notice that there is a radially foliated
disk $D$ in $T$ cut off from $\cD$ by the arc of the singular leaf
of $T$ with both endpoints on $K$.  Hence there must be a trivial
loop in the braid representation of $K$, and we can modify the
braid and the disk $\cD$ in the manner illustrated. We call this
modification \underline{destabilization} via a type $(a)$ vertex.
The effect on the tiling of $\cD$ is that the $aa$-tile $T$ and
its type $(a)$ vertex $v$ are deleted, while the tiling outside
$T$ is unaltered.  The braid index is decreased by one. See the
right sketch in Figure \ref{ch_disc}(i).

\

\noindent$\bullet${\bf Exchange moves.}  The destabilization move that we just described can be
accomplished if a  vertex of valence 1 exists in our tiling.  Exchange moves are based upon the existence of a
vertex of valence 2.  An example of the type of move we would like to achieve on the braid itself was shown in sketch (ii) of Figure \ref{ch_disc}(ii), and also in in the second row from the top in Figure \ref{figure:templates}. There are two different types of tiling patterns, and hence two
different types of local embeddings of $\cD$, which support exchange moves:

\

\noindent {\bf Exchange Move Type (bb):} We now consider a vertex $v$ with sign $(\pm,\mp)$  and
type $(b,b)$.  In this situation, two $bb$-tiles are adjacent along consecutive
$b$-arcs.  The canonical embeddings of
the tiles in this case are shown in Figure~\ref{tiles}, and the overall modifications to $\cD$ are
suggested by Figure~\ref{ch_disc}(ii). The end result is a new surface $\cD'$ whose decomposition has
at least two fewer vertices (in particular, the vertex $v$ is deleted) and at least two fewer regions than the
decomposition of $\cD$. $\|$

\

\noindent{\bf Exchange Move, Type (ab)}  An exchange move of type (ab) is possible whenever the tiling of
$\cD$ has a vertex $v$ of valence 2, type (a,b), and sign $(\pm, \mp)$.  Such a vertex occurs only when two
$ab$-tiles are adjacent along corresponding $a$ and $b$ edges, as in Figure \ref{figure:vertex type}.  The isotopy of
$\partial
\cD = K$ is achieved by pushing a subarc $\alpha$ of $K$ across a disk which is contained in the union of two
$(ab)$ tiles
$T_1, T_2$, to a new arc $\alpha'$.  The effect on the tiling is
that the
$ab$-tiles $T_1$ and $T_2$ are deleted and any adjacent $ab$-tiles (resp. $bb$-tiles) become $aa$-tiles (resp.
$ab$-tiles).

\

In order to prove the MTWS for the unknot, we must locate the
places on the tiling of $\cD$ where the three moves just described
can be made.  We begin with two lemmas.  The first of the two is
implicit in \cite{Bennequin}, but we refer the reader to
\cite{BirFink} for an explicit proof.

\begin{lemma}\label{star}{\rm \cite{Bennequin}}
Let $v$ be a vertex of type $(b,b, \ldots, b)$. Then the set of all singularities $s$ which lie on a
singular leaf ending at $v$ contains both positive and negative singular points.
\end{lemma}

\begin{lemma}{\rm \cite{BM-SLVCB-V}}\label{euler}  Suppose that the disk $\cD$ is nontrivially
tiled and has no vertices of valence $1$.  Then the tiling of $\cD$ contains a vertex of type (ab), (bb), or
(bbb).
\end{lemma}

\noindent{\bf Proof.}  The proof is essentially an Euler
characteristic argument.  A tiling $\cT$ of the disk $\cD$
corresponds to a cell decomposition of the $2$-sphere in the
following way.  Let $V$ be the number of vertices in $\cT$, let
$E$ be the number of boundary arcs in $\cT$, which we will think
of as \underline{edges}, and let $F$ be the number of tiles in $\cT$.
Now let $\Sigma$ be the $2$-sphere obtained by collapsing
$\partial \cD = K$ to a point. This gives a cell decomposition of
$\Sigma$ with $V+1$ 0-cells (vertices), $E$ 1-cells (edges), and
$F$ 2-cells (faces).  We know that the Euler characteristic $\chi
(\Sigma) = 2$, and therefore we have that $V - E + F = 1.$  We
also know that every face in our cell decomposition of $\Sigma$
has 4 edges, and that every edge has 2 adjacent faces, and so $E =
2F$.  We therefore have the following equation.
\begin{equation}\label{euler1}
2V - E = 2.
\end{equation}

Now let $V(\alpha, \beta)$ denote the number of vertices in $\cT$ with $\alpha$ adjacent $a$-arcs and
$\beta$ adjacent $b$-arcs.  If $v$ is the valence of such a vertex, then $v = \alpha + \beta$, and so
$V(\alpha, v-\alpha)$ denotes the number of vertices in $\cT$ with valence $v$ and $\alpha$ adjacent
$a$-arcs.  Since by hypothesis $v \neq 1$, we can now count the number of vertices in $\cT$ using the
following summation:
$$
V = \sum_{v = 2}^\infty \sum_{\alpha = 0}^v V(\alpha, v - \alpha)
$$
Now write $E = E_a + E_b$, where $E_a$ denotes the number of edges in $\cT$ of type a, and similarly
for $E_b$.  Since an $a$-edge is incident at one vertex in $\cT$, and a $b$-edge is incident at two vertices
in $\cT$, we can write:
\begin{eqnarray}
E_a &=& \sum_{v = 2}^\infty \sum_{\alpha = 0}^v \alpha V(\alpha, v - \alpha) \\
E_b &=& \frac{1}{2}\sum_{v = 2}^\infty \sum_{\alpha = 0}^v (v - \alpha) V(\alpha, v - \alpha)
\end{eqnarray}
Substituting into Equation~\ref{euler1}, we now have
$$
\sum_{v = 2}^\infty \sum_{\alpha = 0}^v (4 - v - \alpha) V(\alpha, v - \alpha)=4.
$$

Observe now that if $v \geq 4$, the coefficient $(4 - v - \alpha) \leq 0$.  Writing out the terms
corresponding to $v = 2$ and $v = 3$ gives, respectively, $2 V(0,2) + V(1,1)$ and $V(0,3) - V(2,1) - 2
V(3,0)$.  Each $V(\alpha, \beta) \geq 0$ by definition.  Hence we can move terms around in the above equation
so that each term is nonnegative, as follows:
$$
2 V(0,2) + V(1,1) + V(0,3) = 4 + V(2,1) + 2 V(3,0) + \sum_{v = 4}^\infty \sum_{\alpha = 0}^v (4 - v -
\alpha) V(\alpha, v - \alpha).
$$
The right-hand side of the equation is clearly $\geq 4$, and therefore at least one term on the
left-hand side is nonzero.  Since $V(0,2), V(1,1)$, and $V(0,3)$ record the number of type $(bb)$, type $(ab)$, and
type $(bbb)$ vertices, respectively, the lemma is proved.   \endpf

We are now ready to prove Theorem \ref{theorem:mtws}, the MTWS in
the case of the unknot.

\

\noindent  {\bf Proof of Theorem 4.} We start with an arbitrary
closed braid representative $K$ of the unknot. Our closed braid
$K$ is the boundary of a disk $\cD$ which admits a tiling $\cT$.
Note that if the complexity of the tiling $c(\cD, \cT)$ is equal
to $(1,0)$, then our disk is radially foliated by $a$-arcs and $K$
is the standard embedding of the unknot.  Therefore we assume
$c(\cD, \cT) > (1,0)$.  We will use induction on $c(\cD, \cT)$ to
show that after a finite sequence of our three moves, $\cD$ is
radially foliated by $a$-arcs.

We first observe that if there exists a vertex of type (a) in
$\cT$, then we can destabilize along this vertex, thereby reducing
$c(\cD, \cT)$.  Hence we can asume that all vertices in $\cT$ have
valence at least $2$.  By Lemma~\ref{euler}, $\cD$ must have a
vertex of type $(ab), (bb)$, or $(bbb)$.  If there is a type $(bbb)$
vertex, then at least two of the adjacent tiles must have the same
sign.  In this case we can do a change in foliation, replacing our
$(bbb)$ vertex with a $(bb)$ vertex.  Therefore we may assume that
$\cD$ contains a vertex of type $(ab)$ or type $(bb)$.

If $\cD$ has a type $(ab)$ vertex with sign (+,+) or (-,-), then we
can apply a change in foliation which replaces this vertex with a
type $(a)$ vertex.  We now destabilize along this type $(a)$ vertex
in order to reduce the complexity.  If $\cD$ has a type $(ab)$ vertex
with sign (+, -), then we can do an exchange move of type $(ab)$,
thereby reducing complexity.  If $\cD$ has a type $(bb)$ vertex,
then it is an interior vertex, and by Lemma~\ref{star}, it has
sign (+,-).  Then we can do an exchange move of type $(bb)$, which
reduces complexity.

Thus if $c(\cD, \cT) > (1,0)$, we can always reduce until $c(\cD, \cT) = (1,0)$, and the theorem is proved.
$\|$

\begin{example}{\bf [Morton's irreducible 4-braid]}  {\rm We now give an example
illustrating the necessity of exchange moves.  In
\cite{Morton1983}, Morton gave the following example of an
\underline{irreducible} braid:
$$X  = \sigma_3^{-2} \sigma_2 \sigma_3^{-1} \sigma_2
\sigma_1^3 \sigma_2^{-1} \sigma_1 \sigma_2^{-1},$$ i.e., a braid
whose closure represents the unknot but cannot be isotoped to the
unknot in the complement of the braid axis. Now, it is clear that
if a braid in ${\bf B}_n$ admits a factorization of the form $w_1
\sigma_{n-1} w_2 \sigma_{n-1}^{-1}$, where $w_1$ and $w_2$ are
words in $\sigma_1, \ldots, \sigma_{n-2}$, then its closure admits
an exchange move. As observed in \cite{BM-SLVCB-V}, the conjugate
braid $(\sigma_3 \sigma_2 \sigma_1) X (\sigma_1^{-1}
\sigma_2^{-1} \sigma_3^{-1})$ is isotopic to the braid $$X' =
(\sigma_2^{-2} \sigma_1 \sigma_2^{-1}) \sigma_3 (\sigma_2^3
\sigma_1^{-1} \sigma_2) \sigma_3^{-1}.$$  Thus $X'$ admits an
exchange move to obtain the new braid word $(\sigma_2^{-2}
\sigma_1 \sigma_2^{-1}) \sigma_3^{-1} (\sigma_2^3 \sigma_1^{-1}
\sigma_2) \sigma_3$, which can be cyclically rewritten as
$$
X'' = (\sigma_2 \sigma_3 \sigma_2^{-2}) \sigma_1
(\sigma_2^{-1} \sigma_3^{-1} \sigma_2^3) \sigma_1^{-1}.
$$
If we now interchange the axis of our braid with the point at
infinity, we can perform another exchange move to obtain
$$
X''' = (\sigma_2 \sigma_3 \sigma_2^{-2}) \sigma_1 ^{-1}
(\sigma_2^{-1} \sigma_3^{-1} \sigma_2^3) \sigma_1
$$
which can in fact be isotoped to the unknot in the complement of
its braid axis \cite{BM-SLVCB-V}.

Wright has given a foliated disk corresponding to this braid in
\cite{Wright2000}.  The reader is cautioned that Wright uses a
different notation convention for braid words than that used in
this paper (see p. 98 of \cite{Wright2000}).  We also refer the
reader to \cite{BRBV2002} for further examples.} $\spadesuit$
\end{example}

We note that Theorem \ref{theorem:mtws} has been used as the basis for an algorithm for recognizing the unknot. See \cite{BH1998}.  This algorithm has been put on a computer \cite{BRBV2002}, however it needs more work before it can become a practical tool for recognizing the unknot.

\begin{openproblem}
\label{openproblem:recognize unknot} At this writing the
development of a sound and practical algorithm for unknot
recognition remains one of the major open problems in low
dimensional topology.  $\clubsuit$
\end{openproblem}

With regard to Problem \ref{openproblem:recognize unknot}, we note that Theorem \ref{theorem:mtws} proves the existence of a monotonic simplification process which begins with an arbitrary closed n-braid representative of the unknot and ends with a 1-braid representative, but the complexity function $c(\cD, \cT)$ which gives instructions for the process is `hidden' in the tiling of the surface $\cD$.  We need it in order to know  when complexity-reducing destabilizations and exchange moves are possible. The reader who is interested in this problem might wish to consult \cite{BH1998}, where a somewhat different   approach suggests itself.  Instead of working with the tiled surface, one may work with the `extended boundary word' of \cite{BH1998}.  The latter is a closed braid which is obtained from an arbitrary closed braid representative of the given knot  by threading in additional 1-braids, and it seems likely that it will give an alternative monotonic reduction process.  Thus we suggest:

\begin{openproblem}
\label{openproblem:extended boundary word}
Investigate the monotonic reduction process of Theorem \ref{theorem:mtws} from the point of view of the
extended boundary word of {\rm \cite{BH1998}}. $\clubsuit$
\end{openproblem}

\subsection{The Markov Theorem Without Stabilization, general case}
\label{subsection:MTWS}

\begin{figure}[htpb!]
\centerline{\includegraphics[scale=.7, bb=66 63 537 687] {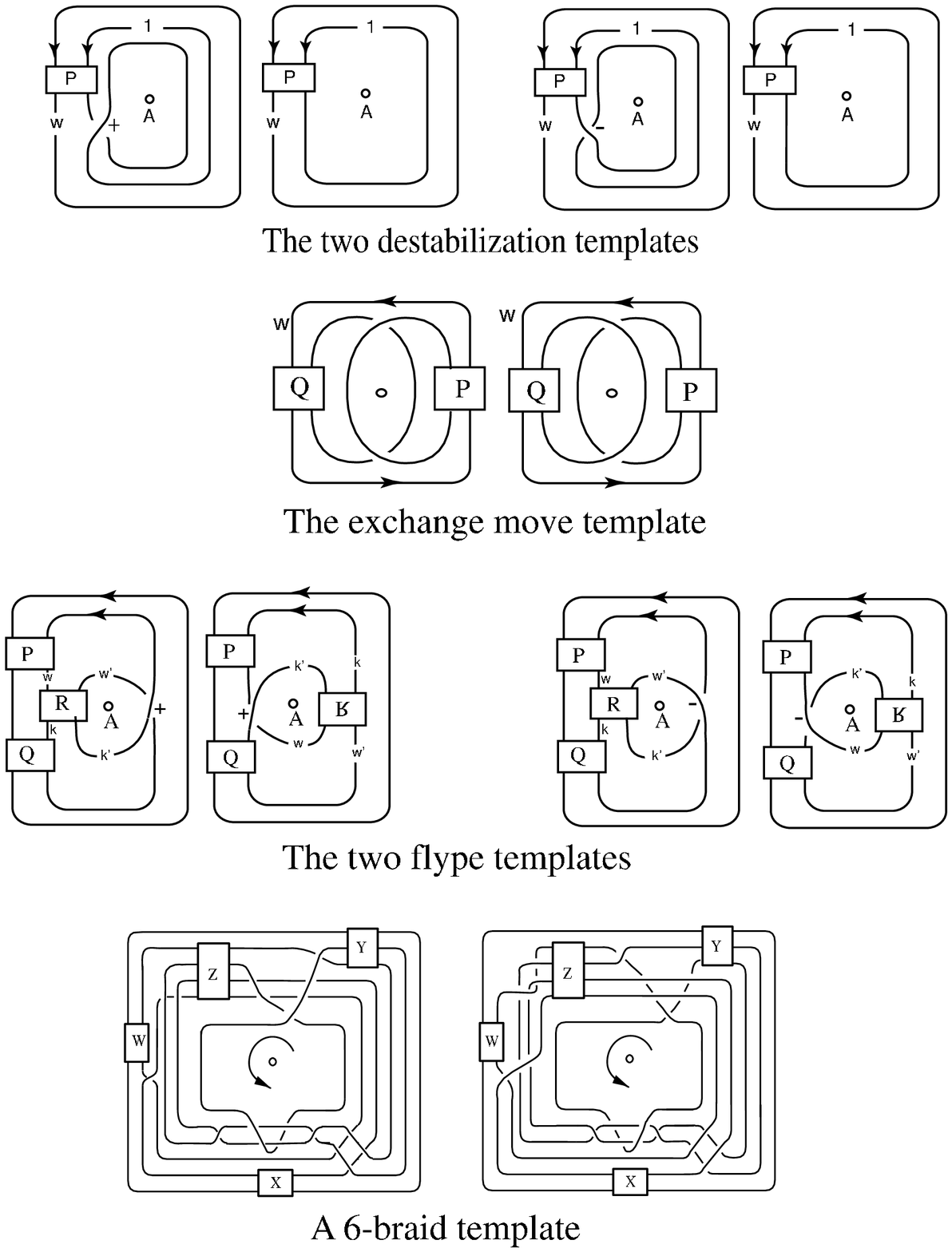}}
%\centerline{\includegraphics[scale=.7] {templates.eps}}
\caption{Examples of templates} \label{figure:templates}
\end{figure}

In this section we state the generalized version of Theorem \ref{theorem:mtws} which was established in \cite{BM-stab-I} for arbitrary closed braid representatives of arbitrary knots and links.   The moves which are needed are of course much more complicated than in the case of the unknot. They are described in terms of `block strand diagrams' and  `templates'.   The concept of a block-strand diagram is fairly easy to understand from the 8 examples in Figure \ref{figure:templates}.
The first important feature of a block-strand diagram is that
after an assignment of a braided tangle to each block, it becomes
a closed braid that represents a specific knot or link.  Thus each
block strand diagram determines infinitely many closed braid
representatives of (presumably) infinitely many knots and links.
A \underline{template} is a pair of block-strand diagrams, both of
which represent the same knot or link, when we make the same
braiding assignments to corresponding blocks.  Note that there are two destabilization templates (they differ in the sign of the `trivial loop' that is removed). There are also two flype templates, which differ in the sign of the single crossing which is outside all the braid boxes.  Templates always have an associated braid index, namely the braid index of any knot or link of minimum braid index that they carry. Destabilization templates occur for every braid index $\geq 2$. Since the exchange move can be realized by braid isotopy when the braid index is $\leq 3$, it does not play a role  until braid index $\geq 4$.  The flype templates occur for braid index $\geq 3$. The method by which
the 6-braid template in Figure \ref{figure:templates} was
constructed is described in the manuscript \cite{BM-stab-I}.

The
two block strand diagrams in a template are always related by a
sequence of Markov moves,  however  the sequence may be quite
complicated, and so the isotopy that takes the left diagram  to
the right diagram is in general not obvious.  This fact is illustrated by the 6-braid template of Figure \ref{figure:templates}.  (See \cite{BM-stab-I} for the isotopy, as an explicit sequence of Markov moves.)

With this brief introduction, we are able to state the Markov
Theorem Without Stabilization in the general case:

\begin{theorem}{\rm \cite{BM-stab-I}}
\label{theorem:MTWS}
Let $\cB$ be the collection of all braid isotopy classes of closed braid representatives of oriented
knot and link types in oriented 3-space. Among these, consider the subcollection
$\cB(K)$ of representatives of a fixed link type $K$.  Among these, let $\cB_{min}(K)$ be the
subcollection of representatives whose braid index is equal to the braid index of $K$.  Choose any
$X_+\in \cB(K)$  and any $X_-\in \cB_{min}(K)$.  Then there is a
complexity function  which is associated
to $X_+,X_-$, and for each braid index $m$ a finite set $\cT(m)$ of templates is introduced, each
template determining a move which is non-increasing on braid index, such that the following hold:
First, there is are initial sequences which modify $X_- \to X_-'$ and $X_+ \to X_+'$:
$$X_- = X_-^1 \to \cdots\to X_-^p = X_-',   \ \ \ \ \ X_+ = X_+^1\to\dots \to X_+^q = X_+'  $$
Each passage $X_-^j \to X_- ^{j+1}$  is strictly complexity
reducing  and is realized by an exchange move, so that $b(X_-^{j+1}) = b(X_-^j)$.  These moves `unwind' $X_-$, if it is wound up as in the top right sketch in Figure \ref{figure:exchange3}.  Each passage
$X_+^j \to X_+^{j+1}$   is strictly  complexity-reducing and is realized by either an exchange move or a destabilization, so that $b(X_+^{j+1})\leq b(X_+^j)$.   Replacing $X_+$ with $X_+'$ and $X_-$ with $X_-'$, there is an additional
sequence which modifies $X_+'$, keeping $X_-'$ fixed:
$$ X_+'=X^q\to\cdots \to X^r = X_-' $$
Each passage $X^j\to X^{j+1}$ in this sequence  is also
strictly complexity-reducing.  It is realized by one of the moves defined by a template
$\cT$ in the finite set $\cT(m)$, where $m=b(X_+)$. The inequality $b(X^{j+1})\leq b(X^j)$ holds
for each $j = q,\dots,r-1$ and so also for each $j=1,\dots,r-1$.
\end{theorem}

The proof of Theorem \ref{theorem:MTWS} uses the braid foliation techniques that were used in the proof of Theorem \ref{theorem:mtws}, but in a more complicated setting.  Instead of looking at a foliated embedded disc which is bounded by a given unknotted closed braid, we are given two closed braids, $X_+$ and $X_-$, and an isotopy that takes $X_+$ to $X_-$.  The trace of the isotopy sweeps out an annulus, but in general it is not embedded.  The proof begins by showing that the given isotopy can be split into two parts, over which we have some control. An intermediate link $X_0$ which represents the same link type $K$ as $X_+$ and $X_-$  is constructed, such that the trace of the isotopy from $X_+$ to $X_0$ is an embedded annulus $\cA_+$.  Also the trace of the isotopy from $X_0$ to $X_-$ is a second embedded annulus $\cA_-$. The union of these two embedded annuli $\cT\cA = \cA_+\cup\cA_-$  is an immersed annulus, but its self-intersection set is controlled, and is a finite number of clasp arcs.  The main tool in the proof of Theorem \ref{theorem:MTWS} is the study of the braid foliation of the  immersed annulus $\cT\cA$.   In the next section,  i.e.  $\S$\ref{section:braids and contact structures}, we will see how Theorem \ref{theorem:MTWS} was used to settle a long-standing problem about contact structures.

 \subsection{Braids and contact structures}
 \label{section:braids and contact structures}

In this section we describe how Theorem \ref{theorem:MTWS}  was
used in \cite{BM-stab-II}  to settle a problem  about contact
structures on $\reals^3$ and $S^3$.

\begin{figure}[htpb!]
\centerline{\includegraphics[scale=.7, bb=53 489 558 698] {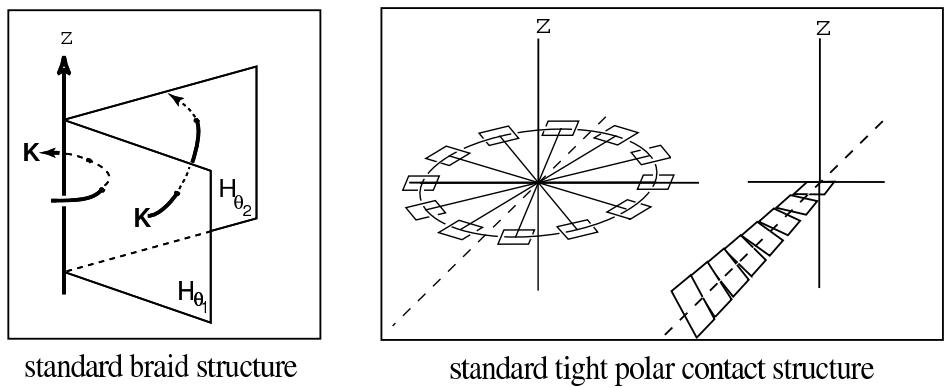}}
%\centerline{\includegraphics[scale=.7] {braid-contact.eps}}
\caption{ The standard braid structure and tight contact structure on $\reals^3$} \label{figure:braid-contact}
\end{figure}

Let $\axis$ be the $z$-axis in $\reals^3$, with standard cylindrical coordinates $(\rho,\theta, z)$ and let ${\bf H}$ be the collection of all half-planes $H_\theta$ through $\axis$. The pair $(\axis,{\bf H})$ defines the  \underline{standard braid} \underline{structure} on $\reals^3$.   Using the same cylindrical coordinates, let $\alpha$ be the 1-form $\alpha = \rho^2 d\theta + dz$.  The kernel $\xi$  of the 1-form $\alpha$ defines  a \underline{contact structure} on $\reals^3$.  We can visualize $\xi$ by imagining that there is  a 2-plane (spanned by $\partial/\partial x$ and $\partial/\partial \rho$) attached to every point in $\reals^3$. Figure \ref{figure:braid-contact} shows both the braid structure and the polar contact structure,  for comparison.  The family of 2-planes that define $\xi$ twist (to the left) as one moves along the $x$-axis  from $0$ to $\infty$.  The family is invariant under rotation of 3-space about the $z$-axis and under translation of 3-space along rays parallel to the $z$-axis.  Its salient feature is that it is totally non-integrable, that is there is no surface in $\reals^3$ which is everywhere tangent to the 2-planes of $(\xi)$ in any neighborhood of any point in $\reals^3$. (Of course this makes it hard to visualize). The twisting is generic in the sense that, if $p$ is a point in a contact 3-manifold $M^3$, then in every neighborhood of $p$ in $M^3$ the contact structure is locally like the one we depicted in Figure \ref{figure:braid-contact}.

Let $K$ be a knot (for simplicity we restrict to knots here, but
everything works equally well for links)  which is parametrized by
cylindrical coordinates $(\rho(t), \theta(t), z(t))$, where $t \in
[0,2\pi]$.    Then, as defined in
$\S$\ref{subsection:closedbraids}, $K$ is a closed braid if
$\rho(t)>0$ and $d\theta/dt > 0$ for all  $t$.  On the other hand,
$K$ is a \underline{Legendrian} (resp. \underline{transversal})
knot if it is everywhere (resp. nowhere) tangent to the 2-planes
of $\xi$.  In the Legendrian case this means that on $K$ we have
$d\theta/dt = (-1/\rho^2)(dz/dt)$.  In the transversal case we
require that $d\theta/dt > (-1/\rho^2)(dz/dt)$ at every point of
$K(t)$. We are interested here primarily in the transversal case.

The \underline{total twist} of the contact structure is the number of multiples of $\pi/2$ as one traverses the positive real axis from the origin to $\infty$.  The case when the total twist angle is $\pi$  is known as the \underline{standard} (polar) contact structure.  We
call it $\xi_\pi$.  While the braid structure is very different
from the standard polar contact structure near $\axis$, for large
values of $\rho$ the 2-planes in $\xi_\pi$ are very close to the
half-planes $H_\theta$ of the braid structure.  In fact, assume
that $X$ is a closed braid which represents a knot $K$, and that
$X$ bounds a Seifert surface of minimum genus which supports a
foliation as in $\S$\ref{subsection:mtws}.   Assume further that
the closed braid $X$ intersects every 2-plane in the contact
structure transversally. Then, in the complement of a tubular
neighborhood of the braid axis, from the point of view of a
topologist, the braid foliation and the foliation induced by the
contact structure will be  `the same'.

Recall that Theorem \ref {theorem:Alexander's theorem} (Alexander's Theorem) was first proved in 1925.  Sixty years later, Bennequin adapted Alexander's original proof (which is different from the proof given in this article) to the setting of transversal knots in \cite{Bennequin}, where he showed that  every transversal knot is isotopic, through transversal knots, to a closed braid.
In 2002  Orevkov and Shevchishin extended Bennequin's ideas and proved a version of Theorem \ref{theorem:MT} (Markov's Theorem)  which holds in the transversal setting:

\begin{theorem}
\label{theorem:TMT}
{\rm \cite{O-S2003}:} Let $TX_+,TX_-$ be closed
braid representatives of the same oriented link type $K$ in oriented
3-space.  Then there exists a sequence of closed braid representatives
of $\cT\cK$:
$$ TX_+ = TX_1 \to TX_2 \to \cdots \to TX_r = TX_1 $$
such that, up to braid isotopy, each $TX_{i+1}$ is obtained from $TX_i$ by a single positive stabilization or destabilization.
\end{theorem}
Is the Transverse Markov Theorem  really different from the Markov
Theorem?  Are there transversal knots which are isotopic as
topological knots but are not transversally isotopic?  To answer
this question we take a small detour and review the contributions
of Bennequin in \cite{Bennequin}.

Why did topologists begin to think about contact structures, and
analysts begin to think about knots?  While we might wish that
analysts suddenly became overwhelmed with the beauty of knots,
there was a more specific and focused reason.   At the time that
Bennequin did his foundational work in \cite{Bennequin} it was not
known  whether a 3-manifold could support more than one isotopy
class of contact structures.  Bennequin answered this question in
the affirmative, in the case of contact structures on $\reals^3$
or $S^3$ which were known to be homotopic to the standard one.
His tool for answering it was highly original, and it had to do
with braids and knots.

Let $TK$ be a transversal knot. By the above, we can assume
without loss of generality that it is a closed braid.  Let $\cT
\cK$ be its transversal knot type, i.e., its knot type under
transversal isotopies, and let $[\cT \cK]_{top}$ be its
topological knot type.  Choose a representative $TX$ of $\cT\cK$,
which (by Bennequin's transversal version of Alexander's theorem)
is always possible.  Choose a Seifert surface $F$ of minimal
genus, with $TX = \partial F$.  Bennequin studied the foliation of
$F$ which is induced by the intersections of $F$ with the plane
field determined by $\xi_\pi$.     Let $n(TX)$ be the braid index
and let $e(TX)$ be the algebraic crossing number of a projection
of $TX$ onto the plane $z=0$.  Both can be determined from the
foliation.  Bennequin found an invariant of $\cT\cK$, given by the
formula $\beta(\cT\cK) = e(TX) - n(TX)$.  Of course if he had
known Theorem \ref{theorem:TMT}  the proof that $\beta(\cT\cK)$ is
an invariant of $\cT\cK$ would have been trivial, but he did not
have that tool.  He then showed a little bit more: he showed that
$\beta(\cT\cK)$ is bounded above by the negative of the Euler
characteristic of $F$ in $\xi_\pi$.  He then showed that this
bound fails in one of the contact structures $\xi_{>\pi}$. In this
way he proved that the contact structures  $\xi_{>\pi}$ cannot be
isotopic to $\xi_\pi$.

To knot theorists, Bennequin's proof should seem intuitively
natural, because the invariant $\beta(\cT\cK)$ is a self-linking
number of a representative $TX\in \cT\cK$ (the sense of push-off
being determined by $\xi$), and the more twisting there is the
higher this number can be.  For an explanation of the
self-linking, and lots more about Legendrian and transversal knots
we refer the reader to John Etnyre's excellent review article
\cite{Etnyre2003}. The basic idea is that $TX$ bounds a Seifert
surface, and this Seifert surface is foliated by  the plane field
associated to $\xi$. Call this foliation the
\underline{characteristic foliation}. Near the boundary, the
characteristic foliation is transverse to the boundary. The
Bennequin invariant is the linking number of $TX, TX'$, where
$TX'$ is is a copy of $TX$, obtained by pushing $TX$ off itself
onto $F$, using the direction determined by the characteristic
foliation of $F$.

Bennequin's paper was truly important.   Shortly after it was
written Eliashberg showed in \cite{Eliashberg1989} that the
phenomenon of an infinite sequence of contact structures related
to a single one of minimal twist angle occured generically in
every 3-manifold, and introduced the term `tight' and
`overtwisted' to distinguish the two cases.   Here too, there is a
reason that will seem natural to topologists. In 2003 Giroux
proved \cite{Giroux} that every contact structure on every closed,
orientable 3-manifold $M^3$ can be obtained in the following way:
represent $M^3$ as a branched covering space of $S^3$, branched
over a knot or link, and lift the standard and overtwisted contact
structures on $S^3$  to $M^3$.

Returning to knot theory, the invariant $\beta(\cT\cK)$ allows us to answer a fundamental question: is the equivalence relation on knots that is defined by transversal isotopy really different from the equivalence relation defined by topological isotopy?   The Bennequin invariant will be used to answer this question in the affirmative:
\begin{theorem}
\label{theorem:transversal isotopy}{\rm \cite{Bennequin}} There
are infinitely many distinct transversal knot types associated to
each topological knot type.
\end{theorem}

\pf  Choose a transversal knot type $\cT\cK$ and a closed braid
representative $TX_0$.  Stabilizing the closed braid $TX_0$ once
negatively (recall the definition of positive/negative
stabilizations given in $\S$\ref{subsection:Markov's Theorem},
Figure~\ref{figure:stab-destab}), we obtain the transverse closed
braid $TX_1$, with $e(TX_1) = e(TX_0)-1$ and $n(TX_1) = n(TX_0) +
1$, so that $\beta(TX_1) = \beta(TX_0) - 2$. Iterating, we obtain
transverse closed braids $TX_2, TX_3, \dots$, defining transverse
knot types $\cT\cK_1, \cT\cK_2, \cT\cK_3, \cdots$,  and no two
have the same Bennequin invariant. Since stabilization does not
change the topological knot type, the assertion follows.  \endpf

At this writing, it is an open problem to find computable
invariants of $\cT\cK$ which are not determined by
$[\cT\cK]_{top}$ and $\beta(\cT\cK)$.   A hint that the problem
might turn out to be quite subtle was in the paper \cite{F-T1997}
by Fuchs and Tabachnikov, who proved that while ragbags filled
with polynomial and finite type invariants of transversal knot
types $\cT\cK$ exist, based upon the work of Arnold in
\cite{Arnold}, they are all determined by $[\cT\cK]_{top}$ and
$\beta(\cT\cK)$. Thus, the seemingly new invariants that many
people had discovered by using Arnold's ideas were just a fancy
way of encoding $[\cT\cK]_{top}$ and $\beta(\cT\cK)$.

This leads naturally to a question:  Are there computable invariants of transversal knots which are {\it not} determined by  $[\cT\cK]_{top}$ and $\beta(\cT\cK)$?  A similar question arises in the setting of Legendrian knots.  Each Legendrian knot $\cL\cK$ determines a topological knot type $[\cL\cK]_{top}$, and just as in the transverse case it is an invariant of the Legendrian knot type. There are also two numerical invariants of $\cL\cK$: the Thurston-Bennequin invariant $tb(\cL\cK)$ (a self-linking number) and the Maslov index $M(\cL\cK)$ (a rotation number).   So until a few years ago the same question existed in the Legendrian setting, but the Legendrian case has recently been settled by Yuri Chekanov:

\begin{theorem} {\rm  \cite{Chekanov}} There exist distinct Legendian knot types  which have the same topological knot type$[\cL\cK]_{top}$, and also the same Thurston-Bennequin invariant $tb(\cL\cK)$ and Maslov index $M(\cL\cK)$.
\end{theorem}

The analogous result for transversal knots proved to be quite difficult, so to begin to understand whether something could be done via braid theory, the first author and Nancy Wrinkle asked an easier question, which they answered in part in \cite{Birman-Wrinkle}: are there transversal knot types which {\it are determined}  by their topological knot type and Bennequin number?  This question lead to a definition: a transversal knot type $\cT\cK$ is \underline{transversally simple} if it is determined by
$[\cT\cK]_{top}$ and $\beta(\cT\cK)$.  So the question is: are there transversally simple knots?  The manuscript \cite{Birman-Wrinkle} gives a purely topological (in fact braid-theoretic) criterion which enables one to answer the question affirmatively, adding one more piece of evidence that topology and analysis walk hand in hand.  To explain it, recall the destabilization and exchange moves, and the two flypes, depicted in Figure \ref{figure:templates}.    A topological knot or link type  $K$ is said to be \underline{exchange reducible} if an arbitrary closed braid representative $X$ of $K$ can be changed to an arbitrary representative of minimum braid index
by braid isotopy, positive and negative destabilizations and exchange moves.  We have:

\begin{theorem} {\rm  \cite{Birman-Wrinkle}}
\label{theorem:exchange-reducibility}
If a knot type $K$ is exchange-reducible, then any transversal knot type $\cT\cK$ which has $[\cT\cK]_{top} = K$  is transversally simple.
\end{theorem}
This theorem was used to give a new proof of a theorem of
Eliashberg \cite{Eliashberg1993}, which asserts that the unlink is
transversally simple, and also (with the help of
\cite{Menasco2001}) to prove the then-new result that most
iterated torus knots are transversally simple.

The rest of this section will be directed at explaining the main result of \cite{BM-stab-II}:
\begin{theorem}{\rm \cite{BM-stab-II}}
\label{theorem:negative flype examples}
There exist transversal knot types which are not transversally simple. Explicitly, the  transverse closed $3$-braids
$TX_+ =  \sigma_1^5 \sigma_2^4 \sigma_1^6\sigma_2^{-1}$ and
  $TX_- = \sigma_1^5 \sigma_2^{-1} \sigma_1^6\sigma_2^4$
determine transverse knot types $\cT\cK_+, \cT\cK_-$ with $(\cT\cK_+)_{top} = (\cT\cK_-)_{top}$ and $\beta(\cT\cK_+) = \beta(\cT\cK_-)$, but  $\cT\cK_+ \not= \cT\cK_-$.
\end{theorem}

\noindent {\bf Sketch of the proof of Theorem
\ref{theorem:negative flype examples}.} See  \cite{BM-stab-II} for
all details.
The examples in Theorem
\ref{theorem:negative flype examples} were obtained by choosing all the weights in the negative flype template to be 1, and assigning
explicit 2-braids to the blocks $P, Q, R$ of the negative flype
template of Figure \ref{figure:templates}.  If
the weights are all chosen to be 1, the blocks $P,Q,R$ are 2-braids and (except in very special cases)  the flype templates have braid index 3.  First one must show that the examples satisfy the
conditions of the theorem.  The topological knot types defined by
the closed 3-braids $ \sigma_1^5 \sigma_2^4
\sigma_1^6\sigma_2^{-1}$ and $\sigma_1^5 \sigma_2^{-1}
\sigma_1^6\sigma_2^4$  coincide because they are carried by the
block strand diagrams for the negative flype template of Figure
\ref{figure:templates}.  The  Bennequin invariant can be computed
as the exponent sum of the braid word (14 in both cases) minus the
braid index (3 in both cases).  So the examples have the required
properties.

The hard part is the establishment of a special version of Theorem
\ref{theorem:MTWS} which is applicable to the situation.  Its
special features are as follows:

\be
\item Both $X_+$ and $X_-$ have braid index 3.
\item Since it
is well known that exchange moves can be replaced by braid isotopy
for 3-braids,   the first two sequences in Theorem
\ref{theorem:MTWS} are vacuous, i.e. $X_\pm =X_\pm^\prime$.

\item Because of the special assumption just noted, the templates
that are needed, in the topological setting, can be enumerated
explicitly: they are the positive and negative destabilization and the positive and negative
flype templates. No others are needed.
\ee

It is proved in
\cite{BM-stab-II} that if $X_-$ and $X_+$ are transversal closed
braids $TX_+$ and $TX_-$, then the isotopy that takes $TX_+$ to
$TX_-$ may be assumed to be transversal.
 So, suppose that a transversal isotopy exists from the transverse
closed braid $TX_+$ to the transverse closed braid $TX_-$.  Then
there is a $3$-braid template that
carries the braids $ \sigma_1^5 \sigma_2^4
\sigma_1^6\sigma_2^{-1}$ and $\sigma_1^5 \sigma_2^{-1}
\sigma_1^6\sigma_2^4$.  This is the first key fact that we use
from Theorem \ref{theorem:MTWS}.  Instead of having to consider
{\it all} possible transversal isotopies from $TX_+$ to $TX_-$, we
only need to consider those that relate the left and right
block-strand diagrams in one of the four 3-braid templates.  So the braids in question are carried either by one of the
two destabilization templates or by one of the two flype
templates.   If it was one of the destabilization templates, then
the knots in question could be represented by 2 or 1-braids, i.e.
they would be type $(2,n)$ torus knots or the unknot, however  an
easy argument shows that the knots in Theorem
\ref{theorem:negative flype examples} are neither type $(2,n)$
torus knots or the unknot.  The positive flype templates are ruled
out in different way:  topologically, our closed braids admit a
negative flype, so if they are also carried by the positive flype
template they admit flypes of both signs.  However, the manuscript
\cite{K-L1999} gives conditions under which a closed 3-braid
admits flypes of both signs, and the examples were chosen
explicitly to rule out that possibility.

We are reduced to isotopies that are supported by the negative flype template.  It is straightforward to show that the obvious isotopy is not transversal, but maybe there is some {\it other} isotopy which is transversal. Here a key fact about the definition of a template is used (and this is a second very strong
aspect of the MTWS).  If such a transversal  isotopy exists, then it exists for every knot or  link defined by a fixed choice of braiding assignments to the blocks.  Choose the braiding assignments
$\sigma_1^3, \sigma_2^4, \sigma_1^{-5}$ to the blocks
$P,R,Q$.  This braiding assignment gives a 2-component link
$L_1\sqcup L_2$ which has two distinct isotopy classes of closed 3-braid representatives. If $L_1$ is the
component associated to the left strand entering the block $P$, then
$\beta(L_1)=-1$ and $\beta(L_2) = -3$ before the flype, but after the flype the
representative will be $\sigma_1^3 \sigma_2^{-1} \sigma_1^{-5} \sigma_2^{4}$, with
$\beta(L_1)=-3$ and $\beta(L_2) = -1$. However, by Eliashberg's isotopy extension theorem (Proposition
2.1.2 of \cite{Eliashberg1993}) a  transversal isotopy of a knot/link extends to an ambient
transversal isotopy  of the 3-sphere.  Any transversal isotopy of $L_1 \sqcup L_2$
must preserve the $\beta$-invariants of the components. It follows that no such transversal isotopy exists, a contradiction of our assumption that $TX_+$ and $TX_-$ are transversally isotopic. \endpf

Other examples of a similar nature were discovered by Etnyre and Honda \cite{E-H2003} after the proof of Theorem \ref{theorem:negative flype examples} was posted on the arXiv.  Their methods are very different from the proof that we just described (being based on contact theory techniques rather than topological techniques), but are equally indirect. They do not produce explicit examples, rather they present a bag of pairs of transverse knots and prove that at least one pair in the bag exists with the properties given by Theorem \ref{theorem:negative flype examples}.  Therefore we pose, as an important open problem:

\begin{openproblem}
{\rm Find new computable invariants of transversal knot types.  Here `new' means an invariant which is not determined by the topological knot type ${\cal TX}$ and and the Bennequin invariant $\beta({\cal TX})$.  `Computable' means that it should be computable from either a closed braid diagram or some other representation of the transversal knot.  Braid groups seem to be a natural setting for investigating this problem. }
$\clubsuit$
\end{openproblem}

\begin{remark}
\label{remark:open book structures}
{\rm The connections between closed braids and contact structures does not end with transversal knots.
There are fundamental relationships between open book structures on 3-manifolds and contact structures on those manifolds, both untwisted and twisted. See \cite{Giroux} for an introduction to this interesting new area, and see \cite{Etnyre2003} for a review of the mathematics and a discussion of many open problems waiting to be investigated.}
\end{remark}

\newpage
\section{Representations of the braid groups}
\label{section:representations of the braid groups}

Before the discovery of Hecke algbra representations of the braid group (discussed in this section) very little was known about finite dimensional but infinite representations of ${\bf B}_n$, except for the ubiquitous Burau representation. That matter changed dramatically in 1987 with the publication of \cite{Jones1987}.  Suddenly, we had more knot invariants and with them more braid group representations than anyone could deal with, and the issue became one of organizing them. However, we shall not attempt to give a comprehensive overview of the rich
theory of representations of braid groups in this section.
Instead, we focus here on the representations of ${\bf B}_n$ which
have played the greatest roles in the development of that theory:
the Burau representation, the Hecke algebra representations, and, more recently, the Lawrence-Krammer  representation.

\subsection{A brief look at representations of $\Sigma_n$}
\label{subsection:symm}

The fact that the representation theory of ${\bf B}_n$ is rooted in
the representation theory of the symmetric group $\sgn$ is a
consequence of the surjection ${\bf B}_n \rightarrow \sgn$, given by mapping the elementary braid $\sigma_i$ to the transposition $s_i = (i,i+1)$.   While the kernel of this homomorphism is very big (it's the entire pure braid group ${\bf P}_n$), it nevertheless turns out that a great deal can be learned about representations of ${\bf B}_n$ by studying the collection of irreducible representations of $\sgn$, and attempting to lift them to reprentations of ${\bf B}_n$ by `deforming' them, and hoping that  you get something new.  For the record, we now note that the group algebra $\complexes\sgn$ has generators $1, s_1,\dots,s_{n-1}$ and defining relations:
\begin{equation} \label{equation:group algebra, symmetric group}
s_i  s_k = s_k s_i \   {\rm if} \  |i-k| \geq 2, \  \  s_i s_{i+1} s_i =  s_{i+1} s_i s_{i+1}, \ s_i^2 = 1,
\end{equation}
where $1 \leq i \not= k \leq n-1.$  As a vector space, $\complexes\sgn$ is spanned by $n!$ reduced words \cite{Bourbaki}  in the transpositions $s_i$:
\begin{equation}
\label{equation:basis for group algebra of symmetric gp}
\{(s_{i_1} s_{i_1 - 1} \ldots s_{i_1 - k_1})(s_{i_2} s_{i_2 - 1}
\ldots s_{i_2 - k_2}) \cdots (s_{i_r} x_{i_r - 1} \ldots s_{i_r -
k_r}) \}
\end{equation}
 where $1 \leq i_1 < i_2 < \cdots < i_r \leq n-1 $ and $i_j - k_j \geq 1$.

\begin{figure}[htpb!]
\centerline{\includegraphics[scale=.3,  bb=232 322 379 469] {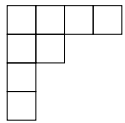}}
%\centerline{\includegraphics[scale=.3] {youngdiagram.eps}}
\caption{The Young diagram corresponding to the partition
$4+2+1+1$.} \label{figure:youngdiagram}
\end{figure}

Irreducible representations of $\sgn$ are
parametrized by  \underline{Young diagrams}, which in turn are
parametrized by partitions of $n$. (See \cite{Ful-Har} for one of
the many good discussions of this subject in the literature.) A
Young diagram consists of stacked rows of boxes, aligned on the
left, with the number of boxes in each row strictly nonincreasing
as you go from top to bottom.  The number of boxes in a row
corresponds to a term in a given partition of $n$.  For example,
the Young diagram shown in Figure~\ref{figure:youngdiagram}
corresponds to the partition $8=4+2+1+1$.

As it turns out, the irreducible representations of the symmetric
group are in 1-1 correspondence with Young diagrams, or
equivalently with partitions of $n$.  There is a row-column
symmetry, and we adopt the convention that the Young diagram which
has one row of $n$ boxes, which comes from the partition $n = n$,
corresponds to the trivial representation. Then the Young diagram
which has $n$ rows consisting of 1 box per row, which comes from
the partition $n = 1 + 1 + \cdots + 1$, corresponds to the parity
representation, mapping each element of $\sgn$ to its sign. The
first interesting representation of a braid group occurs when $n =
3$, corresponding to the partition $3 = 2 + 1$. More generally, we
consider the partition $n = (n - 1) + 1$, or the Young diagram
with two rows, the top row having $(n -1)$ boxes and the bottom
having a single box.  This Young diagram corresponds to the
\underline{standard representation} $\sgn \rightarrow
GL_n(\complexes)$, given by the usual action of $\sgn$ on
$\complexes^n$ by permuting basis vectors.  Thus if $s_i$ denotes
the transposition $(i \hskip .1in i+1) \in \sgn$, the standard
representation sends $s_i$ to the following:
$$
I_{i-1} \oplus \left( \begin{array}{cc}
0 & 1\\
1 & 0
\end{array}
\right) \oplus I_{n-i-1}
$$
where $I_k$ denotes the $k \times k$ identity matrix.

The standard representation is reducible; it is easy to see that
it fixes a $1$-dimensional subspace, namely, the span of the sum
of the basis vectors.  The complementary $(n-1)$-dimensional
subspace is the set of all points $(z_1, \ldots, z_n) \in
\complexes^n$ such that $z_1 + \ldots + z_n = 0$. This subspace
corresponds to the irreducible representation given by the Young
diagram in question.  One can use the \underline{ hook length formula} to
calculate directly that the dimension of the representation
corresponding to this Young diagram is indeed $n - 1$ (we refer
the reader again to \cite{Ful-Har} for an explanation of this
formula).

\subsection{The Burau representation and polynomial invariants of knots.}\label{subsection:The Burau
representation}

Burau first introduced his representation of the braid group in
1936 \cite{Burau1936}.  Much later, it was realized that it could be thought of as a deformation of the standard representation of $\sgn$ corresponding to the partition $n = (n-1) + 1$.  For many years it was the
focus of the representation theory of braid groups.   We define
the \underline{Burau representation} $\rho: {\bf B}_n \rightarrow
GL_n(\ints [t, t^{-1}])$ as follows:
$$
\sigma_i \mapsto I_{i-1} \oplus \left( \begin{array}{cc}
1-t & t\\
1 & 0
\end{array}
\right) \oplus I_{n-i-1}
$$
Note that substituting $t = 1$ gives back the
representation (of ${\bf B}_n$ factoring through $\sgn$), and this is why we say that it is a deformation of the standard representation of  $\sgn$.  Like the
representation of $\sgn$, the Burau representation splits into a
$1$-dimensional representation and an $(n-1)$-dimensional
irreducible representation known as the \underline{reduced Burau
representation} which we denote by $\bar{\rho}: {\bf B}_n \rightarrow
GL_n(\ints [t, t^{-1}])$ as follows:
$$
\sigma_i \mapsto I_{i-2} \oplus \left( \begin{array}{ccc}
1 & -t & 0\\
0 & -t & 0\\
0 & -1 & 1
\end{array}
\right) \oplus I_{n-i-2}
$$
where the $-t$ in the middle of the $3 \times 3$ matrix is always
in the $(i,i)^{th}$ spot.

It has been known for a long time that the Burau representation is
faithful for $n \leq 3$ (see \cite{Birman1974}, for example), and
for many years the representation was held to be a reasonable
candidate for a faithful linear representation of ${\bf B}_n$ for all
$n$.  However, in 1991, Moody showed that for $n \geq 9$, the
Burau representation is {\em not} faithful \cite{Moody1991}.  Long
and Paton later improved Moody's result to $n \geq 6$
\cite{Lo-Pa}. Bigelow further improved this to $n \geq 5$
\cite{Bigelow1999}. At the time of this writing, the case $n = 4$
remains open.

Despite such results, the Burau representation continues to play
an important role in the study of representations of ${\bf B}_n$.  It
was known classically that the Alexander polynomial $\Delta_K(t)$
of a knot or link $K$ can be calculated directly from the image
under the reduced Burau representation of a braid $X$  such that
$b(X)$ represents $K$, as follows (see \cite{Birman1974}, e.g.,
for a proof based upon Theorem \ref{theorem:MT}, the Markov Theorem):
\begin{equation} \label{equation:alexander from burau}
\Delta_{b(X)}(t) = \frac{\textrm{det}( \bar{\rho} (X) - I_{n-1}
)}{1 + t + \cdots + t^{n-1}}.
\end{equation}
Thus the Alexander polynomial of the closed braid associated to the open braid $X$, i.e. $\Delta_{b(X)}(t)$, is a rescaling of the characteristic polynomial of the the image of $X$ in the reduced representation.
In what follows, we shall see how a property of the Burau
representation motivated the definition of a variant on the
two-variable HOMFLY polynomial \cite{HOMFLY}, a knot invariant of
which both the Alexander polynomial and the Jones polynomial are
specializations. Further, in $\S$\ref{subsubsection:A topological
interpretation of the Burau representation} we shall give a
topological interpretation of the Burau representation which
naturally leads to the definition of a faithful linear
representation of ${\bf B}_n$ known as the Lawrence-Krammer
representation, which will be our focus in $\S$\ref{subsection:the
Lawrence-Krammer representation}.

\subsection{Hecke algebras representations of braid groups and polynomial invariants of
knots}\label{subsubsection:Hecke algebras and polynomial
invariants of knots} A simple calculation, together with the
Cayley-Hamilton theorem, shows that the image of each of our
braid group generators under the Burau representation,
$\rho(\sigma_i)$, satisfies the characteristic equation $x^2 = (1
- t) x + t$ and thus has two distinct eigenvalues.  This prompted
Jones to study all representations $\rho: {\bf B}_n\rightarrow
GL_n(\complexes)$ which have at most two distinct eigenvalues
\cite{Jones1987}. Let $x_i = \rho (\sigma_i)$.  Then for all $i$,
$x_i$ must satisfy a quadratic equation of the form $x_i^2 + a x_i
+ b = 0$.  By rescaling, we may assume that one of the eigenvalues
is $1$ and eliminate one of the variables, e.g., $a = -(1+b)$.
Note that by rewriting our quadratic equation and making the
substitution $b = -t$ we regain the characteristic equation from
the Burau representation. However, the convention in the
literature seems to be to rescale our representation by $(-1)$ so
that the equation takes the form  $x_i^2 = (t - 1) x_i + t.$
With this motivation, we define the \underline{Hecke algebra}
$H_n(t)$ to be the algebra with generators
$1, x_1, \ldots, x_{n-1}$
and defining relations as follows:
\begin{equation} \label{equation:Hecke algebra relations}
x_i  x_k =x_k x_i \  \ {\rm if} \  \ |i-k| \geq 2, \  \  x_i x_{i+1} x_i =  x_{i+1} x_i x_{i+1}, \ \  x_i^2 = (t-1)x_i + t ,
\end{equation}
where $1 \leq i \not= k \leq n-1.$     Comparing the relations in (\ref{equation:group algebra, symmetric group}) and (\ref{equation:Hecke algebra relations}), we see that $H_n(1) \cong \complexes \sgn$, the group
algebra of the symmetric group.  Hence we can think of $H_n(t)$ as
a `deformation' of $\complexes \sgn$.

The connection between $H_n(t)$ and $\complexes\sgn$ is made even more transparent by noting that as a vector space, $H_n(t)$ is spanned by $n!$ lifts
of a system of reduced words in the transpositions $s_i \in \sgn$.
For example, we can take as a spanning set
\begin{equation}
\label{equation:spanning set for Hecke algebra}
\{(x_{i_1} x_{i_1 - 1} \ldots x_{i_1 - k_1})(x_{i_2} x_{i_2 - 1}
\ldots x_{i_2 - k_2}) \cdots (x_{i_r} x_{i_r - 1} \ldots x_{i_r -
k_r}) \}
 \end{equation}
  where $1 \leq i_1 < i_2 < \cdots < i_r \leq n-1 $ and $i_j - k_j \geq 1$
\cite{Bourbaki}, \cite{Jones1987}.

\begin{remark}
\textnormal{In $\S$\ref{subsubsection:braiding in operator algebras}  we defined an algebra $J_n(t)$ with generators $1,g_1,\dots,g_{n-1}$ and defining relations (\ref{equation:Jones algebra relations}),  the \underline{Jones algebra}. As it turns out, its irreducible summands are in 1-1 correspondence with the irreducible representations of the Hecke algebra that are parametrized by Young diagrams with exactly 2 rows. The Hecke algebra, as defined above, has generators maps $1, x_1,\dots,x_{n-1}$ with defining relations (\ref{equation:Hecke algebra relations}), so that the map $\xi:H_n(t)\to J_n(t)$ that is defined by $\xi(x_i) = g_i$ is a homorphism of algebras.  }
\end{remark}

Our main purpose in this section is to outline Jones' development
in \cite{Jones1987} of a two-variable polynomial knot invariant
arising from representations of the Hecke algebras $H_n(t)$.  This
polynomial is essentially the well-known HOMFLY polynomial, and
includes the Jones polynomial as a specialization. We have just
seen that the Jones algebra is the quotient of the Hecke algebra
$H_n(t)$ by one extra relation (as an aside, we note that this
extra relation is satisfied by the image of the transpositions
generators $s_i = (i, i+1)$ in $\sgn$ under all representations
arising from Young diagrams with at most two rows.)  To pick up a
second variable, we introduce an additional parameter by allowing
a 1-parameter family of traces on Hecke algebras. We now pursue
this point of view, and we also refer the reader to
\cite{delaHarpe} for another exposition which follows the same
point of view.

We begin by defining a function $f: {\bf B}_n \rightarrow H_n(t)$ by
$f(\sigma_i) = x_i$.  The function $f$ is well-defined on reduced
words in the generators $\sigma_i$ and commutes with the natural
inclusions ${\bf B}_{n - 1} \subset {\bf B}_n$ and $H_{n-1}(t) \subset
H_{n}(t)$, although in general $f$ fails to be a homomorphism. We
can then apply the following result due to Adrian Ocneanu which appeared
in \cite{HOMFLY} and was proved inductively in \cite{Jones1987}
using the $n!$-element basis given above.

\begin{theorem} {\rm \cite{HOMFLY}, \cite{Jones1987}}
\label{theorem:trace} For each $z \in
\complexes^*$ (and each $t \in \complexes^*$), there exists a
unique trace function $\textnormal{tr}: \cup_{n=1}^\infty H_n(t)
\rightarrow \complexes$ such that
\begin{enumerate}
\item $\textnormal{tr}(1) = 1$
\item $\textnormal{tr}(ab) = \textnormal{tr}(ba)$
\item $\textnormal{tr}$ is $\complexes$-linear
\item $\textnormal{tr}(u x_{n-1} v) = z \hskip .05in \textnormal{tr}(uv)$ for all $u,v
\in H_{n - 1}(t)$.
\end{enumerate}
\end{theorem}

Theorem~\ref{theorem:trace} gives us a one-parameter family of
trace functions on a one-parameter family of algebras.  In fact,
using the properties of the trace function given in theorem it is
possible to compute $\textrm{tr} ( f(X))$ for all $X \in
{\bf B}_n$. (We note the fact that for any $w \in H_n(t)$ such that $w
\notin H_{n-1}(t)$, there is a unique reduced word $w = x_{i_1}
\cdots x_{i_r}$ in which $x_{n - 1}$ appears exactly once
\cite{Jones1987}.) In practice, the third relation of $H_n(t)$ is
quite useful for computing the trace function $\textnormal{tr}$,
both in its original form and in the following:
$$
x_i^{-1} = t^{-1} x_i + (t^{-1} - 1).
$$

\begin{example}\label{example:trace}{\rm
Let $X_1 = \sigma_1^3 \in {\bf B}_2$, and let $X_2 = \sigma_1
\sigma_2^{-1} \sigma_1 \sigma_2^{-1} \in {\bf B}_3$.  Note that
$\hat{X_1}$ is the right-handed trefoil knot and that
$\hat{X_2}$ is the Figure-8 knot.  We invite the reader to
check that $$\textrm{tr}(f(X_1)) = (t^2 - t + 1)z + t(t-1)$$
and that $$\textrm{tr}(f(X_2)) = (3 - t^{-1} - t) t^{-1} z^2 +
(3 - t^{-1} -t) (t^{-1} - 1) z - (2 - t^{-1} - t).  \spadesuit  $$}
\end{example}

We also note that the second property of the trace function given
in Theorem~\ref{theorem:trace} implies that $\textnormal{tr} \circ
f $ is invariant on conjugacy classes in ${\bf B}_n$.  It remains to
tweak the function a bit in order to obtain from a given braid a
two-variable polynomial which is also invariant under
stabilization and destabilization moves as defined in
$\S$\ref{subsection:Markov's Theorem}; such a polynomial will be
Markov-invariant and hence an invariant of the knot type of the
closed braid.

Algebraically, stabilization and destabilization each take the
form $X \rightarrow X \sigma_n^{\pm 1}$, the only
difference being appropriate conditions on the braid $X$.  We
would like to rescale our representation $f$ in such a way that
both versions of stabilization (resp. destabilization) have the
same effect on the trace function. Suppose there exists a complex
number $k$ such that $\textrm{tr}(k x_i)= \textrm{tr}((k
x_i)^{-1}))$. Then we can find a `formula' for $k$ as follows:
\begin{eqnarray*}
k^2 \textrm{tr}(x_i) &=& \textrm{tr}(x_i^{-1}) \\
k^2 z &=& \textrm{tr}(t^{-1} x_i + t^{-1} - 1) \\
k^2 &=& \frac{t^{-1} z + t^{-1} - 1}{z} \\
k^2 &=& \frac{1 + z - t}{tz}
\end{eqnarray*}
Solving this for $z$, we obtain $$ z = -\frac{1-t}{1 - k^2 t}.$$
We set $\kappa = k^2$, and define $f_\kappa: {\bf B}_n \rightarrow
H_n(t)$ by $f_\kappa(\sigma_i) = \sqrt{\kappa} \hskip .05in
\sigma_i$.  Now we have
\begin{eqnarray*}
\textrm{tr}(f_\kappa (\sigma_n)) &=& \sqrt{\kappa} \hskip .05in z \\
&=&  -\sqrt{\kappa}\frac{1-t}{1 - \kappa t}.
\end{eqnarray*}

We would like to define a map  ${\bf B}_n \rightarrow \ints[t^{\pm1},
\kappa^{\pm 1}]$ which is Markov invariant.  At the moment, we
have that
$$
\textrm{tr}(f(w\cdot \sqrt{\kappa} \sigma_n)) =
-\sqrt{\kappa}\frac{1-t}{1 - \kappa t} \textrm{tr} (f(w)) =
\textrm{tr}(f(w \cdot \frac{1}{ \sqrt{\kappa} }\sigma_n^{-1}))
$$
for any $w \in {\bf B}_n$.  Now we simply define
\begin{eqnarray*}
F(X) = F_X(t,\kappa) &=& (-\frac{1}{\sqrt{\kappa}} \cdot
\frac{1 - \kappa t}{1 - t})^{n-1} \textrm{tr}(f_\kappa(X))\\
&=&(-\frac{1}{\sqrt{\kappa}} \cdot \frac{1 - \kappa t}{1 -
t})^{n-1}(\sqrt{\kappa})^E \textrm{tr}(f(X))
\end{eqnarray*}
for $X \in {\bf B}_n$, where $E$ is the exponent sum of $X$ as a word
in $\sigma_1, \ldots, \sigma_{n-1}$.  It is clear that $F(X)$
depends only on the knot type of $b(X)$. We now reparametrize one
last time, setting
\begin{eqnarray*}
l &=& \sqrt{\kappa} \sqrt{t} \\
m &=& \sqrt{t} - \frac{1}{\sqrt{t}}.
\end{eqnarray*}

\begin{figure}[htpb!]
\centerline{\includegraphics[scale=.6,  bb= 120 504 444 562] {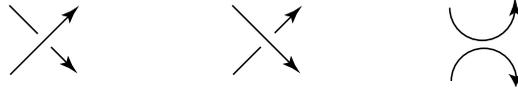}}
%\centerline{\includegraphics[scale=.6] {skein.eps}}
\caption{Crossings for $K_+, K_-$, and $K_0$, respectively, in
the skein relation given in (\ref{equation:skein}). }
\label{figure:skein}
\end{figure}

With this substitution, we obtain a Laurent polynomial in two
variables $l$ and $m$, which we denote $P_{b(X)}(l,m) = P_K
(l,m)$, where $K$ is the (oriented) knot or link type of $b(X)$.
Furthermore, $P_K(l,m)$ satisfies the skein relation
\begin{equation}\label{equation:skein}
m P_{K_0}= l^{-1} P_{K_+} - l P_{K_-} \end{equation} where $K_0,
K_+$, and $K_-$ are oriented knots with identical diagrams except
in a neighborhood of one crossing, where they have a diagram as
given in Figure~\ref{figure:skein} (Proposition 6.2 of
\cite{Jones1987}). Thus by beginning with $P_U = 1$, where $U$
denotes the unknot, it is possible to calculate $P_K$ for any knot
or link $K$ using only the skein relation, which is often simpler
than using the trace function.

\begin{remark}\textnormal{
We note that the polynomial $P_K (l,m)$ obtained in this way is
essentially the same as the two-variable polynomial known as the
HOMFLY polynomial (\cite{HOMFLY}) which is usually reparametrized
as $P_K(i l^{-1}, im)$.}\end{remark}

\begin{example}\label{example:polynomial} \textnormal{
Let $X_1, X_2$ be the braids defined in
Example~\ref{example:trace} whose closures are $K_1 =$ the
right-handed trefoil knot and $K_2 =$ the Figure-8 knot,
respectively. We leave it as an exercise for the reader to check
that
\begin{eqnarray*}
F_{X_1} (t, \kappa) &=& \kappa(1 + t^2 - \kappa t^2) \\
&=& \kappa t (t + t^{-1} - \kappa t) \\
&=& \kappa t (2 - \kappa t + t + t^{-1} -2) \\
&=& \kappa t \left(2 - \kappa t + (\sqrt{t} - \frac{1}{\sqrt{t}})^2 \right) \\
\end{eqnarray*}
and hence we have $$P_{K_1} (l,m) = 2 l^2 - l^4 + l^2 m^2.$$
Similarly, the reader can check that
 \begin{eqnarray*}
F_{X_2} (t, \kappa) &=& \frac{1 - \kappa(1 - t + t^2) +
\kappa^2 t^2} {t \kappa}\\
P_{K_2} (l,m) &=& l^{-2} - m^2 - 1 + l^2.
\end{eqnarray*}
For explicit calculations of $F_{X_2}$ and $P_{K_2}$ using the
trace function, see p. 350 of \cite{Jones1987}.}$\spadesuit$
\end{example}

\begin{remark}
\textnormal{ There is also a 1-variable knot polynomial, the \underline{Jones polynomial}, associated to the algebra $J_n(t)$ generated by $1,g_1,\dots,g_{n-1}$ with defining relations (\ref{equation:Jones algebra relations}).  In the situation of the Jones algebra the trace is unique, whereas in the situation of the Hecke algebra, as we presented it here, there is a 1-parameter family of traces. The 1-variable Jones polynomial was discovered before the 2-variable HOMFLY polynomial.}
\end{remark}

This two-variable knot polynomial has been much studied and
reviewed in the literature.  For the sake of completeness we list
here a few of its noteworthy properties and applications.
\begin{enumerate}
\item  {\bf Connect sums}: $P_{K_1 \# K_2} = P_{K_1} \cdot P_{K_2}$
\item {\bf Disjoint unions}:  $P_{K_1 \amalg K_2} = (\frac{l^{-1} - l}{m}) P_{K_1} \cdot P_{K_2}$.
\item  {\bf Orientation}:  $P_{\bar{K}} = P_K$, where $\bar{K}$ denotes the link obtained by
reversing the orientation of every component of the link $K$.
\item {\bf Chirality}: $P_{\tilde{K}} (l,m) = P_K(l^{-1},-m)$,
where $\tilde{K}$ denotes the mirror image of the link $K$.
\item {\bf Alexander polynomial}:  Note that $F_X(1, \kappa)$ is
not defined. It comes as something of surprise, then, that the
specialization $l = 1, m = \sqrt{t} - \frac{1}{\sqrt{t}}$ gives
the Alexander polynomial $\Delta_K(t)$. Jones shows how to avoid
the singularity by exploiting an alternate method of calculating
the trace function using weighted sums of traces (see
\cite{Wenzl1985} and \cite{KL1979} as well as \cite{HOMFLY} and
\cite{Jones1987}). A by-product of this alternate method is
another derivation of Equation~\ref{equation:alexander from burau}
showing how to calculate $\Delta_K(t)$ from the Burau
representation.
\item {\bf Jones polynomial}: The famous Jones polynomial can be
obtained from the two-variable polynomial by setting
$$
V_K(t) = P(t, \sqrt{t} - \frac{1}{\sqrt{t}}).
$$
Note that we are abusing notation by reusing the variable $t$ here
and above in $\Delta_K(t)$.
\item  {\bf A lower bound for braid index}  At the end of  $\S$\ref{subsection:Alexander's Theorem} we remarked that it is an open problem to determine the braid index of a knot algorithmically. However, the HOMFLY polynomial does give a remarkably useful lower bound, via a famous inequality which is known as the Morton-Franks-Williams inequality. It was proved simultaneously and independently by Hugh Morton in \cite{Morton1986_2} and by John Franks and Robert Williams in \cite{F-W}. While it proved to be sharp on all but 5 of the knots in the standard tables of knots having at most 10 crossings, there are also infinitely many knots on which it fails to be sharp.  The first author came to a new appreciation of the importance of this problem when she was faced with the problem of determining, precisely, the braid index of certain 6-braid knots on which it failed.  Consulting many people, it became apparent that there was essentially no other useful result on this basic question.
\end{enumerate}

\subsection{A topological interpretation of the Burau representation}
\label{subsubsection:A topological interpretation of the Burau
representation} We have seen that the Burau representation is
deeply connected with the topology of closed braids.  We shall now
give a fully topological definition of the Burau representation.
We outline here ideas given in \cite{Bigelow1999}; for full
details see that or Turaev's excellent survey article
\cite{Turaev2000}, among others.

Let $D_n$ be an $n$-punctured disk, which we will think of as a
disk $D$ with $n$ distinguished points $q_1, \ldots, q_n$. Choose
a point $d_0 \in
\partial D_n$ to serve as the basepoint. Then $\pi_1
(D_n,d_0)$ is free on $n$ generators which can be represented by
loops $x_i$ based at $d_0$ travelling counterclockwise about the
puncture $q_i$ for $i = 1, \ldots, n$. We define a surjective map
$\epsilon: \pi_1 (D_n,d_0) \rightarrow \ints$ as follows. Let
$\gamma = x_{i_1}^{n_1} \cdots x_{i_r}^{n_r} \in \pi_1 (D_n,d_0)$.
Then we define the exponent sum $\epsilon (\gamma) = \sum_{i=1}^r
n_i$. The integer $\epsilon (\gamma)$ can be interpreted as the
\underline{total algebraic winding number} of $\gamma$ about the
punctures $\{q_i \}$, i.e., the sum over all $i = 1, \ldots, n$ of
the winding number of $\gamma$ about $q_i$. Now there is a regular
covering space $\tilde{D_n}$ of $D_n$ corresponding to the kernel
of the map $\epsilon$. Since $\epsilon: \pi_1 (D_n,d_0)
\rightarrow \ints$ is surjective, the group of covering
transformations $\textrm{Aut} (\tilde{D_n}) \cong \ints$.  Let
$\tau$ be a generator of $\ints$, and let $\Lambda = \ints
[t,t^{-1}]$. Then $H_1 (\tilde{D_n})$ inherits a $\Lambda$-module
structure from the action of the covering transformations: we
simply set $t \cdot \gamma = \tau_*(\gamma)$, where $\gamma \in
H_1 (\tilde{D_n})$ and $\tau_*$ denotes the induced action on
homology.  We note that as a $\Lambda$-module, $H_1 (\tilde{D_n})$
is free of rank $n - 1$.

In what follows it will be convenient to think of ${\bf B}_n$ as in
$\S$\ref{subsection:bBn and bPn as mapping class groups}, i.e., as
the mapping class group of $D_n$ where $\partial D_n$ is fixed
pointwise, while the punctures may be permuted.  We shall abuse
terminology by not distinguishing between a mapping class and a
diffeomorphism which represents it. Therefore we think of any $X
\in {\bf B}_n$ as a map $X: D_n \rightarrow D_n$. Then $X$ lifts
uniquely to a map $\tilde{X}:\tilde{D_n} \rightarrow \tilde{D_n}$
which fixes the fiber over the basepoint $d_0$ pointwise.
Furthermore, $\tilde{X}$ induces a $\Lambda$-module automorphism
$\tilde{X}_*$ of $H_1(\tilde{D_n})$.   Since $H_1 (\tilde{D_n})$
is a free $\Lambda$-module of rank $n - 1$, we can now define a
map ${\bf B}_n \rightarrow \textrm{GL}_{n-1}(\Lambda)$ by $X \mapsto
\tilde{X}_*$. This map turns out to be equivalent to the reduced
Burau representation defined previously (see \cite{Formanek94} for
a classification of linear representations of the braid group of
degree at most $n-1$).

The main idea of Stephen Bigelow's proof of the non-faithfulness
of the Burau representation in the case $n = 5$ is contained in
the following theorem:

\begin{theorem}{\rm \cite{Bigelow1999}}
\label{theorem:burau5} For $n \geq 3$, the Burau representation $\bar{\rho}: {\bf B}_n
\rightarrow \textnormal{GL}_{n-1}(\Lambda)$ is not faithful if and
only if there exist arcs $\alpha, \beta$ embedded in $D_n$
satisfying:
\begin{enumerate}
\item $\partial\alpha = \{q_1, q_2\}$ and $\partial\beta = \{d_0, q_3 \}$ or $\{q_3, q_4 \}$.
\item $\alpha$ intersects $\beta$ nontrivially  (more precisely,
there exists no isotopy rel endpoints which carries $\alpha$ off
$\beta$)
\item For some choice of lifts $\tilde{\alpha}, \tilde{\beta}$, we have $\sum_{k \in \ints}(t^k \tilde{\alpha}, \tilde{\beta}) t^k =
0$, where $(x,y)$ denotes the algebraic intersection number of two
(oriented) arcs in $\tilde{D_n}$.
\end{enumerate}
\end{theorem}
We note that the case $\partial\beta =\{d_0, q_3 \}$ follows from
Theorem~1.5 of \cite{Lo-Pa}.  Bigelow has produced an explicit
example of arcs $\alpha $ and $\beta$ satisfying the criteria of
Theorem~\ref{theorem:burau5} in the case $n = 5$ (see p. 402 of
\cite{Bigelow1999}).  It follows that the Burau representation of
${\bf B}_n$ is not faithful for $n \geq 5$.  It has been known for
many years that it {\em is} faithful for $n=3$.

\begin{openproblem}
{\rm Is the Burau representation of ${\bf B}_4$ faithful? $\clubsuit$
}
\end{openproblem}

In the next section, we will explore another representation of
${\bf B}_n$, with a topological interpretation analogous to that of the
Burau representation.

\subsection{The Lawrence-Krammer representation}
\label{subsection:the Lawrence-Krammer representation}

In 1990, Ruth Lawrence introduced a family of representations of
${\bf B}_n$  corresponding to Young diagrams with two rows which arise
out of a topological construction of representations of the Hecke
algebras $H_n(t)$ \cite{Lawrence1990}. Later, Krammer gave an
entirely algebraic definition of one of these representations:
$$\lambda: {\bf B}_n \rightarrow \textrm{GL}_r(A)$$
where $A = \ints[t^{\pm 1}, q^{\pm 1}]$ and showed that it was
faithful for $n = 4$ \cite{Krammer1999}.
The representation $\lambda$ has become known as the
Lawrence-Krammer representation. Shortly after Krammer's result
appeared, Bigelow was able to use topological methods to show that
$\lambda$ is faithful for all $n$ \cite{Bigelow2001}. Krammer
later gave an algebraic proof of the same result
\cite{Krammer2002}.  Therefore we now have:

\begin{theorem} {\rm \cite{Bigelow2001}, \cite{Krammer2002} }
\label{theorem:linearity} The map $\lambda$ is a faithful
representation of ${\bf B}_n$, and hence ${\bf B}_n$ is a linear group for
all $n$. \end{theorem} We observe that Krammer's result in
$\cite{Krammer1999}$ also implies that $\textrm{Aut}(F_2)$, the
group of automorphisms of a free group of rank $2$, is linear.

We shall not present here a full proof of
Theorem~\ref{theorem:linearity}; our aim is to define the
representation $\lambda$, to present the key ingredients in
Bigelow's proof: `forks', `noodles', and the pairing between them,
and to describe Bigelow's characterization of the kernel of
$\lambda$ in terms of these ingredients.

As mentioned previously, Bigelow's topological definition of
$\lambda$ is somewhat analogous to the topological definition of
the Burau representation given in the previous section.  This
time, we begin the construction with a certain configuration space
of a punctured disk rather than the punctured disk itself.  In the
notation of $\S$\ref{subsection:bBn and bPn via configuration
spaces}, we let $C = \cC_{0,2}(D_n)$, where $D_n$ is the
$n$-punctured disk described in the previous section. Elements of
$C$ will be denoted $\{z_1, z_2\} = \{z_2, z_1\}$, where $z_1 \neq
z_2 \in D_n$. Now we choose two distinct points $d_1, d_2 \in
\partial D_n$, and let $c_0 = \{d_1, d_2\}$ be the basepoint of
$C$.  We are now in a position to define a map $\phi: \pi_1(C,
c_0) \rightarrow \ints \times \ints$ as follows. Let $[\gamma] \in
\pi_1 (C, c_0)$. Then the loop $\gamma(s)$ in $C$ can be expressed
as $\gamma(s) = \{ \gamma_1(s), \gamma_2(s) \}$, where $\gamma_i$
is an arc in $D_n$ for $i = 1,2$.

We define two integers $a,b$ as follows. \begin{eqnarray*} a &=&
\frac{1}{2 \pi i} \sum_{j = 1}^n \int_{\gamma_1} \frac{dz}{z -
q_j} + \int_{\gamma_2}\frac{dz}{z - q_j}
\\
b &=& \frac{1}{\pi i} \int_{\gamma_1 - \gamma_2} \frac{dz}{z}
\end{eqnarray*}
Since $\gamma$ is a loop in $C$, we either have that both
$\gamma_1$ and $\gamma_2$ are both loops, i.e., $\gamma_1(0) =
\gamma_1(1)$ and $\gamma_2(0) = \gamma_2(1)$, or else, as Bigelow
puts it, the $\gamma_i$ are arcs which `switch places', i.e.,
$\gamma_1(1) = \gamma_2(0)$ and $\gamma_2(1) = \gamma_1(0)$.  In
the first case, Bigelow observes that we can interpret $a$ as the
sum of the total algebraic winding numbers of $\gamma_1$ and
$\gamma_2$ about the punctures $q_i$ and $b$ as twice the winding
number of $\gamma_1$ and $\gamma_2$ about each other.  In the
second case, the product $\gamma_1 \gamma_2$ is a loop, and so we
still have a nice interpretation of $a$ as the total algebraic
winding number of $\gamma_1 \gamma_2$ about the punctures $q_i$.
Bigelow also points out that in this case the condition $\gamma_1
- \gamma_2)(1) = -(\gamma_1 - \gamma_2)(0)$ implies that $b$ is an
odd integer.   It is worth noting that Turaev \cite{Turaev2000}
interprets the integer $b$ in a different way (in both cases). The
composition of the map $\pi_1(C) \rightarrow S^1$ defined by
$$
\gamma \mapsto \frac{\gamma_1 (s) - \gamma_2(s)}{|\gamma_1 (s) -
\gamma_2(s)|}
$$
with the usual projection from $S^1$ onto $\reals P^1$ sends our
loop $\gamma$ to a loop in $\reals P^1$.  Let $\bar{\gamma}$
denote the homology class of this loop in $H_1 (\reals P^1)$.  Up
to a choice of generator $u$ for $H_1 (\reals P^1) \cong \ints$,
we have $\bar{\gamma} = u^b$.  We define the map $\phi$ by setting
$\phi (\gamma) = q^a t^b \in \ints \times \ints$, thought of as
the free abelian group with basis $ \{ q,t \}$.

Now we proceed more or less as in the case of the Burau
representation.  Let $\tilde{C}$ be the regular cover of $C$ which
corresponds to the kernel of $\phi$ in $\pi_1 (C, c_0)$.  Let $A$
denote the ring $\ints [t^{\pm 1}, q^{\pm 1}]$. The homology
groups of $\tilde{C}$ naturally inherit an $A$-module structure
via the action of $\textrm{Aut} (\tilde{C})$.

It is clear that any homeomorphism $X: D_n \rightarrow D_n$
induces a homeomorphism $X' : C \rightarrow C$ defined by
$X'( \{z_1,z_2 \} = \{ X(z_1), X(z_2) \}$.  Thus
$X'$ necessarily fixes the basepoint $c_0$. If $X'_*$
denotes the induced action on $\pi_1 (C)$, then it is an easy
exercise to check that $X'_*$ preserves the values of the
integers $a$ and $b$ defined above.  Then $X'$ lifts uniquely
to a map $\tilde{X}: \tilde{C} \rightarrow \tilde{C}$ such
that the fiber over the basepoint $c_0$ is fixed pointwise.
Moreover, since $\tilde{X}$ commutes with the action of
$\textrm{Aut} (\tilde{C})$, we have that the induced action on
$H_2(\tilde{C})$, denoted $\tilde{X}_*$, is an $A$-module
automorphism.  (Note that $\tilde{C}$ is a 4-manifold.) Now $H_2
(\tilde{C})$ is a free $A$-module of rank $r = \left(
\begin{array}{c} n\\2 \end{array} \right)$ (Theorem 4.1 of
\cite{Bigelow2001}) and thus we can think of $\tilde{X}_* \in A$.
We therefore define the \underline{Lawrence-Krammer
representation} as follows:
\begin{eqnarray*}
\lambda: {\bf B}_n &\rightarrow& \textrm{GL}_r(A)\\
X &\mapsto& \tilde{X}_*
\end{eqnarray*}
If we choose $q$ and $t$ to be algebraically independent in
$\complexes$, we obtain a faithful representation $\lambda: {\bf B}_n
\rightarrow GL_r (\complexes)$.

Recall that ${d_1, d_2} = c_0$ is the basepoint we have chosen for
the configuration space $C$. A \underline{fork} is a tree $F$
embedded in the disk $D$ such that
\begin{enumerate}
\item $F$ has four vertices: $d_1, z, q_i$ and  $q_j$,
\item the three edges of $F$ share $z$ as a common vertex,
\item $F \cap \partial D_n = d_1$, and
\item $F \cap \{ q_1, \ldots, q_n \} = \{ q_i, q_j \}$.
\end{enumerate}
We note that Krammer's algebraic definition of the representation
$\lambda$ depends on the induced action on an $A$-module generated
forks.

The edge of $F$ which contains the vertex $d_1$ is called the
\underline{handle} of $F$.  The union of the other two edges is
called the \underline{tine edge} of $F$, denoted $T(F)$, and may
contain punctures $q_k, k \neq i,j$.   For any fork $F$ we define
a \underline{parallel copy} $F'$ to be a copy of the tree $F$
embedded in $D$, with vertices $\{ d_2, z', q_i$ and $q_j$ such
that $d_1 \neq d_2 \in \partial D_n$, $z \neq z'$, and $F'$ is
isotopic to $F$ rel $\{ q_i, q_j \}$.    Given a fork $F$, we also
define a \underline{noodle} to be an arc $N$ properly embedded in
$D_n$ with $\partial N = \{ d_1, d_2 \}$ and oriented so that its
initial point is $d_1$ and its terminal point is $d_2$.

We next construct surfaces in $\tilde{C}$ associated to forks and
noodles with which we will define a certain pairing on forks and
noodles. For each fork $F$, we choose a parallel copy $F'$ and
first define a surface $S(F)$ as the set of all points of the form
$\{ x, y \} \in C$ with $x \in T(F)\setminus \{q_1, \ldots, q_n\}$
and $y \in T(F')\setminus \{q_1, \ldots, q_n\}$. Let $\zeta$
(resp. $\zeta '$) denote the handle of $F$ (resp. $F'$), oriented
from $d_1$ to $z$ (resp. $d_2$ to $z'$). Define $\tilde{\zeta}$ to
be the unique lift of $\{ \zeta(s), \zeta' (s) \} \in \tilde{C}$
with initial point $\tilde{c_0}$, a fixed point in the fiber over
$c_0$.  Now let $\tilde{S}(F) \subset \tilde{C}$ be the unique
lift of $S(F)$ which contains the terminal point of
$\tilde{\zeta}$.   Similarly, we can associate to each noodle $N$
a surface $\tilde{S}(N) \subset \tilde{C}$. Define $S(N) \subset
C$ as the set of all points of the form $\{ x, y \} \in C$ such
that $x, y \in N, x \neq y$, and let $\tilde{S} (N)$ be the unique
lift of $S(N)$ which contains $\tilde{c_0}$.

Given a noodle $N$ and a fork $F$, we define their
\underline{pairing}, denoted $ \langle N, F \rangle$, to be an
element of the ring $A$ as follows.  We can assume that (up to
isotopy rel endpoints) $T(F)$ and $N$ intersect transversely in a
finite number of points $\{z_1, \ldots, z_r \}$ and similarly that
$T(F')$ intersects $N$ transversely at $\{z_1',\ldots, z_r' \}$.
Then each pair $z_i, z_i'$ cobounds an arc in $N$ which lies
between $T(F)$ and $T(F')$.  Now for each $i,j = 1, \ldots, r$,
there is a unique monomial $m_{i,j} = q^{a_{i,j}} t^{b_{i,j}}$
such that $m_{i,j} \tilde{S}(N)$ intersects the surface
$\tilde{S}(F)$ at a point in the fiber over ${z_i, z_j'} \in C$.
Letting $\epsilon_{i,j}$ denote the sign of that intersection, we
define
$$
\langle N, F \rangle = \sum_{i=1}^r \sum_{j=1}^r \epsilon_{i,j}
q^{a_{i,j}} t^{b_{i,j}}.
$$
Bigelow supplies a method for explicit calculation of the pairing
for a given noodle and fork and shows that it is well defined on
isotopy classes of forks rel $q_i, q_j$.

The proof that $\lambda$ is faithful relies on two important
lemmas (the `Basic Lemma' and `Key Lemma' of \cite{Bigelow2001}),
which we present here.  As usual, we shall abuse notation and not
distinguish between a mapping class and specific representatives.

\begin{lemma}{\rm \cite{Bigelow2001}}\label{lemma:basic}
If $\lambda (X) = 1$, then $\langle N,F \rangle = \langle N,
X(F) \rangle$ for every noodle $N$ and every fork $F$.
\end{lemma}

\begin{lemma}{\rm \cite{Bigelow2001}}\label{lemma:key}
For a noodle $N$ and a fork $F$, the pairing $\langle N,F \rangle
= 0$ if and only if $T(F)$ can be isotoped off $N$ relative to
$q_i, q_j$.
\end{lemma}

\noindent Thus it is the pairing on noodles and forks which gives
the essential characterization of the kernel of $\lambda$.

We end this section by noting two facts about the
Lawrence-Krammer representation. First, Budney has shown that
$\lambda$ is unitary\footnote{Squier had previously shown that the
Burau representation is unitary \cite{Squier1984} .} for an
appropriate choice of $q$ and $t$ (still algebraically
independent) \cite{Budney2002}.  In other words, we have
$$
\lambda: {\bf B}_n \rightarrow U_r (\complexes)
$$
where $\textrm{U}_r = \{X \in \textrm{GL}_r (\complexes)|
\bar{X}^T = X^{-1} \}$. Thus the conjugacy class of $\lambda (X)
\in U_r (\complexes)$ is determined by its eigenvalues, and one
could hope for an efficient solution to the conjugacy problem for
${\bf B}_n$ by passing to the whole of $U_r (\complexes)$. However,
Budney has given examples of non-conjugate braids $X_1, X_2 \in
{\bf B}_n$ such that $\lambda (X_1)$ and $\lambda (X_2)$ are conjugate
in $U_r (\complexes)$ (see Section 4 of \cite{Budney2002}).  We
will return to these examples in a bit more detail in
$\S$\ref{subsection:other solutions to the word problem}.

Second, Matthew Zinno has shown a connection between the
Lawrence-Krammer representation and the Birman-Murakami-Wenzl
(BMW) algebra.  The BMW algebra is related to Kauffman's knot
polynomial and can be thought of as a deformation of the Brauer
algebra in the same way that the Hecke algebra can be thought of
as a deformation of the group algebra $\complexes \sgn$ (see
\cite{BirmanWenzl}, \cite{Murakami1987}).  Braid groups map
homomorphically into the BMW algebra, giving rise to irreducible
representations of ${\bf B}_n$. Zinno has identified a summand of the
BMW algebra which corresponds exactly to the Lawrence-Krammer
representation:

\begin{theorem} {\rm \cite{Zinno2000}}
\label{theorem:zinno}
The Lawrence-Krammer representation of ${\bf B}_n$ is equivalent to
the $(n-2) \times 1$ irreducible representation of the BMW
algebra.
\end{theorem}
It follows immediately from Theorem~\ref{theorem:zinno} that the
Lawrence-Krammer representation $\lambda$ is irreducible and also
that the regular representation of the BMW algebra is faithful.

\subsection{Representations of other mapping class groups}
\label{subsection:representations of more general mapping class groups}
The first theorem that we proved in this review article was Theorem \ref{theorem:isomorphism between bbn and MCG of punctured disc}, which asserted that the braid group  ${\bf B}_{n-1}$ has a natural interpretation as the mapping class group ${\cal M}_{0,1,n-1}$ of the $n$-times punctured disc. One wonders, then, whether the linearity of ${\bf B}_n$ extends to a more general statement about the linearity of other mapping class groups ${\cal M}_{g,b,n}$?

The manuscripts of M. Korkmaz \cite{Korkmaz} and of S. Bigelow and R. Budney \cite{Bigelow-Budney} exploit the very special connection between the braid group ${\bf B}_{n-1} = {\cal M}_{0,1,n-1}$ and the mapping class group ${\cal M}_{0,0,n}$, and between  the mapping class groups ${\cal M}_{0,0,2g+2}$ and ${\cal M}_{g,0,0}$, to produce failthful finite dimensional representations of   ${\cal M}_{0,0,n}$ for every $n$, also of ${\cal M}_{2,0,0}$, and finally of a particular subgroup of ${\cal M}_{g,0,0}$ for every $g\geq 3$.  In the case of ${\cal M}_{2,0,0}$, the basic fact used by both Korkmaz and Bigelow-Budney is a theorem proved by Birman and Hilden in 1973 (see \cite{Birman1974}), which asserts that the mapping class group ${\cal M}_{2,0,0}$ is a $\ints_2$ central extension of ${\cal M}_{0,0,6}$.  This theorem is special to genus 2, which is the only group among the mapping class groups  ${\cal M}_{g,0,0}$ which has a center.  For $g\geq 3$ the so-called `hyperelliptic involution' which generates the center of  ${\cal M}_{2,0,0}$ generalizes to an involution whose centralizer is a subgroup of infinite index in ${\cal M}_{g,0,0}$. The same circle of ideas yield faithful matrix representations of those subgroups.

We now explain how Bigelow and Budney used the Lawrence-Krammer
representation of the braid groups to obtain faithful matrix
representations of ${\cal M}_{0,0,n.}$. If one considers any
surface $S_{g,b,n}$ and caps one of the boundary components by a
disc, one obtains a geometrically induced \underline{disc-filling}
homomorphism $d_\star: {\cal M}_{g,b,n-1} \to {\cal M}_{g,b-1,n}$.
Its kernel is the Dehn twist about the distinguished boundary
component (see  $\S$2.8 of \cite{Ivanov});  a distinguished point in the interior of the new disc
becomes the new fixed point.
Applying these ideas to the braid group, one then sees that there
is a natural  homomorphism $d_\star: {\bf B}_{n-1} \to {\cal
M}_{0,0,n}$.   It turns out that image$(d_\star)$ is not the full
group ${\cal M}_{0,1,n}$, but the stabilizer of the new fixed
point,  and  kernel$(d_\star)$ is the infinite cyclic subgroup of
${\bf B}_{n-1}$ that is generated by the braid $C =
(\sigma_{n-1}\sigma_{n-2}\cdots\sigma_2\sigma_1)^n$, a full twist
of all of the $n$ braid strands.  Since the image $\lambda(C)$ in
the Lawrence-Krammer representation is a scalar matrix (the
diagonal entries are $q^{2(n-1)}t^2$ in the representation as it
is defined in \cite{Bigelow2001}), one obtains a faithful
representation of image$(d_\star)$ by rescaling the
Lawrence-Krammer matrices, setting $t^{-2} = q^{2(n-1)}$.  The
group image$(d_\star)$ is of finite index in ${\cal M}_{0,0,n}$,
which therefore is a linear group. The dimension of the explicit
representation of ${\cal M}_{0,0,n}$ constructed in
\cite{Bigelow-Budney} is  $n\left(\begin{array}{c} {n-1}\\2
\end{array} \right)^2$.  It leads to a related representation of
dimension 64 of ${\cal M}_{2,0,0}$.    In this regard we note
that, while Korkmaz uses the identical geometry, he uses less care
with regard to dimension, and his representation of ${\cal
M}_{2,0,0}$ has dimension $2^{10}3^55^3$, which is very much
bigger than 64.

\begin{openproblem} {\rm Is there a faithful finite dimensional matrix representation of the mapping class group
${\cal M}_{g,b,n}$ for any values of the triplet $(g,b,n)$ other
than $(0,1,n), (0,0,n)$, $(2,0,0), (1,0,0)$ and $(1,1,0)$?  A folklore conjecture is
that most of the mapping class groups are, in fact
linear.    See \cite{Brendle-Hamidi} for evidence in this regard.
New ideas seem to be needed to construct candidates. } $\clubsuit$
\end{openproblem}

In view of the very large dimension of the representation of ${\cal M}_{2,0,0}$ which we just discussed, we note that there is an interesting 5-dimensional representation of the same group which occurs as one of the summands in the Hecke algebra representation of ${\bf B}_n$, namely the one belonging to the Young diagram with 2 rows and 3 columns. It is discussed in \cite{Jones1987}. At this writing its kernel does not seem to be known, although it is known to be infinite.

\subsection{Additional representations of ${\bf B}_n$.}
\label{subsection:additional representations} We end our discussion of representations of the braid groups by describing a construction which yields infinitely many finite
dimensional representations of ${\bf B}_n$, each one over a ring
$\complexes[t_1,t_1^{-1}, t_2,t_2^{-1},\dots,t_k,t_k^{-1}]$, where $t_1,\dots,t_k$ are parameter,  for
some $k\geq 1$. The construction includes all the summands in the
Temperley-Lieb algebra and the Lawrence-Krammer representation
too, and in addition infinitely many presumably new faithful
representations of ${\bf B}_n$. It was first described in
\cite{BLM1994}, generalizing ideas in \cite{Lawrence1990}. It is
due to Moody, with details first worked out by
Long in \cite{Long1994}.

We will be interested in the braid group ${\bf B}_n$ on n-strands,
but to describe our construction  it will be convenient to regard
${\bf B}_n$ as a subgroup of ${\bf B}_{n+1}$. Number the strands
in the latter group as $0,1,\dots,n$. Let ${\bf B}_{1,n}\subset
{\bf B}_{n+1}$ be the subgroup of braids in ${\bf B}_{n+1}$ whose
associated permutation fixes the letter 0. Its relationship to
${\bf B}_n$ is given by the  group extension
\begin{equation}
\label{short exact sequence defining B(1,n)}
 1\to {\bf F}_n \to {\bf B}_{1,n}\to {\bf B}_n \to 1,
\end{equation}
where the homomorphism ${\bf B}_{1,n}\to {\bf B}_n$ is defined by
pulling out the zero$^{th}$ braid strand.  There is a cross
section which is defined by mapping ${\bf B}_n$ to the subgroup of
braids on strands $1,\dots,n$ in ${\bf B}_{1,n}$. Therefore we may
identify ${\bf B}_{1,n}$ with ${\bf F}_n\semi{\bf B}_n$.   The
semi-direct product structure arises when we regard ${\bf B}_n$ as
a subgroup of the automorphism group of a free group.  The action
of ${\bf B}_n$ on $F_n$ is well known and is given in
\cite{Artin1925}, also  \cite{F-vB1962}, and also in
\cite{Birman1974}. Thinking of $F_n$ as the fundamental group
$\pi_1(S_{0,1,n})$ of the $n$-times punctured disc, the action of
the elementary braid $\sigma_i$ is given explicitly by:
\begin{equation}
\label{equation:action of Bn on Fn} \sigma_i{\bf x_j}\sigma_i^{-1}
= \cases {{\bf x_{i+1}} & if  $j=i$;\cr
 {\bf x_{i+1}}^{-1}{\bf x_i}{\bf x_{i+1}} & if  j=i+1\cr
{\bf x_j} & otherwise}
\end{equation}
Since ${\bf B}_n$ is a subgroup of ${\bf B}_{n+1}$ we see In this way that the groups ${\bf B}_n,  \ {\bf F}_n$ and
also ${\bf F}_n\semi{\bf B}_n$ are all subgroups of ${\bf
B}_{n+1}$.

In order to describe the idea behind the construction we recall
the notion of homology or cohomology of a space with coefficients
in a flat vector bundle. Suppose that $X$ is a manifold and that
we are given a representation
 $\rho : \pi_1(X) \rightarrow GL(V) $.
This enables us to define a flat vector bundle $E_\rho$: Let
$\tilde{X}$ be the universal covering of $X$. The group $\pi_1(X)$
acts
 on $ \tilde{X} \times V $ by
$ g.(\tilde{x}, {\bf v}) = (g.\tilde{x}, \rho(g).{\bf v})$. Then
$E_\rho$ is the quotient of $ \tilde{X} \times V $ by this action.
 We now form the cohomology
groups of 1-forms with coefficients in $E_\rho$, denoting these by
$H^1(X ; \rho)$
 or $H^1_c(X ; \rho)$ for compactly supported cochains.
In order to get an action of the braid groups, we use
(\ref{equation:action of Bn on Fn}). The action gives a canonical
way of forming a split extension $F_n  \semi {\bf B}_n$. It turns
out that in order to get an action on the twisted cohomology group
 what is required is exactly
a representation of this split extension.  Since ${\bf B}_{1,n}$
is a subgroup of ${\bf B}_{n+1}$, any representation of the latter
will of course do the job, and that is why we think of ${\bf B}_n$
as a subgroup of ${\bf B}_{n+1}$.

\begin{theorem} {\rm \cite{BLM1994}}
\label{theorem:splitconstruct} Given a representation $\rho : F_n
\semi {\bf B}_n \rightarrow GL(V) $ we may construct another
representation $ \rho_t^+ : {\bf  B}_n \rightarrow H^1_c(
S_{0,1,n} ; \rho) $ where $t$ is a new parameter.  In particular,
given any representation $\rho : {\bf }B_{n+1} \rightarrow GL(V)$,
we may construct a representation $\rho_t^+ : {\bf B}_n
\rightarrow H^1_c( S_{0,1,n} ; \rho) $.
\end{theorem}

This works in exactly the way one might expect.
 The representation restricted
to the free factor gives rise to the local system on the punctured
disc and thus the twisted cohomology group and the
 compatibility condition provided
by the split extension structure gives the
 braid group action.

A comment is in order concerning Theorem
\ref{theorem:splitconstruct}.  Although the theorem is stated
abstractly, (\ref{equation:action of Bn on Fn}) gives a concrete
recipe which enables one to write down the description of
$\rho_t^+$ given $\rho$.  Moreover, as we have noted, any time
that we have a representation of ${\bf B}_{n+1}$ we also have one
of the semi-direct product, which we identify with the subgroup
${\bf B}_{1,n}$.

The theorem shows that given a $k$ parameter representation of
the braid group, the construction yields a $k+1$ parameter
representation, apparently in a nontrivial way. For example, if
one starts with the (zero parameter) trivial representation of
$F_n  \semi {\bf B}_n$, the theorem produces the Burau
representation, and starting with the Burau representation, the
construction produces the Lawrence-Krammer representation
\cite{Long1994}. However the role of this extra parameter is not
purely to add extra complication, it also adds extra structure.
For there is a natural notion of what it should mean for a
representation of a braid group to be unitary (See
\cite{Squier1984}, for example) and the results of Deligne-Mostow
and Kohno imply:

\begin{theorem} {\rm \cite{Long1994}}
\label{theorem:unitary} In the above notation, if $\rho$ is
unitary, then
 for generic values of s,
so is $\rho_s^+$.
\end{theorem}

\begin{openproblem}{\rm This problem is somewhat vague. It begins with a suggestion
that Long's construction be studied in greater detail, and goes on
to ask whether (a wild guess) all finite dimensional unitary
matrix representations of ${\bf B}_n$ arise in a manner which is
related to the construction of Theorem  \ref{theorem:splitconstruct}?}
$\clubsuit$
\end{openproblem}

\newpage

\section{The word and conjugacy  problems in the braid groups}
\label{section:the word and conjugacy problems in the braid groups}
Let $G$ be any finitely generated group.  Fix a set of generators for $G$.
  A \underline{word} $\omega$ is a word in the given set of generators and their inverses.
 The element of $G$ that it represents will be denoted $[\omega]$, and its conjugacy class will be denoted $\{\omega\}$.   We consider two problems. The \underline{word problem} begins with $\omega,\omega^\prime\in G$, and asks for an algorithm that will decide whether $[\omega] = [\omega^\prime]$ ?  The  \underline{conjugacy problem} asks  for an algorithm to decide whether $\{\omega\} = \{\omega^\prime\}$, i.e. whether there exists $\alpha$ such that $\{\omega^\prime\} = \{\alpha^{-1}\omega \alpha\}$ ? A sharper version asks for a procedure for finding $\alpha$, if it exists.  In both cases we are interested in the complexity of the algorithm, and ask whether it is polynomial in either braid index $n$ or word length $|\omega|$ or both?
The word problem and the conjugacy are two of the three classical
`decision problems' first posed by Max Dehn \cite{Dehn1911}.  In
the case we consider here of $G = {\bf B}_n$, the set of generators we
choose will be those of either the classical or the new
presentation.  Artin gave the first solution to the word problem
for the braid group in 1925 \cite{Artin1925}, and a number of
fundamentally different solutions to the word problem in ${\bf
B}_n$ exist today, with at least two of them being polynomial in
both $n$ and $|\omega|$.  We focus here on an approach due to
Frank Garside \cite{Garside} and improved on by a number of
others, and briefly describe other methods at the end of this
section.

On the other hand, we know of only one definitive solution to the conjugacy problem, namely the combinatorial solution that was discovered by Garside.   He found a finite set of conjugates of an arbitrary element $[\omega]\in {\bf B}_n$,  the `Summit Set' ${\bf S}_\omega$ of $\omega$,  with the properties that if $\{\omega\} = \{\beta\}$, then ${\bf S}_\omega = {\bf S}_\beta$, whereas if $\{\omega\} \not= \{\beta\}$, then ${\bf S}_\omega \cap {\bf S}_\beta = \emptyset$.  If ${\bf S}_\omega = {\bf S}_\beta$, his methods also find $\alpha$ such that $\beta = \alpha\omega\alpha^{-1}$.  While his algorithm has been improved in major ways over the years, at this writing it is exponential in both $n$ and $|\omega|$.  The principle difficulties may be explained in the following way:

\bi
\item [($\star$)]  The entire set ${\bf S}_\omega$ and a single element in ${\bf S}_\beta$ must be computed to decide whether $\omega$ and $\beta$ are or are not conjugate.   While the calculation of a single element can now be done rapidly, in the general case ${\bf S}_\omega$ has unpredictable size.  The combinatorics that determine the size of ${\bf S}_\omega$ are particularly subtle and difficult to understand, and at this writing only partial progress has been made, in spite of much effort by experts.
\ei

The improvements that have been made over the years have included
the replacement of ${\bf S}_\omega$ by a proper subset $S_\omega,$
the `Super Summit Set',  via the work of ElRifai and Morton in
\cite{El-M1994}.   However, while  $S_\omega$ is very much smaller
than ${\bf S}_\omega$,  the principle difficulty $(\star)$ is
unchanged.  Very recently Gebhardt  found a still smaller subset,
the `Ultra Summit Set $U_\omega$', to replace $S_\omega$
\cite{Gebhardt}. The difficulty is unfortunately the same as it
was for the summit set and the super summit set, but can be made
more specific: the thing that one needs to understand is the
number of `closed cycling orbits' in $U_\omega$ and the length of
each orbit, and how distinct orbits are related.   Sang Jin Lee
\cite{Lee-2004} has given examples to show that the number of
orbits and their size can be arbitrarily large, and that it is in
no way clear how distinct orbits are related, except in special
cases. More work remains to be done.

It was shown in Theorem \ref{theorem:isomorphism between bbn and
MCG of punctured disc}, proved earlier in this article, that there
is a faithful action of the braid group ${\bf B}_n$ on the
$n$-times punctured disc.  Investigating this action (in the more
general setting of the action of the mapping class group of a
2-manifold on the 2-manifold),  Thurston proposed in the early
1980's a very different approach to the conjugacy problem which is
based upon the dynamical properties of that action.  These ideas
were investigated in Lee Mosher's PhD thesis \cite{Mosher-1983}.
The curious fact is that while the dynamical picture quickly
yields several very interesting class invariants not easily seen
in the combinatorial picture, namely a graph (known as a `train
track' ) which is embedded in the surface and is invariant under
the action of properly chosen representatives of $\{\omega\}$, and
a real number $\lambda$, however  they are not enough to determine
conjugacy.    Generic elements in the mapping class group of a
surface have an invariant train track, but it is not unique. It is
known that there are finitely many possible train tracks
associated to a given element, but their enumeration remains an
unsolved problem, and that is what gets in the way of a nice
solution to the conjugacy problem.  In fact, when one begins to
understand the details, the entire picture suggests difficulties
very much like those in $(\star)$ again.

We have just one more remark of a general nature, before we
proceed to review all these matters in more detail.  One of the
reasons that the Garside approach to the conjugacy problem is
interesting is because the techniques that were developed for the
braid groups revealed unexpected structures in ${\bf B}_n$ that
generalize to related unexpected structures in Artin groups (see $\S$\ref{subsection:Artin groups}), and
in the less well-known class of `Garside groups'.  See $\S$ \ref{subsection:generalization:from bBn to Garside groups}.
There has been major activity in recent years regarding these
structures, and there are also many open problems.  The same can
also be said for the Thurston approach, as it too generalizes from
an action of ${\bf B}_n$ on the n-times punctured disc  to actions
of mapping class groups on curves on surfaces, discussed in $\S$\ref{subsection:the conjugacy problem:the dynamic approach}.   All of these matters, and other related ones,  will be discussed below.

\subsection{The Garside approach, as improved over the years}
\label{subsection:Garside} We assume, initially,  that the group
${\bf B}_n$ is defined by the classical presentation
(\ref{equation:classical presentation}), with generators
$\sigma_1,\dots,\sigma_{n-1}$.   In this subsection we describe
the solution to the word and conjugacy problem which was
discovered by Garside in 1968 (see \cite{Garside}) and
subsequently sharpened and expanded, in many ways, through the
contributions of others, in particular Thurston \cite{Ep1992},
ElRifai and Morton \cite{El-M1994}, Birman, Ko and Lee
(\cite{BKL1998} and \cite{BKL2001}), Franco and Gonzales-Meneses
\cite{F-GM} and most recently Gebhardt \cite{Gebhardt}.  Other
relevant papers are \cite{GM2003} and \cite{GM-W2003}.

In (G1)-(G6) below we describe  a constructive method for finding  a normal form for words. In (G7)-(G10) we describe how to find a unique finite set of words in normal form that characterizes $\{\omega\}$.
\be
\item[{\bf (G1)}]  Elements of ${\bf B}_n$ which can be represented by braid words which only involve positive powers of the $\sigma_i$ are called \underline{positive braids}.  A key fact which was observed by Garside in \cite{Garside} is that the presentation (\ref{equation:classical presentation}) not only defines the group ${\bf B}_n$, it also defines a monoid ${\bf B}_n^+$. He then went on to prove that this monoid embeds in the obvious natural way in ${\bf B}_n$, in the strong sense that two positive words in ${\bf B}_n$ define the same element of ${\bf B}_n$ if and only if they also define the same element in the monoid ${\bf B}_n^+$.  The monoid of positive braids is particularly useful in studying ${\bf B}_n$ because the defining relations all preserve word length, so that the number of candidates for a positive word which represents a positive braid is finite.  However, all this is of little use unless we can show that the monoid of positive braids is more than just a small and very special subset of ${\bf B}_n$.  In (G3) below we show that this is indeed the case.

\begin{figure}[htpb!]
\centerline{\includegraphics[scale=.6,  bb=72 334 503 513] {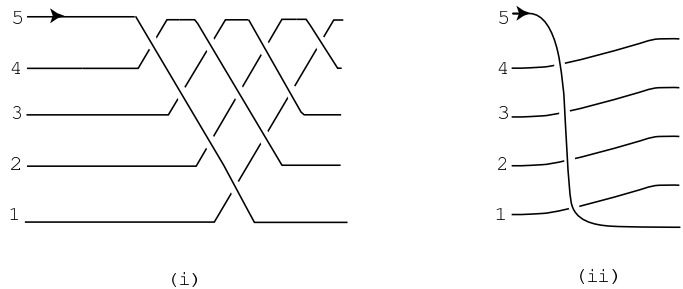}}
%\centerline{\includegraphics[scale=.6] {Delta.eps}}
\caption{(i) The 5-braid $\Delta_5$. (ii) The 5-braid $\delta_5$. }
\label{figure:Delta}
\end{figure}

\item [{\bf (G2)}]  A special role is played by the \underline{Garside braid} $\Delta$ in ${\bf B}_n$.  It is a positive half-twist of all the braid strands, and is defined by the braid word
$$\Delta_n = (\sigma_{n-1}\sigma_{n-2}\cdots\sigma_1)(\sigma_{n-1}\sigma_{n-2}\cdots\sigma_2)\cdots(\sigma_{n-1}\sigma_{n-2})(\sigma_{n-1}).$$
It is illustrated in sketch (i) of Figure \ref{figure:Delta}, for $n=5$.
The  square of $\Delta_n$ (a full twist of all the braid strands) generates the infinite cyclic center of ${\bf B}_n$.  The inner automorphism $\tau:{\bf B}_n\to{\bf B}_n$ which is defined by $\tau(X) = \Delta^{-1}X\Delta$  is the symmetry which sends each $\sigma_i$ to $\sigma_{n-i}$.  Of course $\Delta$ itself is invariant under this symmetry.
Later we will show that almost all of the combinatorics which we are describing hold equally well for the new presentation (\ref{equation:new presentation}).  The role of $\Delta$ in the classical presentation is replaced by that of the braid $\delta=\sigma_{n-1,n}\cdots\sigma_{2,3}\sigma_{1,2}$ which is depicted in Figure \ref{figure:Delta}(ii). See $\S$\ref{subsection:the new presentation} below.

\item [{\bf (G3)}]  The Garside braid $\Delta$ is very rich in elementary braid transformations. In particular, for each $i=1,\dots,n-1$   there exist (non-unique) positive braids $L_i$ and $R_i$ such that $\Delta = L_i \sigma_i = \sigma_i R_i$.   This implies that $\sigma_i^{-1} = \Delta^{-1} L_i = R_i \Delta^{-1}$ for each $i=1,\dots,n-1$, so that an arbitrary word $\omega= \sigma_{\mu_1}^{\epsilon_1}\cdots\sigma_{\mu_r}^{\epsilon_r}, \ \ \ \epsilon_q = \pm 1,$ can be converted to a word which uses only positive braid generators, at the expense of inserting arbitrary powers of $\Delta^{-1}$.  Since, by (G2) we have $\beta\Delta^\epsilon = \Delta^\epsilon\tau(\beta)$ for every braid $\beta$, all powers of $\Delta^{-1}$ can be moved to the left (or right).  It follows that $[\omega]$  is also defined by a word of the form $\Delta^rT$ where $T$ is positive.   This representation is non-unique.

\item[{\bf (G4)}] To begin to find unique aspects, Garside observes that there is a maximum value, say $i$, of $r$ with the property that the braid defined by $\omega$ is also represented by a word $\Delta^iZ$ where $Z$ is positive and $i$ is maximal for all braids of this form. Note that this means that $Z$ cannot be written in the form $Z_1\Delta Z_2$ for any positive words $Z_1,Z_2$, because if it could, then it would be possible, by (G2), to rewrite $\Delta^iZ$ as  $\Delta^{i+1} \tau(Z_1)Z_2$, with $\tau(Z_1)Z_2\in {\bf B}_n^+$, contradicting the maximality of $i$.  The integer $i$ is known as the \underline{infimum} of $[\omega]$, and written $inf(\omega)$. It is an invariant of $[\omega]$.
Garside's complete invariant of $[\omega]$ is $\Delta^iZ_0$, where among all positive words $Z_0,Z_1,\dots,Z_q$ that define the same element as $Z$, one chooses the unique word $Z_0$ whose subscript array (as a positive word in  $\sigma_1,\dots,\sigma_{n-1}$) is lexicographically minimal.  This is Garside's solution to the word problem in ${\bf B}_n$.  It is exponential in both $n$ and $|\omega|$ because there exist braids which admit very many of the elementary braid transformations in $(\ref{equation:classical presentation})$.   To make this very explicit, we note that by  \cite{Stanley1984} the number of positive words of letter length $(n)(n-1)/2$ in the standard braid generators  that represents $\Delta_n$ is the same as the number of standard Young tableaux of shape $(n,n-1,\dots,1)$. It is given by the hook length formula   as
$ ((n)(n-1)/2)) !  / (1^{n-1} 3^{n-1} 5^{n-2}\cdots (2n-3)).$   Calculating, we learn that $\Delta_5$ can be represented by 768 positive words, all of letter length 10, in the $\sigma_i's$, giving a good idea of the problem.    Moreover, this formula shows that the number of words increases exponentially with $n$. Of course the elements that are of interest to us, that is the positive braids defined by  words like $Z$ above, do not contain $\Delta$, however there is not reason why they cannot contain, for example, arbitrarily high powers of large subwords of $\Delta$.

\item[{\bf (G5)}]  Let $l$ and $r$ be positive words such that there is a factorization of $\Delta$ in the form $\Delta = lr$. We call $l$ a \underline{left divisor}  of $\Delta$ and $r$ a \underline{right divisor} of $\Delta$.  Let ${\cal P}$ (respectively ${\cal P'}$) be the collection of all left (resp. right) divisors of $\Delta$.  Note that we already used the fact that each of the generators $\sigma_i$ belongs to both ${\cal P}$ and ${\cal P'}$.  As it turns out, the sets ${\cal P}$ and ${\cal P'}$ coincide.  The set ${\cal P}$ will be seen to play an essential role in Garside's solution to the conjugacy problem, but before we discuss that we describe how to use it to get a very fast solution to the word problem.  The solution to the word problem was first discovered by Adyan \cite{Adyan} in 1984,  but Adyan's work was not well-known in the West. It was not referenced by either Thurston (see Chapter 5 of \cite{Ep1992}) or ElRifai and Morton  \cite{El-M1994}, who rediscovered Adyan's work 8-10 years later, and added to it.   Nevertheless, we use  \cite{El-M1994} and \cite{Ep1992} as our main source rather than \cite{Adyan},  because the point of view in those papers leads us more naturally to recent generalizations.

\item [{\bf (G6)}]  A key observation is that there are $n!$
braids in ${\cal P}$, and that these braids are in 1-1
correspondence with the permutations of their end-points, under
the correspondence defined by  sending each $\sigma_i$ to the
transposition $(i,i+1)$.  For this reason ${\cal P}$ is known as
the set of  \underline {permutation braids}.  This, combined with
the fact that the crossings are always positive,  makes it
possible to reconstruct any permutation braid from its
permutation, and so to obtain a unique element, even though its
representation as a word is highly non-unique.  Permutation braids
have a key property: after the braid is tightened, any two strands
cross at most once, positively.   To understand the importance of
this fact, we note that any positive 5-braid in which 2 strands
cross at most once, whose associated permutation is $(1,2,3,4,5)
\to (5,4,3,2,1)$ must be $\Delta_5$, a criterion which is an
enormous improvement over searching through the 768 distinct
positive words that the hook length formula showed us represent
$\Delta_5$, giving a good idea of the simplification that Adyan, Thurston
and ElRifai and Morton discovered.

The final step in finding a very rapid solution to the word
problem in ${\bf B}_n$ is a unique way to factorize
the braid $P$ in the partial normal form $\Delta^iP$ of (G4) above
as a unique product of finitely many permutation braids.   If the
strands in $P$ cross at most once, then $P\in \cal{P}$ and we may
choose any positive braid word, say $l_1$,  which has the same
permutation as $P$ as its representative.  If not, set $P =
l_1l_1^\star$ where $l_1$ is a positive braid of maximal length in
which two strands cross at most once and $l^\star$ is the rest of
$P$. If no two strands in $l_1^\star$  cross more than once, set
$l_2=l_1^\star$, and stop. If not, repeat, setting $l_1^\star=
l_2l_2^\star$, where $l_2$ has maximal length among all positive
braids whose strands cross at most once, and so forth to obtain $P
= l_1l_2\cdots l_s$ where each $l_i\in{\cal P}$ and each $l_i$ has
maximal length for all factorizations of $l_{i-1}^{-1}\cdots
l_1^{-1} L$ as $l_il_i^\star$ with $l_i\in{\cal P}$ and
$l_i^\star$ positive. This representation is unique, up to the
choices of the words which represent $l_1,\dots,l_s$.  As we
remarked earlier, each $l_i$ is determined uniquely as an element
in ${\bf B}_n^+$ by the permutation of its strands.

The factorization of an arbitrary braid word as $\Delta^iP =
\Delta^i l_1l_2\cdots l_s$ is the \underline{left greedy} normal
form, where the term `left greedy' suggests the fact that each
$l_i$ is a maximal permutation subbraid relative to the positive
braid to its right. The left greedy normal form solves the word
problem. The  integer $i$ was already defined to be $inf(\omega)$.  The integer $i + s$ is known as the
\underline{supremum} of $\omega$, written
$sup(\omega)$. The integer  $s$ is the
\underline{canonical length} $L(\omega)$.

\begin{remark} {\rm The Adyan-Thurston-ElRifai-Morton solution to the word problem is shown in \cite{BKL1998}  to be ${\cal O}(|\omega |^2 n \rm{log} n)$, where $|\omega|$ is  the word length of the initial representative of the braid $\omega$, as a word in the standard generators of ${\bf B}_n$.}
\end{remark}
\begin{example}\label{example:leftgreedy} {\rm We illustrate how to find the left greedy normal form for a positive braid via  an example.  Assume that we are given the positive braid
$$P = \sigma_1 \sigma_3 \sigma_2^2 \sigma_1 \sigma_3^2 \sigma_2 \sigma_3 \sigma_2.$$   The first step is to factor $P$  as a product of  permutation braids by working along the braid word from left to right  and inserting the next partition whenever two strands are about to cross a second time since the beginning of the current partition.  The reader may wish to draw a pictures to go with this example.  This gives the factorization:
$$(\sigma_1\sigma_3\sigma_2)(\sigma_2\sigma_1\sigma_3)(\sigma_3\sigma_2\sigma_3)(\sigma_2)$$
Next, start at the right end of the braid word,  and ask whether the elementary braid relations can be used to move crossings from one partition to the partition on its immediate left, without forcing an adjacent pair of strands in a factor to cross.  This process is repeated three times. The first time one applies the braid relations to factor 3, and then to adjacent letters in factors 2 and 3, to increase the length of factor 2 at the expense of decreasing the length of factor 3:
$$(\sigma_1\sigma_3\sigma_2)(\sigma_2\sigma_1\sigma_3\sigma_2)(\sigma_3\sigma_2)(\sigma_2.)$$
In fact one can push one more crossing from factor 3 to factor 2:
$$(\sigma_1\sigma_3\sigma_2)(\sigma_2\sigma_1\sigma_2\sigma_3\sigma_2)(\sigma_2)(\sigma_2),$$
Now it is possible to move a crossing from  factor 2 to factor 1, giving:
$$(\sigma_1\sigma_3\sigma_2\sigma_1)(\sigma_2\sigma_1\sigma_3\sigma_2)(\sigma_2)(\sigma_2),$$
No further changes are possible, we have achieved left greedy normal form. If it had happened that $P$ contained $\Delta_4$, our factorization would reveal it because $\Delta_4\in{\cal P}$.   Each canonical factor represents a braid in ${\cal P}\subset {\bf B}_4^+$ and so has length strictly less than $|\Delta_4| = 6$.  The words which represent the canonical factors are non-unique. Their associated permutations determine the canonical factors uniquely. }   $\spadesuit$  \end{example}

 \item[{\bf (G7)}]  We pass to the conjugacy problem, which builds on the solution just given to the word problem.   From now on a word $\omega$ be always be assumed to be in left greedy normal form.  Notice that the abelianizing map $e:{\bf B}_n \to \ints$ has  infinite cyclic image, and that $e(\omega)$ is the exponent sum of a representing word in the $\sigma_i$. Clearly $e(\omega)$ is an invariant of both word and conjugacy class. It follows that  there are only finitely many braids $\omega^\prime$ in the conjugacy class $\{\omega\}$ which have left greedy normal form  with $inf(\omega^\prime) >$ $inf(\omega)$, because  $e(\Delta) = (n)(n-1)/2$ and any increase in $inf(\omega)$ must come at the expense of a corresponding decrease in word length of the positive word that remains after all powers of $\Delta$ have been pushed to the left.  Let $Inf(\omega$) be the maximum  value of $inf (\omega')$  for all braids $\omega'$ in the conjugacy class $\{\omega\}$.   Assume that a representative $\omega_1$ of $\{\omega\}$ is given and that it has left greedy normal form $\Delta^I L_1L_2\dots L_S$, where $I = Inf(\omega)$ and $S$ is minimal for all braids in normal form which are conjugate to $\omega$.   Let $Sup(\omega)$  be the integer $I+S$.  Then $Inf(\omega)$ and $Sup(\omega)$ are class invariants of $\omega$.   The ElRifai-Morton \underline{super summit set} $S_\omega$ is the collection of all elements in left-greedy normal form which realize $Inf(\omega)$ and $Sup(\omega)$.   It is a major improvement over Garside's \underline{summit set} ${\bf S}(\omega)$, which is the larger set of all elements which realize $Inf(\omega)$ but not $Sup(\omega)$, with left greedy form replaced by the subscript ordering described in (G4) above.

One might wonder how the normal forms which we just described for
words and conjugacy classes relate to length functions on ${\bf
B}_n$.  In this regard, R. Charney has introduced a concept of
`geodesic length' in \cite{Charney}.  As above, let ${\cal P}$ be
the set of all permutation braids. Charney defines the
\underline{geodesic length} of a braid $\omega$ to be the smallest
integer $K = K(\omega)$ such that there is a word
$L_1^{\epsilon_1}L_2^{\epsilon_2}\cdots L_K^{\epsilon_K}$, where
each $\epsilon_i = \pm 1$, which represents the conjugacy class
$\{\omega\}$, with each $L_i\in {\cal P}$ and each $L_i^{-1}$ the
inverse of a word in ${\cal P}$.   If $\Delta^u L_1L_2\cdots L_s$
is an arbitrary element in the super summit set of $\omega$, it is
not difficult to show that the geodesic length of $\{\omega\}$ is
the maximum of the 3 integers $(s+u, -u, s)$.  We remark that
Charney's geodesic length is defined in  \cite{Charney} for all
Artin groups of finite type.  This and other ways in which the
Garside machinery generalizes to  other classes of groups,
including all Artin groups of finite type,  will be discussed
(much too briefly) in   $\S$\ref{subsection:generalization:from
bBn to Garside groups}.

\begin{openproblem}
\label{openproblem:Krammer's proof of faithfulness} {\rm In the
manuscript  \cite{Krammer2002} Krammer gives several proofs of the
faithfulness of the Lawrence-Krammer representation of ${\bf
B}_n$.  One of his very interesting proofs shows that if
$[\omega]\in{\bf B}_n$, then it is a very simple matter to read
Charney's geodesic length from the matrix representation of
$\omega$. Since the unique element of geodesic length zero is the
identity, it follows that the representation is faithful.  Our
suggestion for future work is to  investigate the Garside solution
to the conjugacy problem and its improvements (to be described
below) via the Lawrence-Krammer representation of ${\bf B}_n$.  We
mention this because we feel that this aspect of Krammer's work
has received very little attention.} $\clubsuit$
\end{openproblem}

\item [{\bf (G8)}] Here is a fast constructive procedure
for finding $Inf(\omega)$ and $Sup(\omega)$, due to ElRifai and
Morton \cite{El-M1994} and to Birman, Ko and Lee in \cite{BKL2001}.  Starting with any element in $\{\omega\}$,  define the
\underline{cycling} of $\omega=\Delta^il_1l_2\cdots l_k$ to be the
braid $c(\omega) = \tau^{-i}(l_1)\omega\tau(l_1)$ and the
\underline{decycling} of $\omega$ to be the braid
 $d(\omega) = l_k\omega (l_k)^{-1}$.  Note that  $c(\omega)$
 and $d(\omega) \in \{\omega\}$.  Putting $c(\omega)$ and
$d(\omega)$  into left greedy normal form, one obtains braids
which have at least $i$ powers of $\Delta$, and possibly more
because it can happen that after cycling the braid
$\Delta^il_1l_2\cdots l_k$ will change to one whose left greedy
normal form is $\Delta^{i^\prime}l_1'l_2\cdots l_k'{s^\prime}$
with $i^\prime\geq i$.  This would, of course, reduce $s$ if it
increases $i$.  In \cite{El-M1994} ElRifai and Morton  proved that
if $inf(\omega)$ is not maximal for the conjugacy class, then it
can be increased by repeated cycling.  Also, if $sup(\omega)$ is
not minimal, then it can be decreased by decycling.

They could not say, however, how many times one might have to
cycle or decycle before being sure that no further improvement was
possible.  The solution to that problem was found by Birman, Ko
and Lee in \cite{BKL2001}. It was shown that if $inf(\omega)$ is
not maximal for the conjugacy class if it will be increased after
fewer than $(n)(n-1)/2$ cyclings, and similarly if $sup(\omega)$
is not minimal it will be decreased after fewer than $(n-1)(n)/2$
decyclings. One then has a tool for increasing $inf(\omega)$  and
decreasing $sup(\omega)$, and also a test which tells,
definitively, when no further increase or decrease is possible.

\begin{remark} {\rm The fact that everything we do to compute Inf($\omega)$ also applies to the computation of Sup($\omega)$  is not surprising because ElRifai and Morton showed that $Sup(\omega)$ = -Inf$(\omega^{-1})$.}
\end{remark}

\item [{\bf (G9)}]  We now come to Gebhardt's very new work.  Following the steps given in (G8) above, one will have on hand the summit set $S(\omega)$, i.e. the set of all braids in the conjugacy class $\{\omega\}$ which have left greedy normal form $\Delta^I L_1 L_2 \cdots L_S$, where $I $= $Inf(\omega)$ is maximal for the class and $I+S$ = $Sup(\omega)$ is minimal for the class. This set is still very big.  Gebhardt's improvement is to show that it suffices to consider only the subset of braids which in a closed orbit under cycling.  This finite set of words is called the \underline{ultra summit set} $ U_\omega$. It has been proved by Volker Gebhardt \cite{Gebhardt} that $\{\omega\} = \{ \omega'\}$  if and only if $U_\omega $= $U_{\omega^\prime}$.  Note that Gebhardt does not need to use decycling, he proves that it suffices to consider the closed orbits under cycling.

\item[{\bf (G10)}]  To compute $U_\omega$, Gebhardt also shows that if $\omega_i, \omega_j\in$ $U_\omega$, then there is a finite chain $\omega_i = \omega_{i,1}\to\omega_{i,2}\to\cdots\to\omega_{i,q}=\omega_j$ of braids, with each $\omega_{i,j}$ in $U_\omega$ such that each $\omega_{i,t}$ is obtained from $\omega_{i,t-1}$ by conjugating by a single element in ${\cal P}$.  Thus the following steps suffice to compute $
U_\omega$, after one knows a single element $\rho\in U_\omega$: One first computes the conjugates of $\rho$ by the $n!$ elements in ${\cal P}$.  One then puts each into left greedy normal form, and discards any braid that either does not (a) realize $Inf$ and $Sup$, or (b) realizes $Inf$ and $Sup$ but is not in a closed orbit under cycling, or (c) realizes $Inf$ and $Sup$ and is in a closed orbit under cycling but is not a new element in $U_\omega$.  Ultimately, the list of elements so obtained closes to give  $U_\omega$.

Note that in doing this computation, one not only learns, for each $\omega_i,\omega_j\in$ $U_\omega$, that $\{\omega_j\} = \{\omega_j\}$, but one also computes an explicit element $\alpha$ such that $\omega_j = \alpha^{-1}\omega_i\alpha$.

An inefficient part of this computation is the constant need
to access the $n!-2$ non-trivial elements in ${\cal P}$.
A more efficient process is known (see (G11) below), but to describe it, we need new notions, which will be introduced after we discuss a wider class of groups, known as Garside groups.\ee

\subsection{Generalizations: from ${\bf B}_n$ to Garside groups}
\label{subsection:generalization:from bBn to Garside groups} The
first person to realize that the structure described in (G1)-(G10)
is not restricted to braids was Garside himself  \cite{Garside},
but his generalizations were limited to examples.  Soon after his
paper was published,  the ideas were shown to go through,
appropriately modified, in all \underline{finite type} Artin
groups, i.e., Artin groups whose associated Coxeter group is
finite, by Breiskorn and Saito \cite{B-S} and by Deligne
\cite{Deligne}. See \cite{Michel} for a survey (and reworking) of the results first proved in \cite{B-S} and \cite{Deligne}.  (See $\S$\ref{subsection:Artin groups} for
definitions of Artin and Coxeter groups in general.)  They used
explicit properties of finite reflection groups in their proof,
but Paris and Dehornoy were thinking more generally and defined a
broader class which they called `Garside groups'.  The class included
all finite type Artin groups.  Over the last
few years, several tentative definitions of the term `Garside
group', referring to various classes of groups that generalize the
braid group in this way, were proposed and appear in the
literature (see \cite{D-P}, \cite{Picantin2001}, for example)
before the one that we give below was agreed upon by many,
although the search for the most general class of such groups
continues and the definition of `Garside groups' is likely to
continue to be in flux for some time.   As we proceed through the
definitions of `Garside monoids', `Garside structures', and
`Garside groups', the reader can look to the classical
presentation of ${\bf B}_n$  for examples.

Given a finitely generated monoid $G^+$ with  identity $e$, we can
define a partial orders on its elements.  Let $a, b \in G^+$. We
say that $a \prec b$  if $a$ is a left divisor of $b$, i.e. there
exists $c\in M, \  c\not= e$, with $ac = b$. Also   $b\succ a$ if
there exists $c$ such that $b = ca$.  Caution: the two orderings
are really different, that is, $a\prec b$ does not imply that
$b\succ a$. An interesting example in the positive braid monoid
${\bf B}_n^+$ is the partial order  induced on the elements in the
sets ${\cal P}$  of left divisors of $\Delta_n$, which gives it
the structure of a lattice.   The reader who wishes to get a
feeling for the ordering might wish to construct the lattice in
the cases $\Delta_3$ and $\Delta_4$.  The lattice for $\Delta_n$
has $n!$ elements.

Now given $a,b\in G^+$ we can define in a natural way the (left)
greatest common divisor of $a$ and $b$, if it exists, written $d =
a\wedge b$, as follows: $d \prec a, d\prec b$ and if, for any $x$,
it happens that $x\prec a$ and $x\prec b$ then $x\prec d$.
Similarly, we define the (left)  least common multiple of $a,b$,
denoted $m = a\vee b$,  if $a\prec m, b \prec m$ and if for any
$x$ it happens that $a\prec x$ and $b\prec x$ then $m\prec x$.

It turns out that in the case of the braid monoid ${\bf B}_n^+$
the left partial ordering extends to a right-invariant ordering on
the full braid group ${\bf B}_n$, a matter which we will discuss
in $\S$\ref{subsection:other solutions to the word problem} below.

We continue with our description of the features of the monoid
${\bf B}_n^+ \subset {\bf B}_n$ which will lead us to define more
general monoids $G^+$ and their associated groups $G$.  An element
$x\in G^+$ is an \underline{atom} if $x\not= e$ and if $x$ has no
proper left or right divisors.  For example, the generators
$\sigma_1,\dots,\sigma_{n-1}$ are the atoms in the monoid $G^+$.
Note that  by (G3) above the atoms in ${\bf B}_n^+$ are left and
right divisors of $\Delta$, also they generate ${\bf B}_n^+$.

 While we have not had occasion to introduce a key property of the Garside braid $\Delta$ before this, we do so now:   As noted earlier, its sets ${\cal P}$ and ${\cal P'}$  of left and right divisors of $\Delta$ coincide.  This property was used in the proofs of the facts that we described in (G1)-(G10) above, and is part of a long story about symmetries in the braid group.

 A  monoid $G^+$ is  a \underline{Garside monoid}  if:
\be
\item $G^+$ is generated by its atoms,
\item For every $x\in G^+$ there exists an integer $l(x)>0$ such that $x$ cannot be written as a product of more than $l(x)$ atoms,
\item  $G^+$ is left and right cancellative, and every pair of elements in $G^+$ admits a left and also a right least common multiple and greatest common divisor,
\item There is an element $\Delta \in G^+$, the \underline{Garside element},  whose left divisors and right divisors coincide, also they form a lattice, also each generates $G^+$.
\ee All of the data just listed is called a \underline{Garside
structure}.  It follows (via the work of Ore, described in Volume
1 of the book  \cite{C-P1961}) that every Garside monoid embeds in
its group of fractions, which is defined to be a
\underline{Garside group}.  It is not difficult to prove that the
Garside element in a Garside monoid is the least common multiple
of its set of atoms.  An interesting property is that a Garside
monoid admits a presentation  $\langle S \ | \ R \rangle$, where
for every pair of generators $x,y$ there is a relation of the form
$x\cdots = y\cdots$  that prescribes how to complete $x$ and $y$
on the right in order to obtain equal elements.

Examples abound. While we shall see that the braid groups are
torsion-free, there are examples of Garside groups which have
torsion. Whereas the relations in the braid monoid all preserve
word length (which results in various finiteness aspects of the
algorithms that we described) there are examples of Garside
monoids in which this is not the case.  As noted above, every
Artin group of finite type is a Garside group.  Torus knot groups
are Garside groups, as are fundamental groups of complements of
complex lines through the origin \cite{D-P}.  See
\cite{Picantin2001} for additional examples.

Our reason for introducing Garside groups is that the solutions to
the word problem and conjugacy problem in ${\bf B}_n$ described in
$\S$\ref{subsection:Garside}, suitably modified, generalize to the
class of Garside groups (\cite{D-P},  \cite{Picantin2001}), as do
the simplifications of (G10) which we now describe.   These
improvements were first discovered by Franco and Gonzales-Meneses
\cite{F-GM}, and later improved by  Gebhardt \cite{Gebhardt}.

We return to the braid group, with the partial ordering of ${\bf B}_n$ on hand:
\be
\item [{\bf (G11)}]   We already learned that if we begin with an arbitrary braid $\omega$, then after a bounded number of  cyclings and decyclings we will obtain a braid in $\{\omega\}$ which realizes $Inf(\omega)$ and $Sup(\omega)$. Continuing to cycle, we will arrive at an element, say $\rho$,  in  $U_\omega$.  The remaining task is the computation of the full set $U_\omega$, and the method described in (G10) is inefficient.  The difficulty is that, starting with $\rho\in U_\omega$  one is forced to  compute its conjugates   by {\it all} the $n!$ elements in ${\cal P}$, even though many of those will turn out to either not realize $Inf(\omega)$ and/or $Sup(\omega)$, and so will be discarded, whereas others will turn out to be duplicates of ones already computed.  Moreover, this inefficient step is done repeatedly.

%%JOAN -- There was a typo below that I fixed.
The good news is that, following ideas first introduced by N.
Franco and J. Gonzalez-Meneses in \cite{F-GM}, Gebhardt proves  in
Theorem 1.17 of \cite{Gebhardt},  that  if  $\rho\in U_\omega,$
and if $a,b \in {\cal P}$, with $a^{-1}\rho a$ and $b^{-1}\rho b
\in  U_\omega$, then $c^{-1} \rho c \in U_\omega$, where $c =
a\wedge b$.    See \cite{F-GM}, and then \cite{Gebhardt},  for a
systematic procedure that allows one to use this fact to find all
the orbits in $U_\omega$ efficiently.  In this regard we remark
that the work of Gebhardt is very new.  The major open problem
that remains is to improve it to an algorithm  which will be
polynomial in the word length of $\omega$: \ee

\begin{openproblem}
\label{openproblem:structure of the USS} {\rm The bad news is
that, like the Garside's summit set and ElRifai and Morton's super
summit set, Gebhardt's ultra summit set also has the key
difficulty $(\star)$. More work remains to be done.  } $\clubsuit$
\end{openproblem}

\subsection{The new presentation and multiple Garside structures}
\label{subsection:the new presentation}

As it turns out,  essentially all of the structure that we just described also exists with respect to the second presentation of ${\bf B}_n$ which was given in $\S$\ref{section:introduction} of this paper.  This is the main result of \cite{BKL2001}.  The fundamental braid $\Delta_n$ is replaced by
$\delta_n = \sigma_{n-1,n} \sigma_{n-2,n-1}\cdots \sigma_{2,3} \sigma_{1,2}.$     It is proved in \cite{BKL2001} that the associated monoid ${\bf B}_n^+$ embeds in ${\bf B}_n$, and that all the results described in (G1)-(G11) above have curious (and very surprising) variations which hold in the new situation.

For example, the elements in the set of left divisors of
$\delta_n$ are in 1-1 correspondence with a set of permutations,
namely permutations which are products of non-interlacing
descending cycles.  The set of left divisors of $\delta_n$ turns
out to have  order equal to the $n^{th}$ Catalan number,
$(2n!)/(n!(n+1)!)$, whereas the left divisors of $\Delta$ have
order $n!$.  The fact that the order of the set of permutation
braids is smaller in the new presentation than in the classical
presentation had led to the hope by the authors of
\cite{BKL2001}, when they first discovered the new presentation,
that it would result in a polynomial algorithm  for the conjugacy
problem, but the Catalan numbers grow exponentially with index,
and once again $(\star)$ proved to be a fundamental obstacle.
Thus, while the new presentation is extremely interesting in its
own right, and does lead to faster word and conjugacy algorithms,
the improvement in that regard does not address the fundamental
underlying difficulties.

\begin{openproblem}
\label{openproblem:other Garside monoids in Bn?}
{\rm Curiously, it appears very likely that the classical and new
presentations of ${\bf B}_n$ are the only positive presentations
of ${\bf B}_n$ in which the Garside structure exists, although
that has not been proved at this time, and is an interesting open
problem.  For partial results in this direction, see
\cite{K-H2002}.  } $\clubsuit$
\end{openproblem}

There is a different aspect of the dual presentations which we mention now, which involves a small detour.  Before we can explain it, recall that one of our earliest definitions of the braid group ${\bf P}_n$, given in $\S$\ref{subsection:bBn and bPn via configuration spaces} was as the fundamental group $\pi_1(\cC_{0,\hat{n}}, {\vec p})$,  of the space formed from $\complexes^n$ by deleting the hyperplanes along which two or more complex coordinates coincide.  Of course this gives a natural cell decomposition for  $\cC_{0,\hat{n}}$ as a union of (open) cells of real dimension $2n$. The braid group ${\bf B}_n$, as we defined it in $\S$\ref{subsection:bBn and bPn via configuration spaces},  is the fundamental group of the quotient $\cC_{0,n} = \cC_{0,\hat{n}}/\Sigma_n$, where the symmetric group $\Sigma_n$ acts on $\cC_{0,\hat{n}}$ by permuting coordinates.   In the interesting manuscript \cite{Fox-Neuwirth}, Fox and Neuwirth used this natural cell decomposition of $\cC_{0,\hat{n}}$ to find a presentation for ${\bf B}_n$, arriving at the classical presentation (\ref{equation:classical presentation}) for the braid group ${\bf B}_n$.  See Section C of Chapter 10 of \cite{Burde-Zieschang} for a succinct presentation of the results in \cite{Fox-Neuwirth}.   Fox and Neuwirth also prove that $\cC_{0,\hat{n}}$ is aspherical, and use this to give the first proof that ${\bf B}_n$ is torsion-free. See $\S$\ref{subsection:braid groups are torsion-free} below for a different and very easily understood proof of that same fact.

%%JOAN -- I thought we ought to give the reference to the subsection in which the lattice discussion takes place that you refer to in the next paragraph.  I added it and reworded a little bit to incorporate this addition.

Around the same time that the first author, together with Ko and
Lee, wrote the paper \cite{BKL1998},  which introduced the dual
presentation of ${\bf B}_n$, Thomas Brady was thinking about other
complexes which, like the one just described might serve as a
$K(\pi,1)$ for the group ${\bf B}_n$.  See \cite{Brady}, which
describes the construction of a finite $CW$-complex $K_n$ of
dimension $n-1$ which is homotopy equivalent to $C_{0,n}$, and
which is defined combinatorially, using the partial ordering on
$\Sigma_n$, as described in
$\S$\ref{subsection:generalization:from bBn to Garside groups} in
the discussion of the lattice of simple words in ${\bf B}_n$.
Brady then used his complex to determine a presentation for ${\bf
B}_n$, arriving at the dual presentation (\ref{equation:new
presentation}).

%%JOAN -- with the addition of the new first sentence of this paragraph, we now have "As it turned out" twice in the following paragraph.  We could just delete one of them.

As it turned out, the new presentation was important for other
reasons too. The braid groups inherit a second Garside structure
from the new presentation (\ref{equation:new presentation}).  The
fact that the same was true for many other Garside groups played a
major role in the discovery of the appropriate definitions.   The
Garside structure on the braid group arising from the new
presentation ends up revealing even more structure: it is dual to
that coming from the standard presentation, in the sense of an
action on a complex and a dual complex.  As it turned out, all
finite type Artin groups of finite type also have dual 'Garside
structures, as proved by Bessis \cite{Bessis} (this reference also
contains the details of this dual structure). See also
\cite{Brady-Watt}. We also refer the reader to \cite{Picantin2001}
for explicit presentations for the dual monoids associated to the
Artin groups of finite type, and also for an interesting table
that gives the number of simple elements defined by the left
divisors  of the Garside element in the classical and dual monoids
for the finite type Artin groups.

\begin{openproblem}
\label{openproblem:multiple Garside structures on Garside groups}
{\rm This one is a very big set of problems. It is not known
whether all Garside groups have dual presentations, in fact, it is
also not known how many distinct Garside structures a given
Garside group may have. } $\clubsuit$
\end{openproblem}

\subsection{Artin monoids and their groups}
\label{subsection:Artin groups}  There is another class of groups
which is intimately related to the braid group and its associated
monoid,  but it is much less well understood than the Garside
groups.  Let $S = \{ u,\dots,\dots,t, \dots,v\}$ be a finite set.
A \underline{Coxeter graph} $\Gamma$ over  $S$ is a graph whose
vertices are in 1-1 correspondence with the elements of $S$. There
are no edges joining a vertex $s$ to itself.  There may or may not
be an edge joining a vertex $s$ to a vertex $t\not= s$.  Each pair
of vertices $(s,t) = (t,s)$ is labeled by a non-negative integer
$m(s,t)$.  There are 2 types of labels: The label $m(s,t)$ is $2$
if $\Gamma$ does not have an edge that joins $s$ and $t$, and it
is $\in \{3,4,\dots,\infty\}$ if there is an edge joining $s$ and
$t$.

The \underline{Artin group} $A(\Gamma)$ associated to $\Gamma$ has
generators $\{\sigma_s \  \  s\in S\}$.  There is a relation for
each label $m(s,t)< \infty$, namely  $\sigma_s
\sigma_t\sigma_s\sigma_t\cdots  = \sigma_t\sigma_s\sigma_t\sigma_s
\cdots,$ where there are $m(s,t)$ terms on each side,  $2 \leq
m(s,t) < \infty$.
 The \underline{Coxeter group} $C(\Gamma)$ associated to the Artin group $A(\Gamma)$  is obtained
by adding the relations $\sigma_s^2 = 1$ for every $s\in S$.  As
previously mentioned, we say that the group $A(\Gamma)$ has finite
type if its associated Coxeter group is finite. The braid group
${\bf B}_n$ is an example of an Artin group of finite type; its
associated Coxeter group being the symmetric group $\Sigma_n$.
Note that by definition the defining relations in the
presentations that we just described for Artin groups all preserve
word length.

One of the many interesting properties of Artin groups is that for
every Artin group there is an associated monoid $A^+$, and just as
every Garside monoid embeds in its group, it was proved by Paris
in \cite{Paris2002} that every Artin monoid injects in its group.
The injectivity property holds in the following strong sense too:
If two positive words $P, P'$ represent the same element of $A$,
then they represent the same element of $A^+$.  Paris's proof is
completely different from Garside's proof of injectivity in the
case of the braid group, which is the basis for the known proofs
of the same fact for Garside groups.  On the other hand, the other
properties that are needed to obtain a Garside structure may or
may not hold, for example it is definitely not true that every
element in an Artin group can be written in the form $NP$, where
$N$ is negative and $P$ is positive.

\begin{remark}{\rm
An interesting special case of Artin groups are the
\underline{right-angled} Artin groups. They are Artin groups in
which the defining relations are all of the form $\sigma_s\sigma_t
= \sigma_t\sigma_s$.  For example, the right-angled Artin group associated to
the braid group is defined by the presentation:
\begin{equation}
\label{equation:right angle braid group}
 \langle \sigma_1,\dots,\sigma_{n-1}\ | \  \sigma_i \sigma_k = \sigma_k \sigma_i \  \ {\rm if} \  \ |i-k| \geq 2  \     \rangle
\end{equation}}
\end{remark}

\begin{openproblem}
\label{openproblem:which Artin groups are linear?}
{\rm In $\S$\ref{subsection:the Lawrence-Krammer representation} we
explained the fairly recent proof that the braid groups are
linear.  This leads one to ask, immediately, whether the same is
true for its natural generalizations, e.g., Garside groups and Artin groups. It turns out that,
like the braid groups, all right-angled Artin groups have faithful
matrix representations \cite{Digne}, \cite{DJ2000}. Of course any
Artin group that injects as a subgroup of a related braid group is
also linear, e.g. Artin groups of type $B_n$.  It was proved recently by A. Cohen and D. Wales \cite{C-W} and simultaneously by F. Digne \cite{Digne} that all finite type Artin groups are linear. Going beyond that, the matter seems to be wide
open and interesting.} $\clubsuit$
\end{openproblem}

%% JOAN -- I changed "it's" to "it is" here.
\begin{remark}
{\rm  For reasons of space, we have not included any significant
discussion of the vast literature on Artin groups and associated
complexes on which they act. It is a pity to omit it,  because it
is a major area, and much of it had its origins in work on ${\bf
B}_n$.}
\end{remark}

\subsection{Braid groups and public key cryptography}\label{subsection:Braid
groups and public key cryptography} The problem which is the focus
of of `public key cryptography' was mentioned, very briefly, in
$\S$\ref{subsubsection:braid groups and public key cryptography}.
The basic issue is how to send information, securely, over an
insecure channel. The solution is always to use some sort of code
whose main features are known to the sender and recipient, but
which cannot be deduced by a viewer who lacks knowledge of the
shared keys. To the best of our knowledge, all solutions to this
problem rest on the same underlying idea: they make use of
problems which have a precise answer, which is known to both the
sender and recipient,  but one which is deemed to be so difficult
to compute that it is, in effect, unavailable to a viewer, even
though the viewer has all the necessary data to deduce it.  The
earliest such schemes were based upon the difficulty of
factorizing large integers into a product of primes. The
individuals who wish to exchange data over a public system are
denoted $A$ and $B$.  In a vastly oversimplified version, Both $A$
and $B$ have agreed, privately, on the choice of a prime number
$p$. The sender chooses another prime $q$ and transmits the
product  $n=pq$.  A viewer may learn $n$, but because of the
difficulty of factorizing $n$ into primes cannot deduce $p$ and
$q$. The recipient, who knows both $p$ and $n$, has no problem
computing $q$.

 A more recent approach is due to W. Diffie and M. Hellman \cite{Diffie-Hellman}. As before, both $A$ and $B$ have agreed, privately, on the choice of a prime number $p$ and a generator $g$ of the finite cyclic group $\ints/p\ints$.  A chooses a number $a$ at random and computes $g^a ({\rm mod} \ p)$, which she sends to $B$ on the public channel.  As for $B$ he chooses a number $b$ at random, computes $g^b (({\rm mod}  \ p)$ and sends in to $A$ on the public channel. Since $A$ now knows both $g^b$ and $a$ she can compute $(g^a)^b = g^{ab} ({\rm mod}) \ p$.  Similarly, since $B$ now knows both $g^a$ and $b$, he can compute
$(g^b)^a = g^{ba} ({\rm mod} \ p) = g^{ab} ({\rm mod} \ p)$.  So
both know $g^{ab} ({\rm mod} \ p)$. As for the viewer, he knows
both $g^a ({\rm mod} \ p)$ and $g^b ({\rm mod} \ p)$, however,
because of the known difficulty of computing discrete logs, the
crucial information $g^{ab} ({\rm mod} \ p)$ is in effect
unavailable to the viewer.

A rather different set of ideas was proposed in \cite{AAG1999} and
\cite{KLC}, and this is where the braid group comes in.  The
security of a system that is based upon braid groups relies upon
the assumption  that the word and conjugacy problems in the braid
group have both been solved, but the conjugacy problem is
computationally intractable whereas the word problem is not.
However, as we have seen, that matter seems to be wide open at
this moment.  We have described the underlying mathematics behind
the ElRifai-Morton solution to the word problem, and the best of
the current solutions to the word and conjugacy problem, namely
that of Gebhardt. The solution to the word problem is used in the
conjugacy problem. The reason that the conjugacy problem is not
polynomial is that we do not understand enough about the structure
of the summit set, the super summit set and the ultra summit set.
Our strong belief is that these matters will be settled.   The
assumption that our current lack of understanding of aspects of
the mathematics of braids means that they cannot be understood
seems unwarranted.

\subsection{The Nielsen-Thurston approach to the conjugacy problem in ${\bf B}_n$}
\label{subsection:the conjugacy problem:the dynamic approach} In
this subsection we consider the conjugacy problem in ${\bf B}_n$
from a new point of view.  We return again to the interpretation
given in Theorem 1 of $\S$\ref{subsection:bBn and bPn as mapping
class groups} of the braid group ${\bf B}_n$ as the mapping class
group ${\cal M}_{0,1,n}$ of the punctured disc $S_{0,1,n} =
D^2_n$.  In this context, then, we emphasize that the term
\underline{braid} refers to a mapping class, that is, an isotopy
class of diffeomorphisms.

The Nielsen-Thurston classification of mapping classes of a
surface is probably the single most important advance in this
theory in the last century, and we review it here. For simplicity,
we shall focus for now on the mapping class group of a closed
surface ${\cal M}_g = {\cal M}_{g,0,0}$ in the notation of
$\S$\ref{subsection:bBn and bPn as mapping class groups}.  We note
that the groups ${\cal M}_g$ have trivial center for all $g\geq 3$.

We make the tentative definition that an element of ${\cal M}_g$ is `reducible' if it
preserves up to isotopy a family of disjoint nontrivial curves on
the punctured disc. Such a family of curves is known as a
\underline{reduction system}. (Throughout this discussion, we use
`curve' to refer to the isotopy class of a curve.  By `disjoint
curves', we mean two distinct isotopy classes of curves with
respective representatives which are disjoint.)
It was proved in \cite{BLM1983} that for a reducible map $\phi$, there exists a \underline{essential
reduction system}, denoted ERS($\phi$). A curve $c \in
\rm{ERS}(\phi)$ if and only if
\begin{enumerate}
\item There exists an integer $k$ such that $\phi^k(c) = c$.
\item If a curve $x$ has nonzero geometric intersection with $c$,
then $\phi^m (x) \neq x$ for all integers $m$.  In particular, any curve $x$ which has nonzero geometric intersection with a curve $c\in$ ERS($\phi$) is not in ERS($\phi)$.
\end{enumerate}
It is proved in \cite{BLM1983} that the curves in ERS$(\phi)$ are contained in every reducing system for $\phi$, and are a minimal reduction system for $\phi$. Also, the system of curves ERS($\phi)$ is unique. Keeping all this in mind, we now define an element $\phi\in {\cal M}_g$ to be \underline{reducible} if it fixes the curves in an essential reduction system, i.e. in ERS$(\phi)$.

A mapping class $\phi \in {\cal M}_g$ is \underline{periodic} if
$\phi^k = 1$ in ${\cal M}_g$.

We note that a reducible mapping class always contains a
representative diffeomorphism which fixes a given reduction system
setwise.  Likewise, Nielsen showed that a periodic mapping class
always contains a representative which is periodic as a
diffeomorphism.
Finally, a mapping class which is neither reducible nor periodic
is \underline{pseudo-Anosov} (abbreviated as PA).

Braid groups, as the mapping class group of a punctured disc,
admit a similar classification.  The group ${\bf B}_n$ is
torsion-free, but by analogy with the above case of mapping class
groups of closed surfaces of genus at least two, which are
centerless, we say that an element $\omega \in{\bf B}_n$ is
\underline{periodic} if for some integer $k$, $\omega^k$ is
isotopic to a full Dehn twist on the boundary of the disc $D^2_n$.
In terms of generators and relations this is equivalent to
saying that the braid is a root of $\Delta^2$, where $\Delta$ is
the Garside braid.  (Recall from $\S$\ref{subsection:Garside} that
the square of $\Delta_n$ generates the center of ${\bf B}_n$.)

The definitions of reducible braids and pseudo-Anosov (PA) braids
require no alteration other than replacing the surface $S_g$ with
the punctured disc.

\begin{remark}{\rm
It is clear for ${\bf B}_n$ and ${\cal M}_g$ alike, we also have a
classification of elements as reducible, periodic and pseudo-Anosov up to conjugacy.}
\end{remark}

Each of these three possibilities reveals new structure, so we
consider them one at a time.  A major reference for us is the
paper \cite{BGN1995}, which gives an algorithm for finding a
system of reducing curves if they exist, and for recognizing
periodic braids. This paper did not receive much attention at the
time that it was written because it had the misfortune to be
written simultaneously and independently with the ground-breaking
and much more general papers of Bestvina and Handel \cite{BH1995}.
However, as is often the case, one learns very different things by
examining a particular case of a phenomenon in detail, and by
proving a broad generalization of the same phenomenon, and that is
what happened here.

\paragraph{Periodic braids.}  The first author learned from \cite{BGN1995}
that the classification of periodic braids had been solved via the
work of Kerekjarto (1919) and Eilenberg (1935), who proved that up
to conjugacy a periodic braid is a power of either $\delta$ or
$\alpha$, where $\delta = \sigma_1\sigma_2\dots\sigma_{n-1}$ and
$\alpha = \sigma_1\sigma_2\dots\sigma_{n-1}\sigma_1$. One
verifies, using the elementary braid relations in the classical
presentation, that $\delta^n = \alpha^{n-1} = \Delta^2$.  To
visualize the assertion for $\delta$ geometrically, place the
punctures at polar angles $2\pi k/n$ around a circle of radius
$r<1$, and think of $\delta$ as a rotation of $D^2_n$ of angle
$2\pi/n$. To visualize $\alpha$, do the same, only now place one
of the punctures at the origin and arrange the remaining $n-1$
punctures symmetrically at polar angles $2\pi k /(n-1)$.

\paragraph{Reducible braids.} Let $\omega = \Delta^I L_1 L_2 \cdots L_s$
be a reducible mapping class in left greedy normal form, as
described in $\S$\ref{subsection:Garside}.  Our initial model for
$D^2_n$ will be the unit disc, with the $n$ punctures arranged
along the real axis, and placed symmetrically so that they divide
the interval $[-1,1]\subset\reals$ into $n+1$ equal line segments.
It is shown in \cite{BGN1995} that one can choose a braid
$\omega'$ conjugate to $\omega$ and in the super summit set of
$\omega$ which fixes the simplest possible family of closed
curves, namely a family $C$ of geometric ellipses whose centers
are on a horizonal `axis' through the punctures, and which are
chosen in such a way that the axis bisects the discs that the
curves in $\rm{CRS}(\omega')$ bound.

\begin{openproblem}
\label{openproblem;reducible braids}
{\rm Is it true that, if a reducible braid is in its ultra summit set,  then its invariant multicurves are
can always be chosen to be geometric ellipses?} $\clubsuit$
\end{openproblem}

The remaining braiding may then be thought of as going on inside
tubes, also the tubes may braid with other tubes. With this model
it should be intuitively clear that the choice of such a
representative for the conjugacy class of a braid which permutes
the tubes in the required fashion may require very different
choices (again up to conjugacy) for the braids which are inside
the tubes, and that the sensible way to approach the problem is to
cut the initial disc open along the reducing curves and focus on
the periodic or PA maps inside the tubes. This is, of course, a
very different approach from the one that was considered in
$\S$\ref{subsection:Garside}   above, where no such considerations
entered the picture.

\paragraph{Pseudo-Anosov braids.} In the PA case (which is the generic
case), there is additional structure, and now our description
becomes very incomplete.  In this case \cite{Thurston1988} there
exist two projective measured foliations ${\cal F}^u$ and ${\cal
F}^s$, which are preserved by an appropriate representative $w$ of
$\{\omega\}$.   Moreover the action of $w$ on   ${\cal F}^u$ (the
unstable foliation) scales its measure by a real number
$\lambda>1$, whereas the action on ${\cal F}^s$ (the stable
foliation) scales its measure by $1/\lambda$.   These two
foliations and the scaling factor $\lambda$ are uniquely
determined by the conjugacy class of $\{\omega\}$, however the
triplet $({\cal F}^u, {\cal F}^s, \lambda )$ does not determine
$\{\omega\}$.  To explain the missing pieces, we replace the
invariant foliations by an invariant `train track', with weights
associated to the various branches.

In the case of pseudo-Anosov mapping classes acting on a
once-punctured surface, a method for enumerating the train tracks
is in Lee Mosher's unpublished PhD thesis \cite{Mosher-1983}.  But
Mosher's work was incomplete and remained unpublished for many
years, even as many of the ideas in it were developed and even
expanded, leading to some confusions in the literature about
exactly what is known and what remains open.   At this writing a
complete solution to the conjugacy problem for braid groups or
more generally surface mapping class groups, based upon the
Nielsen-Thurston machinery, does not exist in the literature.
Mosher has a partially completed monograph in preparation, {\it
Train track expansions of measured foliations}, which promises to
give such a solution, but for more general surfaces, i.e. for $S_{g,b,n}$
very little seems to be known.   We therefore pose it as an open problem:

\begin{openproblem}
\label{openproblem:conjugacy problem in the MCG} {\rm Investigate
the conjugacy problem in the mapping class groups $\cM_{g,b,n}$
with the goals of (i) pinning down precisely what is known for
various triplets $(g,b,n)$; (ii) describing all cases in which
there is a complete solution; and (iii) describing what remains to
be done in the simplest cases, that is $\cM_{g,0,0}$ and
$\cM_{g,0,1}$.} $\clubsuit$
\end{openproblem}

Essentially nothing is known, at this writing, about the interface
between the dynamic and combinatorial solutions to the conjugacy
problem in the braid groups, and still less in the general case of
the more general mapping class groups ${\cal M}_{g,b,n}$, to which
the entire Thurston machinery applies.  We know of nothing which
even hints at related dynamic structures in Artin groups or
Garside groups.
\begin{openproblem}
\label{openproblem:Garside approach to MCG} {\rm  Does the Garside
approach to the conjugacy problem in $\cM_{0,1,n}$ generalize to a
related approach to the problem in $\cM_{g,1,0}$ for any $g>0$?}
$\clubsuit$
\end{openproblem}

\begin{openproblem}
\label{openproblem:Nielsen-Thurston approach to Artin groups} {\rm
Does the Nielsen-Thurston approach to the conjugacy problem in
$\cM_{0,1,n}$ generalize to a related approach to the problem in
any other Artin group?} $\clubsuit$
\end{openproblem}

\subsection{Other solutions to the word problem}
\label{subsection:other solutions to the word
problem} We review, very briefly, other ways that the word and
conjugacy problems have been solved in the braid groups.  We
restrict ourselves to results which revealed aspects of the
structure of ${\bf B}_n$ that has had major implications for our
understanding of braid groups, even when the implications for the
word and/or conjugacy problems fall short of that criterion.

\paragraph{1. Artin's solution to the word problem.} The earliest
solution to the word problem was discovered in 1925 by Artin, in
his first paper on braids \cite{Artin1925}. It was based upon his
analysis of the structure of the pure braid group, which we
described in $\S$\ref{subsection:definition of bBn via generators
and relations}, and in particular in the defining relations
\ref{equation:presentation bPn} and the sequence given just below
it:
 $$\{1\} \to  {\bf F}_{n-1} \to {\bf P}_n \stackrel{\pi_n^\star}\to {\bf P}_{n-1} \to \{1\}. $$
 The resulting normal form, described in \cite{Artin1925} as \underline{combing} a braid, is well known.  We refer the reader to Artin's original paper for a very beautiful example.  In spite of much effort over the years nobody has managed, to this day, to use related techniques to solve the conjugacy problem in ${\bf B}_n$, except in very special cases.

\paragraph{2. The Lawrence-Krammer representation.}  See
$\S$\ref{section:representations of the braid groups} above.  When
a group admits a faithful matrix representation, there exists a
fast way to solve the word problem.  It is interesting to note
that Krammer's first proof of linearity, in the case of ${\bf
B}_4$, used the solution to the word problem which came from the
presentation (\ref{equation:new presentation}) in
$\S$\ref{section:introduction} to this paper.  See
\cite{Krammer1999}.  Also, when he proved linearity in the general
case, in \cite{Krammer2002},  he gives two proofs. The second
proof, Theorem 6.1 of \cite{Krammer2002}, shows that the
Lawrence-Krammer matrices detect the infimum and supremum of a
braid.  It follows as an immediate corollary (noted in
\cite{Krammer2002})
 that the representation is faithful, because any braid which has non-zero infimum and supremum cannot be the identity braid.

As regards the conjugacy problem, there is a difficulty.  The
image of ${\bf B}_n$ under the isomorphism given by the
Lawrence-Krammer representation yields a group $\cB_n$ which is a
subgroup of the general linear group $GL_m(\ints[q^{\pm 1}, t^{\pm
1}] )$, where $m=(n)(n-1)/2$. To the best of our knowledge it is
unknown at this time how to describe $\cB_n$ in any way that
allows one to test membership in $\cB_n$.  Lacking such a test,
any class invariants which one finds in this way will be limited
in usefulness, because if $\omega\in {\bf B}_n$,  there will be no
way to distinguish between class invariants which arise from
conjugation by elements in  $\cB_n$  from those which arise from
conjugation by more general  elements of $GL_m(\ints[q^{\pm 1},
t^{\pm 1}] )$ from invariants  thus one cannot hope for a complete
solution.
\begin{openproblem} {\rm Investigate the
eigenvalues and the trace of the Lawrence-Krammer matrices.
}$\clubsuit$
\end{openproblem}

\paragraph{3. The Dehornoy ordering.}  A group or a monoid has a
\underline{right-invariant} (resp.
\underline{left}-\underline{invariant}) ordering if there exists a
strict linear ordering of its elements, denoted $<$, with the
property that if $f,g,h\in G$, then $f < g$ implies  $fh < gh$
(resp. $hf < hg$).    To the best of our knowledge nobody had
considered the question of whether ${\bf B}_n$ had this property
before 1982, when Patrick Dehornoy announced his discovery that
the groups ${\bf  B}_n$ admit such an ordering, and that it can be
chosen to be either right-invariant or left-invariant, but not
both.   In the 5-author paper \cite{FGRRW} the Dehornoy ordering
was shown to have the following topological meaning. We restrict
to the the right-invariant case, the two results being essentially
identical. We regard ${\bf B}_n$ as ${\cal M}_{0,1,n}$, where the
surface $S_{0,1,n}$ is the unit disc in the complex plane and the
punctures lie on the real axis in the interval $(-1,1)$.  The
punctures divide $[-1,1]\subset\reals$ into $n+1$ segments which
we label $\xi_1,\xi_2,\dots,\xi_{n+1}$ in order, where $\xi_1$
joins $\{-1\}$ to the first puncture.  Choose
$[\omega_1],[\omega_2]\in {\bf B}_n$, and representatives
$\omega_1,\omega_2$ with the property that $\omega_1([-1,1])$ and
$\omega([-1,1])$ intersect minimally.  Note that
$\omega_k([-1,1])$ divides $S_{0,1,n}$ into two halves, and it
makes sense to talk about the upper and lower half of
$\omega_2([-1,1])$ because $\{-1\}$ and $\{1\}$ are on $\partial
S_{0,1,n}$ and $\omega_k \  | \ \partial S_{0,1,n} $  is the
identity map. Then
 $[\omega_2] >[\omega_1]$ if $\omega_1(\xi_i) = \omega_2(\xi_i)$ for $i=1,\dots,j-1$ and an initial segment of $\omega_2(\xi_j)$ lies in the upper component of $S_{0,1,n} - \omega_1([-1,1])$.

It was proved in \cite{Wiest1999} that the resulting
left-invariant ordering of ${\bf B}_n$ extends the ElRifai-Morton
left partial ordering that we described in
$\S$\ref{subsection:generalization:from bBn to Garside groups}
above.   It is proved in \cite{DDRW} that there are infinitely
many other left orderings.

The subject blossomed after the paper \cite{FGRRW} appeared. In
fact so much work resulted that  there is now a 4-author monograph
on the subject \cite{DDRW}, written by some of those who were the
main contributors, containing an excellent review of what has been
learned during the 10 years since the initial discovery.  We do
not wish to repeat what is readily available elsewhere, especially
because we are not experts, so we point the reader to Chapter 3 of
\cite{DDRW}, where it is shown that the Dehornoy  ordering leads
to a solution to the word problem. Unfortunately, the solution
so-obtained is exponential in $|\omega|$, whereas the solution
described in (G1)-(G6) above is quadratic in $|\omega|$.

\paragraph{4.  B$_{\rm n}$ as a subgroup of Aut (F$_{\rm n}$).} It is
well known that ${\bf B}_n$ has a faithful representation as a
subgroup of Aut(F$_{\rm n}$), for example see \cite{Birman1974}
for a proof.  This of course gives a solution to the word problem,
but (like the Dehornoy solution) it is exponential in $|\omega|$.
On the other hand, it is extremely interesting that the
Nielsen-Thurston machinery, described in $\S$\ref{subsection:the
conjugacy problem:the dynamic approach} above, generalizes to the
entire group Aut(F$_{\rm n})$, giving yet another instance where
the braid group is at the intersection of two rather different
parts of mathematics, and could be said to have pointed the way to
structure in the second based upon known structure in the first.

In this section we have described several solutions to the word
problem, and noted that one of them (the modified Garside approach
of (G1)-(G6)) is ${\cal O}(|\omega|^2(n)(log (n))$, where $|\omega|$
is the letter length of two arbitrary representatives of elements
of ${\bf B}_n$, using the classical presentation
(\ref{equation:classical presentation}) for ${\bf B}_n$.  The word
problem seems to be one short step away from the
\underline{non-minimal braid} problem:  given a word  $\omega$ in
the generators $\sigma_1,\dots,\sigma_{n-1}$ and their inverses,
determine whether there is a shorter word $\omega^\prime$ in the
same generators which represents the same element of ${\bf B}_n$?
For, it is clear that a decision process exists: list all the
words that are shorter than the given one, thin the list by
eliminating as many candidates as possible by simple criteria such
as preserving length mod 2, an obvious invariant, and then test,
one by one, whether the survivors represent the same element of
${\bf B}_n$ as the given word $\omega$?   It therefore seemed
totally surprising to us that in 1991 M. S. Patterson and A. A.
Razborov proved:
\begin{theorem}{\rm \cite {P-R} :}
\label{theorem:non-minimal braids} The non-minimal braid problem
is NP complete.
\end{theorem}
Thus, if one could find an algorithm to decide whether a given
word $\omega$ is non-minimal, and if the algorithm was polynomial
in $|\omega|$, one would have proved that $P = NP$!  To the best
of our knowledge, essentially nothing has been done on this
problem. In this regard we suggest two research problems:

\begin{openproblem} {\rm
The proof that is given in \cite{P-R} is very specific to the
classical presentation (\ref{equation:classical presentation}) for
${\bf B}_n$.    Can it be adapted to the new presentation? To
other presentations?}
$\clubsuit$
\end{openproblem}

\begin{openproblem} {\rm
Investigate the shortest word problem in the braid group ${\bf B}_n$, using the classical presentation. } $\clubsuit$
 \end{openproblem}

In $\S$\ref{subsection:mtws} we discussed and proved the Markov
Theorem Without Stabilization in the special case of the unknot.
See Theorem \ref{theorem:mtws}, which asserts that there is a
complexity measure on closed braid representatives of the unknot,
and using it a sequence of strictly complexity-reducing
destabilizations and exchange moves that reduce a closed braid
diagram for the unknot to a round planar circle. In \cite{BH1998},
theorem \ref{theorem:mtws} was used to develop an algorithm for
unknot recognition.  In \cite{BRBV2002} the first steps were taken
to develop a computer program to realize the resulting algorithm.
The algorithm is far from being practical in even simple cases,
however we are now in a position to describe three open problems,
all closely related to (G1)-(G-11) above, which would lead to an
effective solution to the unknot recognition problem:

\begin{openproblem} {\rm Develop an algorithm that will detect when the conjugacy
class of a closed braid admits a destabilization. }
$\clubsuit$
\end{openproblem}

An $n$- braid
$W$ admits an exchange move if it is conjugate to a braid of the
form $U \sigma_{n-1} V \sigma_{n-1}^{-1}$, where $U$ and $V$
depend only on $\sigma_1,\dots,\sigma_{n-2}$.  Thus, up to
conjugation, $W$ is a product of two reducible braids, one
positively reducible and the other negatively reducible.   By Theorem \ref{theorem:mtws}, it
may be necessary to modify the given conjugacy class by exchange
moves in order to jump from the given class to a new class (if it
exists) which admits a destabilization. Thus there is an unknown
complexity measure on conjugacy classes, with the ones which
admit destabilizations being especially nice.  This leads us to our
second open problem.

\begin{openproblem} {\rm Develop an algorithm
that will detect when the conjugacy class of a closed braid admits
an exchange move, and when a collection of classes are
exchange-equivalent.}
$\clubsuit$
\end{openproblem}

This would still not be a complete tool, because exchange moves
can be either be complexity-reducing or complexity-increasing, and
unless we can tell the difference (that is the content of open
problem 3 below) we are left with the option of trying a sequence
of exchange moves of unpredictable length in our search for
destabilizations.

\begin{openproblem} {\rm The complexity measure that was introduced in
$\S$\ref{subsection:mtws} is `hidden' in the braid foliation of an
incompressible surface whose boundary is the given knot. This is a
highly non-trivial matter, because if we had the foliated surface
in hand, then we would be able to compute its Euler characteristic
and  would know whether we had the unknot. On the other hand, the
braid foliation determines the embedding of its boundary (see
Theorem 4.1 of \cite{BirFink}), so there is no essential obstacle
to `seeing' the complexity measure in the given braid. That is the
essence of our third open problem, which asks that we recognize
how to translate the complexity measure that was used in the proof
of Theorem \ref{theorem:mtws} into a complexity measure on closed
braids which will be able to distinguish exchange moves that
reduce complexity from those which do not.}
$\clubsuit$
\end{openproblem}

\newpage

\section{A potpourri of miscellaneous results}
\label{section:a potpourri of miscellaneous results}
This section is for leftovers--- topics on which there have been interesting new discoveries which did not seem to fit well anywhere else.

\subsection {Centralizers of braids and roots  of braids}
\label{subsection:centralizers and roots  of braids} The
\underline{mixed braid groups} are defined in \cite{GM-W2003} to
be the braids which preserve a given partition of the puncture
points on the disc.  In \cite{GM-W2003} Gonzales-Meneses and Wiest
gave full descriptions of the centralizer of a braid in terms of
semi-direct and direct products of mixed braid groups, and found
sharp bounds on the number of generators of the centralizer of a
braid.

In a related paper \cite{GM2003} Gonzales-Meneses proved that, up to conjugacy, braids have unique roots. That is, if $\omega\in {\bf B}_n$ and if $\alpha,\beta\in {\bf B_n}$ have the property $\alpha^k = \beta^k = \omega$, then $\{\alpha\} = \{\beta\}$.

\subsection{Singular braids,  the singular braid monoid, and the desingularization map}
\label{subsection:singular braids}
The singular  braid monoid ${\bf SB}_n$ is a monoid extension of the braid group ${\bf B}_n$. It was introduced in \cite{Baez1992} and, simultaneously and independently, in \cite{Birman1993}. Its definition was suggested by the mathematics which revolved about  Vassiliev invariants of knots and links.   To define it, we need to describe a presentation for ${\bf SB}_n$, taken from  \cite{Birman1993}.    There are three types of generators, which we call $\sigma_i, \sigma_i^{-1}$ and $  \tau_i, \ \  1\leq i \leq n-1$.  Here the $\sigma_i^{\pm 1}$ are to be thought of as the classical positive and negative elementary braids and the $\tau_i$ are to be thought of as  \underline{elementary singular braids}. The braid $\tau_i$ is obtained from $\sigma_i$ (see Figure \ref{figure:Delta} (ii)) by identifying strands $i$ and $i+1$ as they cross.  Defining relations in ${\bf SB}_n$ are:
$$\sigma_i\sigma_i^{-1} = \sigma_i^{-1}\sigma_i = 1, \ \ \sigma_i\tau_i = \tau_i\sigma_i, $$
$$\sigma_i \sigma_j = \sigma_j\sigma_i , \ \ \sigma_i \tau_j = \tau_j\sigma_i, \ \ \tau_i\tau_j = \tau_j\tau_i  \  \ \ \ \ { \rm   if} \ \ \ \  |i-j| \geq 2 $$
$$\sigma_i\sigma_j\sigma_i = \sigma_j\sigma_i\sigma_j, \ \ \ \ \sigma_i\sigma_j\tau_i = \tau_i\sigma_i\sigma_j \ \ \ \  {\rm if} \ \ \ \  |i-j| = 1. $$

The \underline{desingularization map} is a homomorphism from ${\bf SB}_n$ to the group ring $\ints {\bf B}_n$ of the braid group, defined by $\psi(\sigma_i^{\pm 1}) = \sigma_i^{\pm 1}, \ \ \ \ \ \ \psi(\tau_i) = \sigma_i - \sigma_i^{-1}$.
Birman used this map to  develop the relationship between
Vassiliev invariants and quantum groups'. It was conjectured in
\cite{Birman1993} that the map $\psi$ is an embedding.   After
several proofs of special cases of the conjecture, it was settled
by Luis Paris in the affirmative in 2003 \cite{Paris2003}.

During the 10 year interval after the introduction of the singular
braid monoid, and before the proof of the embedding theorem, it
came as quite a surprise when it was discovered that there was a
new group on the scene--the \underline{singular braid group} of
\cite{FKR1998}.  To this day we are unsure of its significance,
although its existence is unquestioned!  It is most easily defined
via generators and relations, starting with the presentation that
we  just gave for the singular braid monoid, and then adding one
new generator $\bar{\tau}$ (to suggest that it behaves the way
that the inverse of $\tau$ ought to behave). Defining relations
for ${\bf GB}_n$ are all the relations in ${\bf SB}_n$, plus ones
satisfied by the new generator. The latter are  `monoid relations'
which are the same as those  in the singular braid monoid, but
substituting $\bar{\tau}$ for $\tau$,  and two additional
relations, namely $\tau\bar{\tau} = \bar{\tau}\tau = 1.$
Pictorially, one has two types of singular crossings, and they
annihilate one another.  The main result in \cite{FKR1998} is that
${\bf SB}_n$ embeds in ${\bf GB}_n$.

\subsection{The Tits conjecture}
\label{subsection:the Tits conjecture} The Tits conjecture is very
easy to state, in its simplest form. Consider the classical
presentation of the braid group ${\bf B}_n$, i.e. the presentation
(\ref{equation:classical presentation}).  The Tits conjecture, in
the special case of the braid group, is that the subgroup $G$ of
${\bf B}_n$ generated by the elements $T_i = \sigma_i^2$ has the
presentation $\langle T_1,\dots,T_{n-1} \ | \ T_i T_j = T_j T_i \
\ {\rm if} \ \ |i-j|\geq 2 \rangle$. A generalized version of the
conjecture (the generalization relates to the arbitrary choices of
the powers) was proved by Crisp and Paris in 2001, settling a
question which had plagued the experts for many years:

\begin{theorem}  {\rm \cite{C-P2001}}
\label{theorem:Tits conjecture} Let $S$ be a finite set, let
$\Gamma$ be a Coxeter graph over $S$, and let $A(\Gamma)$ be the
associated Artin group, as defined in $\S$\ref{subsection:Artin
groups}, with generating set $\{ \sigma_s,\ \ s\in S \}$.
Associate further to each $s \in S$ an integer $m_s \geq 2$.
Consider the subgroup $G$ of $A(\Gamma)$ generated by the elements
$\{ T_{s,m_s} = \sigma_s^{m_s} \}$.    Then defining relations
among the generators of the subgroup $G$ are `the obvious ones',
namely that $T_{s,m_s}$ and $T_{t,m_t}$ commute in $G$ if and only
if they commute in $A(\Gamma)$.  No other relations are needed.
\end{theorem}

The proof is interesting, because it introduces a technique which is very closely related to the themes
that we have explored in this article.   The basic idea is that, for a key subclass of Artin groups there is a representation $f$ of
$A(\Gamma)$ into the mapping class group ${\cal M}(S)$ of a connected surface $S$ with boundary that is associated to the graph $\Gamma$,
which induces an action of $A(\Gamma)$ on a monoid determined by $\Gamma$.  Since the group $H(\Gamma)$ which is presented in the theorem
has an obvious homomorphism onto $A(\Gamma)$, and the proof shows that the restriction of the action to $H(\Gamma)$ gives the desired
isomorphism.

\subsection{Braid groups are torsion-free: a new proof}
\label{subsection:braid groups are torsion-free}
As we stated in $\S$ \ref{section:introduction}, it was
necessary to make choices in the writing of this review, and our decision was to be guided by the principle of focusing on new results or
new proofs of known results.  During the years since \cite{Birman1974} was written, many people have written to the first author with
questions about braids, and a regular question has been ``Isn't there a {\it simple} proof that the braid groups have no elements of
finite order?"   It is therefore fitting that we end this review with exactly that -- a beautiful simple proof, based upon the discovery,
due to Dehornoy, that the braid groups admit a left-invariant ordering:

\begin{theorem}
\label{theorem:braid groups are torsion free}
The groups ${\bf B}_n,  \ n=1,2,3,\dots$ are all torsion free.
\end{theorem}
\pf  See $\S$\ref{subsection:other solutions to the word problem}
for a discussion of the Dehornoy left-invariant ordering of the
braid groups. Choose any element $g \in {\bf B}_n, \ \  g \not=
1.$ Replacing $g$ if necessary with $g^{-1}$ we may assume that $
1 < g$. Since the ordering is left-invariant we then have that $g<
g^2$ and $g^{-1} < 1$. Iterating, $\cdots < g^{-3} < g^{-2} <
g^{-1}< 1 < g < g^2 < g^3 < \cdots $. \endpf

\newpage

\centerline{{\bf \LARGE{Appendix: Computer programs}} }

\

 In 2004, it is almost as important to know about computer tools as it is to have a guide to the literature, so we
supplement our bibliography with a guide to the computer tools that we know about and which have been useful to us and colleagues.

\paragraph{Changing knots and links to closed braids.} Vogel's proof of his method for changing arbitrary knot diagrams to
closed braid diagrams is ideal for computer programming.  We refer the reader to the URL
http://www.layer8.co.uk/maths/braids/,
for a program, due to Andrew Bartholomew  and  Roger Fenn, which does this and much more.

\paragraph{Garside's algorithm for the word and conjugacy problems.}   Many people have programmed Garside's algorithm
for the word and conjugacy problem, but the one we know best is the program of Juan Gonzalez-Meneses, which can be downloaded from
http:www.personal.us.es/meneses.  The very robust version that the reader will find there computes Garside's normal forms for
braids, and
the ultra-summit set of a braid, a complete invariant of conjugacy.

\paragraph{Nielsen-Thurston classification of mapping classes in $\cM_{0,n+1}$.}  We know of two very useful
computer programs, all based upon the Bestvina-Handel algorithm \cite{BH1995}. Both assume that the boundary of the $n$-times punctured
disc has been capped with a disc, so that they compute in the mapping class group of ${\bf B}_n/ $ modulo its center, rather than in
${\bf B}_n$.  Equivalently, they work with the mapping class group of an $(n+1)$-times punctured sphere, where admissible maps fix the
distinguished point.

The first, due to W. Menasco and J. Ringland, can be downloaded from \\
http://orange.math.buffalo.edu/software.html, by following the link to  ``BH2.1 An Implementation of the Bestvina-Handel Algorithm".
The second, due to T. Hall,  can be downloaded from
http://www.liv.ac.uk/Maths/pure and following the links to the research group on Dynamical Systems, and then to the home page
of Toby Hall.   Both determine whether an input map is pseudo-Anosov, reducible or finite order and both find an invariant train track
and a train track map (which makes it possible to calculate essential dynamics in the class by Markov partition
techniques).  In the reducible case Hall's program provides a set of reducing curves.
Recently Hall updated his program to adapt it to MacIntosh OS-X computers, versions 10.2 and above.

\paragraph{Other software.}  We mention M. Thistlethwaite's Knotscape, because the program has a reputation for being very versatile and
user-friendly.  It is available for download from \\
 http://www.math.utk.edu/~morwen. It accepts as input knots that are defined as closed braids (and also knots defined by the
Dowker code or mouse-drawn diagrams), and locates it in the tables if it has at most 16 crossings.
The program computes numerous invariants, including the Alexander,
Jones, Homfly and Kauffman polynomials, and hyperbolic invariants (assuming
that the knot is hyperbolic).  The hyperbolic routines in Knotscape were
taken with permission from Jeff Weeks's program SnapPea; the procedures for
calculating polynomials were supplied by Bruce Ewing and Ken Millett, and
the procedure for producing a knot picture from Dowker code is part of
Ken Stephenson's Circlepack program.

\end{document}